# Mordukhovich Derivatives and Covering Constant for Set-Valued Metric Projection in General Banach Spaces

Jinlu Li

Department of Mathematics
Shawnee State University
Portsmouth, Ohio 45662 USA
jli@shawnee.edu

**Abstract**. In this paper, we find the explicit representation of the set-valued metric projection operator from $l_1$ to the unit closed ball in $l_1$. By using this representation, we investigate the properties of the Mordukhovich derivatives of the set-valued metric projection operator, which is applied to calculate its covering constant. As a special case, we consider the 2-d Banach space. We find the explicit solutions of the set-valued metric projection operator from the considered 2-d Banach space to its unit closed ball. By these solutions, we will calculate its Mordukhovich derivatives in details.

## 1. Introduction

In the theory of differentiation for single-valued mappings in Banach spaces, the most popular concepts of derivatives may be Gâteaux directional derivative and Fréchet derivative, which are generalization of ordinary derivatives in calculus. These concepts cannot be directly extended to set-valued mappings. However, in set-valued analysis, the theory of Mordukhovich derivatives (coderivatives) of set-valued mappings has played the crucial and fundamental roles in the theories of set-valued and variational analysis in Banach spaces (see [30−33]). The theory of Mordukhovich derivatives has been rapidly developed and has been widely applied to both pure and applied mathematics such as optimization theory, variational inequalities, game theory, economics theory and so forth (see [3−8, 21, 29−33]).

One of the most important applications of Mordukhovich derivatives for set-valued mappings is the well-known Arutyunov Mordukhovich Zhukovskiy Parameterized Coincidence Point Theorem (that is simply called the AMZ Theorem). The AMZ Theorem deals with the parameterized coincidence points of two set-valued mappings in Asplund spaces, in which, the covering constant for one of the considered two mappings is calculated by the Mordukhovich derivatives. For the convenience of readers, we review this theorem below.

Let $(X, \|\cdot\|_X)$ and $(Y, \|\cdot\|_Y)$ be real Banach spaces with topological dual spaces $X^*$ and $Y^*$, and with origins $\theta_X$ and $\theta_Y$, respectively. Let $(S, \tau)$ be a topological space. Let $F(\cdot): X \rightrightarrows Y$ and $G(\cdot, \cdot): X \times S \rightrightarrows Y$ be set-valued mappings. The Arutyunov Mordukhovich Zhukovskiy Parameterized Coincidence Point Theorem (Theorem 3.1 in [3]) is stated as follows.

**(**Arutyunov Mordukhovich Zhukovskiy Parameterized Coincidence Point Theorem**)** *Let the Banach spaces X and Y be Asplund and let P be a topological space. Let* $F$: $X \rightrightarrows Y$ *and* $G(\cdot, \cdot)$: $X \times P \rightrightarrows Y$ *be set-valued mappings. Let* $\bar{x} \in X$ *and* $\bar{y} \in Y$ *with* $\bar{y} \in F(\bar{x})$. *Suppose that the following conditions are satisfied:*

(A1) *The multifunction* $F: X \rightrightarrows Y$ *is closed around* $(\bar{x}, \bar{y})$.

(A2) *There are neighborhoods* $U \subset X$ *of* $\bar{x}$, $V \subset Y$ *of* $\bar{y}$, *and* $O$ *of* $\bar{p} \in P$ *as well as a number* $\beta \geq 0$ *such that the multifunction* $G(\cdot, p): X \rightrightarrows Y$ *is Lipschitz-like on* $U$ *relative to* $V$ *for each* $p \in O$ *with the uniform modulus* $\beta$, *while the multifunction* $p \rightarrow G(\bar{x}, p)$ *is lower/inner semicontinuous at* $\bar{p}$.

(A3) *The Lipschitzian modulus* $\beta$ *of* $G(\cdot, p)$ *is chosen as* $\beta < \hat{\alpha}(F, \bar{x}, \bar{y})$, *where* $\hat{\alpha}(F, \bar{x}, \bar{y})$ *is the covering constant of* $F$ *around* $(\bar{x}, \bar{y})$.

*Then for each* $\alpha > 0$ *with* $\beta < \alpha < \hat{\alpha}(F, \bar{x}, \bar{y})$, *there exist a neighborhood* $W \subset P$ *of* $\bar{p}$ *and a single-valued mapping* $\sigma: W \rightarrow X$ *such that whenever* $p \in W$ *we have*

$$F(\sigma(p)) \cap G(\sigma(p), p) \neq \emptyset \quad \text{and} \quad \|\sigma(p) - \bar{x}\|_X \leq \frac{\text{dist}(\bar{y}, G(\bar{x}, p))}{\alpha - \beta}. \tag{1.1}$$

From Condition (A3) in AMZ Theorem, in order to apply this theorem, we have to calculate the covering constant $\hat{\alpha}(F, \bar{x}, \bar{y})$ for the considered set-valued mapping $F$ at a point $(\bar{x}, \bar{y})$ (which will be reviewed in section 2). Meanwhile, the calculation of $\hat{\alpha}(F, \bar{x}, \bar{y})$ is based on the Mordukhovich derivatives $\widehat{D}^* F(x, y)$ for the point $(x, y)$ being around point $(\bar{x}, \bar{y})$. This leads us to consider the following question regarding to a considered set-valued mapping in a given problem in the analysis of Banach spaces.

Question A: How do we find the explicit Mordukhovich derivatives of set-valued mappings?

To answer Question A, in this paper, we select the metric projection operator in Banach spaces. The reason to choose this operator is that the metric projection operator has been widely applied to many fields such as, approximation theory, fixed point theory, optimization theory, variational analysis, and so forth (see [1, 2, 12, 35, 37].

If the considered Banach spaces are uniformly convex and uniformly smooth, then the metric projection operator (projecting to nonempty closed and convex subsets) is a single-valued mapping, which can be considered as a special case of set-valued mappings. In this case, the Gâteaux, Fréchet and Mordukhovich derivatives of the single-valued metric projection operator have been studied (see [9−11, 13−20, 22−28, 34, 37]. More precisely, when the metric project operator projects to closed balls, closed convex cones, closed and convex cylinders, the explicit solutions of Gâteaux, Fréchet and Mordukhovich derivatives of the single-valued metric projection operator have been proved in [14−16, 18−20, 23−28].

If the considered Banach spaces are not uniformly convex and uniformly smooth, the metric projection operator is a set-valued operator, in general. In this case, it is difficult and complicated to find the explicit solutions of the Mordukhovich derivative of the single-valued metric projection operator.

In this paper, we concentrate to consider the general real Banach space $l_1$. This paper is organized as follows:

In section 2, we review the concepts the Gâteaux and Fréchet derivatives of the single-valued mappings and the concept of Mordukhovich derivatives of the set-valued mappings.

In section 3, we first find the explicit representation of the set-valued metric projection operator from $l_1$ to the unit closed ball in $l_1$. Then, by this representation, we find the explicit solutions of the Mordukhovich derivative of the set-valued metric projection operator with respect to some special cases in $l_1$, by which, we calculate its covering constant.

In section 4, we concentrate to a 2-dimentional Banach space. We first use the results in section 3 to find the representation of the set-valued metric projection operator in 2-d Banach space, which is also projects to the closed unit ball. A directly proof is put in the Appendix. By this representation, will find its explicit solutions of the Mordukhovich derivative. From these solutions, we will find that the problem to calculate the Mordukhovich derivative of the set-valued is great complexity, even though in 2-d Banach space.

Section 5 is put in Appendix, in which we quickly consider the set-valued metric projection operator in a 3-dimentional Banach space. We find that the calculation of Mordukhovich derivative of the set-valued metric projection operator in 3-dimentional Banach spaces are much more complicated than the calculation in 2-d Banach spaces.

## 2. Review of Mordukhovich derivatives of set-valued mappings in Banach Spaces

For given Banach spaces $(X, \|\cdot\|_X)$ and $(Y, \|\cdot\|_Y)$, as usual, let $X^*$ and $Y^*$ denote the topological dual spaces of $X$ and $Y$ equipped with norms $\|\cdot\|_{X^*}$ and $\|\cdot\|_{Y^*}$, respectively. Let $\langle\cdot, \cdot\rangle_X$ denote the paring between $X^*$ and $X$ and let $\langle\cdot, \cdot\rangle_Y$ denote the paring between $Y^*$ and $Y$. For given $x_0 \in X$ and $r > 0$, let $B_X(x_0, r)$ denote the closed ball in $X$ with radius $r$ and centered at $x_0$. For given $y_0 \in Y$ and $r > 0$, we similarly write the closed ball $B_Y(y_0, r)$ in $Y$ with radius $r$ and centered at $y_0$.

Let $D$ be a nonempty subset of $X$. Let $f: D \to Y$ be a single-valued mapping. Let $\bar{x} \in D$ and $v \in X$ with $v \neq \theta_X$. If there is a vector in $Y$, which is denoted by $f'(\bar{x}, v)$ such that

$$f'(\bar{x}, v) = \lim_{t \to 0,\ \bar{x}+tv \in D} \frac{f(\bar{x}+tv)-f(\bar{x})}{t},$$

then, $f$ is said to be Gâteaux directionally differentiable at point $\bar{x}$ along direction $v$. The point $f'(\bar{x}, v) \in Y$ is called the Gâteaux directional derivative of $f$ at point $\bar{x}$ along direction $v$. If there is a continuous and linear mapping $\nabla f(\bar{x}): X \to Y$ such that,

$$\lim_{u \to \theta_X,\ \bar{x}+u \in D} \frac{f(\bar{x}+u)-f(\bar{x})-\nabla f(\bar{x})(u)}{\|u\|_X} = \theta_Y,$$

then $f$ is said to be Fréchet differentiable at $\bar{x}$ and $\nabla f(\bar{x})$ is called the Fréchet derivative of $f$ at $\bar{x}$.

The Mordukhovich derivative for set-valued mappings in Banach spaces forms the foundation of generalized differentiation in set-valued and variational analysis in Banach spaces (See [30−33]). We review the concepts and some properties of Mordukhovich derivative for set-valued mappings. For more details, the readers are referred to [31].

Let $F: D \rightrightarrows Y$ be a set-valued mapping. The graph of $F$ is defined by the following subset in $D \times Y$

$$\text{gph}F = \{(x, y) \in D \times Y: y \in F(x)\}.$$

For $(x, y) \in \text{gph}F$, that is, for $x \in D$ and $y \in F(x)$, the Mordukhovich derivative (which is also called, Fréchet coderivative, Mordukhovich coderivative, or just coderivative) of $F$ at point $(x, y)$ is a set valued mapping $\widehat{D}^*F(x, y): Y^* \rightrightarrows X^*$. For any $y^* \in Y^*$, it is defined by (see Definitions 1.13 and 1.32 in Chapter 1 in [31])

$$\widehat{D}^*F(x, y)(y^*) = \left\{ x^* \in X^*: \limsup_{\substack{(u,v)\to(x,y) \\ u\in D \text{ and } v\in F(u)}} \frac{\langle x^*, u-x\rangle_X - \langle y^*, v-y\rangle_Y}{\|u-x\|_X + \|v-y\|_Y} \leq 0 \right\}$$

$$= \left\{ x^* \in X^*: \limsup_{\substack{(u,v)\to(x,y) \\ (u,v)\in \mathrm{gph}F}} \frac{\langle x^*, u-x\rangle_X - \langle y^*, v-y\rangle_Y}{\|u-x\|_X + \|v-y\|_Y} \le 0 \right\}. \quad (2.1)$$

This can be rewritten as

$$\widehat{D}^*F(x,y)(y^*) = \left\{ x^* \in X^*: \liminf_{\substack{(u,v)\to(x,y) \\ (u,v)\in \mathrm{gph}F}} \frac{\langle y^*, v-y\rangle_Y - \langle x^*, u-x\rangle_X}{\|u-x\|_X + \|v-y\|_Y} \ge 0 \right\}. \quad (2.2)$$

If $(x,y) \notin \mathrm{gph}F$, then, we define

$$\widehat{D}^*F(x,y)(y^*) = \emptyset, \text{ for any } y^* \in Y^*.$$

By the above definition (1.1), $\widehat{D}^*F(x,y): Y^* \rightrightarrows X^*$ is a set valued mapping, which is called the Mordukhovich derivative (or the Mordukhovich coderivative) of $F$ at $(x,y)$.

In particular, let $f: D \to Y$ be a single valued mapping. Let $x \in X$ and $y \in Y$ with $y = f(x)$. The Mordukhovich derivative (coderivative) of $f$ at point $(x, y)$ is a set-valued mapping, denoted by $\widehat{D}^*f(x,y): Y^* \rightrightarrows X^*$, which is defined by, for any $y^* \in Y^*$,

$$\widehat{D}^*f(x,y)(y^*) = \left\{ z^* \in X^*: \limsup_{\substack{(u,f(u))\to(x,y) \\ u\in D}} \frac{\langle z^*, u-x\rangle_X - \langle y^*, f(u)-y\rangle_Y}{\|u-x\|_X + \|f(u)-y\|_Y} \le 0 \right\}.$$

In particular, if $f: D \to Y$ is a continuous single-valued mapping, then, for any $y^* \in Y^*$, we have

$$\widehat{D}^*f(x,y)(y^*) = \left\{ z^* \in X^*: \limsup_{\substack{u\to x \\ u\in D}} \frac{\langle z^*, u-x\rangle_X - \langle y^*, f(u)-y\rangle_Y}{\|u-x\|_X + \|f(u)-y\|_Y} \le 0 \right\}.$$

The following theorem shows the connection between Fréchet derivatives and Mordukhovich derivatives for single-valued mappings, which provides a powerful tool to calculate the Mordukhovich derivatives by the Fréchet derivatives of single-valued mappings between Banach spaces.

**Theorem 1.38 in [31]**. *Let $X$ be a Banach space with dual space $X^*$ and let $f: X \to Y$ be a single-valued mapping. Suppose that $f$ is Fréchet differentiable at $x \in X$ with $y = f(x)$. Then, the Mordukhovich derivative of $f$ at $x$ satisfies the following equation*

$$\widehat{D}^*f(x,y)(y^*) = \{(\nabla f(x))^*(y^*)\}, \text{ for all } y^* \in Y^*. \quad (2.3)$$

One of the most important applications of the Mordukhovich derivatives of single-valued mappings is to define the covering constants, which is recalled as follows. Let $f: X \to Y$ be a single-valued mapping. For given $\bar{x} \in X$ and $\bar{y} \in Y$ with $\bar{y} = f(\bar{x})$, the covering constant for $f$ at $(\bar{x}, \bar{y})$ is defined by

$$\hat{\alpha}(f,\bar{x},\bar{y}) = \sup_{\eta>0} \inf\left\{ \|z^*\|_{X^*}: z^* \in \widehat{D}^*f(x,y)(w^*), x \in \mathbb{B}_X(\bar{x},\eta), y = f(x) \in \mathbb{B}_Y(\bar{y},\eta), \|w^*\|_{Y^*} = 1 \right\}. (2.4)$$

## 3. Mordukhovich Differentiability of the Metric Projection Operator in $l_1$

Throughout this paper, we always let $\mathbb{N}$, $\mathbb{R}$ and $\mathbb{R}_+$ respectively denote the set of nonnegative integers, the set of real numbers and the set of nonnegative real numbers. For any positive integer $n$, let $\mathbb{R}^n$ denote the $n$-d Euclidean space. In particular, $\mathbb{R}^1 = \mathbb{R}$. Let $\mathbb{S}$ denote the linear space of real sequences with term-by-term addition and scalar multiplication.

Let $(l_1, \|\cdot\|)$ denote the ordinary real Banach space of absolutely summable sequences of real numbers, in which the norm $\|\cdot\|$ on $l_1$ is defined by

$$\|x\| = \sum_{n=1}^{\infty}|x_n| < \infty, \text{ for any } x = (x_1, x_2, \dots) \in l_1.$$

The topological dual space of $(l_1, \|\cdot\|)$ is the Banach space $(l_\infty, \|\cdot\|_\infty)$, in which the supremum-norm $\|\cdot\|_\infty$ on $l_\infty$ is defined, for any $y^* = (y_1^*, y_2^*, \dots) \in l_\infty$, by

$$\|y^*\|_\infty = \sup\{|y_1^*|, |y_2^*|, |y_3^*|, \dots\}.$$

The real pairing between $l_\infty$ and $l_1$ is written by $\langle\cdot,\cdot\rangle$ and

$$\langle y^*, x\rangle = \sum_{n=1}^{\infty} x_n y_n^*, \text{ for any } x = (x_1, x_2, \dots) \in l_1 \text{ and } y^* = (y_1^*, y_2^*, \dots) \in l_\infty.$$

Their null element is denoted by $\theta = (0, 0, \dots)$. Let $K$ be the positive cone in $l_1$and

$$K = \{x = (x_1, x_2, \dots) \in l_1 : x_n \geq 0, \text{for } n = 1, 2, \dots\}.$$

$K$ is a nonempty convex and closed pointed cone in $l_1$, which induces a partial order $\preccurlyeq$ on $l_1$. More precisely, $\preccurlyeq$ is defined, for $u = (u_1, u_2, \dots) \in l_1$ and $x = (x_1, x_2, \dots) \in l_1$ by

$$u \preccurlyeq x \quad \text{if and only if } \ u_n \leq x_n, \text{ for } n = 1, 2, \dots.$$

and

$$u \prec x \quad \text{if and only if } \ u_n < x_n, \text{ for } n = 1, 2, \dots.$$

Let $\mathbb{B}$ and $\mathbb{B}^*$ denote the unit balls in $l_1$ and $l_\infty$, respectively, which satisfy

$$\mathbb{B} = \{x \in l_1 : \|x\| \leq 1\}.$$

The set-valued metric projection operator $P_{\mathbb{B}}: l_1 \rightrightarrows \mathbb{B}$ is defined, for any $x = (x_1, x_2, \dots) \in l_1$, by

$$P_{\mathbb{B}}(x) = \{u \in \mathbb{B} : \|x - u\| \leq \|x - v\|, \text{for each } v \in \mathbb{B}\}.$$

Notations: Let $a$ and $b$ be real numbers, as usual, we denote $a \wedge b = \min\{a, b\}$, $a \vee b = \max\{a, b\}$ and $[a, a] = \{a\}$. In particular, for any $x = (x_1, x_2, \dots) \in l_1$, we write $x \wedge 1 = (x_1 \wedge 1, x_2 \wedge 1, \dots) \in l_1$ and $x \vee 1 = (x_1 \vee 1, x_2 \vee 1, \dots)$. For $y^* = (y_1^*, y_2^*, \dots) \in l_\infty$, we write $y^* \vee a = (y_1^* \vee a, y_2^* \vee a, \dots) \in l_\infty$ and $\beta \coloneqq (1, 1, \dots) \in l_\infty$.

In the following theorem, we find the explicit solutions of $P_{\mathbb{B}}: l_1 \to 2^{\mathbb{B}}$.

**Theorem 3.1**. *Let $x = (x_1, x_2, \dots) \in l_1$ with $\|x\| \geq 1$. The set-valued metric projection operator $P_{\mathbb{B}}: l_1 \rightrightarrows \mathbb{B}$ has the following properties.*

$$P_{\mathbb{B}}(x) = \left\{y \in \mathbb{B} : \|y\| = 1, \|x - y\| = \|x\| - 1, \begin{matrix} 0 \leq y_n \leq x_n \wedge 1, & \text{if } x_n \geq 0, \\ x_n \vee (-1) \leq y_n \leq 0, & \text{if } x_n < 0, \end{matrix} \ for\ n \in \mathbb{N}\right\}.$$

*In particular,*

(I) $\frac{x}{\|x\|} \in P_{\mathbb{B}}(x)$, *for any* $x \in l_1$ *with* $\|x\| \geq 1$.

(II) *If* $x \in K$ *with* $\|x\| \geq 1$, *then*

$$P_{\mathbb{B}}(x) = \{y \in \mathbb{B}: \sum_{n=1}^{\infty} y_n = 1, \|x - y\| = \|x\| - 1, \ 0 \leq y_n \leq x_n \wedge 1, \textit{for each } n \ \in \mathbb{N}\}.$$

(III) *If* $x \in -K$ *with* $\|x\| \geq 1$, *then*

$$P_{\mathbb{B}}(x) = \{y \in \mathbb{B}: \sum_{n=1}^{\infty} y_n = -1, \|x - y\| = \|x\| - 1, x_n \vee (-1) \leq y_n \leq 0, \textit{for each } n \ \in \mathbb{N}\}.$$

*Proof*. For a given $x = (x_1, x_2, ...) \in l_1$. If $\|x\| = 1$, then $P_{\mathbb{B}}(x) = \{x\} \coloneqq x$ and this theorem is proved. Hence, we assume that $\|x\| > 1$. Let $z = \frac{x}{\|x\|} = \left\{\frac{x_n}{\|x\|}\right\}_{n=1}^{\infty}$. It is clear that $z \in l_1$ and $\|z\| = 1$, which implies $z \in \mathbb{B}$. We calculate

$$\|x - z\| = \left\|x - \frac{x}{\|x\|}\right\| = \frac{1}{\|x\|}(\|x\| - 1)\|x\| = \|x\| - 1. \tag{3.1}$$

Let arbitrarily $y = (y_1, y_2, ...) \in \mathbb{B}$. Then, it satisfies that $\|y\| = \sum_{n=1}^{\infty} |y_n| \leq 1$. We calculate

$$\|x - y\| = \sum_{n=1}^{\infty} |x_n - y_n| \geq \sum_{n=1}^{\infty} |x_n| - \sum_{n=1}^{\infty} |y_n| \geq \sum_{n=1}^{\infty} |x_n| - 1 = \|x\| - 1. \tag{3.2}$$

On the other hand, by $z \in \mathbb{B}$ and by (3.2), we have

$$y \in \mathbb{B} \quad \Longrightarrow \quad \|x - y\| \geq \|x\| - 1 = \|x - z\|. \tag{3.3}$$

By (3.1) and (3.2) and by definition of $P_{\mathbb{B}}$, we obtain that

$$z = \frac{x}{\|x\|} \in P_{\mathbb{B}}(x).$$

By (3.1) and (3.3), we obtain that

$$y \in P_{\mathbb{B}}(x) \quad \Longrightarrow \quad \|x - y\| = \|x - z\| = \|x\| - 1. \tag{3.4}$$

Next, we prove

$$y \in P_{\mathbb{B}}(x) \quad \Longrightarrow \quad \begin{matrix} 0 \leq y_n \leq x_n, & \text{if } x_n \geq 0, \\ x_n \leq y_n \leq 0, & \text{if } x_n < 0, \end{matrix} \quad \text{for } n = 1, 2, \dots . \tag{3.5}$$

In particular, by (3.5), we have

$$y \in P_{\mathbb{B}}(x) \quad \Longrightarrow \quad |y_n| \leq |x_n|, \quad \text{for } n = 1, 2, \dots . \tag{3.6}$$

To prove (3.5) and (3.6), we define the following partition of $\mathbb{N}$ with respect to $x = (x_1, x_2, ...) \in l_1$ and $y = (y_1, y_2, ...) \in P_{\mathbb{B}}(x)$.

$$A = \{n \in \mathbb{N}: x_n \geq y_n \geq 0 \ \text{ and } \ x_n \geq 0\} \quad \text{and} \quad B = \{n \in \mathbb{N}: x_n \leq y_n \ \leq 0 \text{ and } \ x_n \leq 0\},$$

$$C = \{n \in \mathbb{N}: y_n > x_n \geq 0\} \quad \text{and} \quad D = \{n \in \mathbb{N}: y_n < x_n \leq 0\},$$

$$E = \{n \in \mathbb{N}: y_n < 0 \ \text{ and } \ x_n \geq 0\} \quad \text{and} \quad F = \{n \in \mathbb{N}: y_n > 0 \ \text{ and } \ x_n \leq 0\}.$$

In order to prove (3.5), it equivalently proves $E = \emptyset$, $F = \emptyset$, $C = \emptyset$ and $D = \emptyset$. To this end, we first

prove

$$y \in P_{\mathbb{B}}(x) \quad \Longrightarrow \quad E = \emptyset \quad \text{and} \quad F = \emptyset. \tag{3.7}$$

Proof of (3.7). Suppose $y \in P_{\mathbb{B}}(x)$. Assume, by the way of contradiction, that $E \neq \emptyset$. Then, there is $m \in \mathbb{N}$ such that $y_m < 0$ and $x_m \geq 0$. Define $v = (v, v_2, \ldots) \in l_1$ as follows

$$v_n = \begin{cases} y_n, & \text{if } n \neq m, \\ 0, & \text{if } n = m. \end{cases}$$

It is clear that $\|y\| = \|v\| + |y_m| > \|v\|$, which implies $v \in \mathbb{B}$. Since $x_m \geq 0$, with $y_m < 0$, we have

$$\|x - y\| = \sum_{n \neq m} |y_n - x_n| + |y_m - x_m| > \sum_{n \neq m} |y_n - x_n| + |0 - x_m| = \|x - v\|, \text{ with } v \in \mathbb{B}.$$

This contradicts to the assumption that $y \in P_{\mathbb{B}}(x)$, which proves $E = \emptyset$.

Next suppose $y \in P_{\mathbb{B}}(x)$ and assume, by the way of contradiction, that $F \neq \emptyset$. Then, there is $k \in \mathbb{N}$ such that $y_k > 0$ and $x_k \leq 0$. Define $w = (w_1, w_2, \ldots) \in l_1$ as follows

$$w_n = \begin{cases} y_n, & \text{if } n \neq k, \\ 0, & \text{if } n = k. \end{cases}$$

It is clear that $\|y\| = \|w\| + |y_k| > \|w\|$, which implies $w \in \mathbb{B}$. Since $x_k \leq 0$, with $y_k > 0$, we have

$$\|x - y\| = \sum_{n \neq k} |y_n - x_n| + |y_k - x_k| > \sum_{n \neq m} |y_n - x_n| + |0 - x_k| = \|x - w\|, \text{ with } w \in \mathbb{B}.$$

This contradicts to the assumption that $y \in P_{\mathbb{B}}(x)$, which proves $F = \emptyset$. Hence, (3.7) is proved. Next, under the conditions (3.7), we secondly prove

$$y \in P_{\mathbb{B}}(x) \quad \Longrightarrow \quad C = \emptyset \quad \text{and} \quad D = \emptyset. \tag{3.8}$$

(3.7) implies that $\mathbb{N} = A \cup B \cup C \cup D$. Then, by the assumption that $y \in P_{\mathbb{B}}(x)$ and by (3.4), we calculate

$$\begin{aligned}
&\|x - y\| \\
&= \textstyle\sum_{n=1}^{\infty} |x_n - y_n| \\
&= \textstyle\sum_{n\in A} |x_n - y_n| + \sum_{n\in B} |x_n - y_n| + \sum_{n\in C} |x_n - y_n| + \sum_{n\in D} |x_n - y_n| \\
&= \textstyle\sum_{n\in A} (x_n - y_n) + \sum_{n\in B} (y_n - x_n) + \sum_{n\in C} (y_n - x_n) + \sum_{n\in D} (x_n - y_n) \\
&= \textstyle\sum_{n\in A} x_n - \sum_{n\in A} y_n + \sum_{n\in B} y_n - \sum_{n\in B} x_n + \sum_{n\in C} y_n - \sum_{n\in C} x_n + \sum_{n\in D} x_n - \sum_{n\in D} y_n \\
&= \textstyle\sum_{n\in A} x_n - \sum_{n\in B} x_n - \sum_{n\in C} x_n + \sum_{n\in D} x_n - \sum_{n\in A} y_n + \sum_{n\in B} y_n + \sum_{n\in C} y_n - \sum_{n\in D} y_n \\
&= \textstyle\sum_{n\in A} x_n - \sum_{n\in B} x_n + \sum_{n\in C} x_n - \sum_{n\in D} x_n - 2\sum_{n\in C} x_n + 2\sum_{n\in D} x_n \\
&\quad \textstyle - \sum_{n\in A} y_n + \sum_{n\in B} y_n - \sum_{n\in C} y_n + \sum_{n\in D} y_n + 2\sum_{n\in C} y_n - 2\sum_{n\in D} y_n \\
&= \textstyle\|x\| - 2\sum_{n\in C} x_n + 2\sum_{n\in D} x_n - \|y\| + 2\sum_{n\in C} y_n - 2\sum_{n\in D} y_n \\
&= \textstyle\|x\| - 1 - 2\sum_{n\in C} x_n + 2\sum_{n\in D} x_n + 2\sum_{n\in C} y_n - 2\sum_{n\in D} y_n
\end{aligned}$$

$= \|x\| - 1 - 2\sum_{n\in C}(x_n - y_n) - 2\sum_{n\in D}(y_n - x_n)$

$= \|x\| - 1 + 2\sum_{n\in C}(y_n - x_n) + 2\sum_{n\in D}(x_n - y_n)$

$= \|x\| - 1 + 2\sum_{n\in C}|y_n - x_n| + 2\sum_{n\in D}|x_n - y_n|$

By (3.4), this implies that

$$y \in P_{\mathbb{B}}(x) \quad \Longrightarrow \quad \sum_{n\in C}|y_n - x_n| = 0 \ \text{ and } \ \sum_{n\in D}|x_n - y_n| = 0. \tag{3.9}$$

By the definitions of $C$ and $D$, (3.9) is equivalent to

$$y \in P_{\mathbb{B}}(x) \quad \Longrightarrow \quad C = \emptyset \quad \text{and} \quad D = \emptyset.$$

By the definitions of $A$ and $B$, (3.5) is proved by (3.7) and (3.8). Summarizing (3.4) and (3.5), we have that if $y \in P_{\mathbb{B}}(x)$, then

$$\|x - y\| = \|x\| - 1, \quad \text{and} \quad \begin{matrix} 0 \le y_n \le x_n, & \text{if } x_n \ge 0, \\ x_n \le y_n \le 0, & \text{if } x_n < 0, \end{matrix} \quad \text{for } n = 1, 2, \dots . \tag{3.10}$$

Under conditions (3.4) and (3.8), or (3.10), we prove that

$$y \in \mathbb{B} \ \text{ and } y \in P_{\mathbb{B}}(x) \ \Longrightarrow \ \|y\| = 1. \tag{3.11}$$

Proof of (3.11). Let $y \in \mathbb{B}$. Suppose $y \in P_{\mathbb{B}}(x)$. Let $z \in \mathbb{B}$ be given in (3.1). (3.7) and (3.8) imply that $\mathbb{N} = A \cup B$. By (3.1) and (3.4), and by the definitions of the sets $A$ and $B$, we have

$$\begin{aligned}
&\|x\| - 1 \\
&= \|x - z\| \\
&\ge \|x - y\| \\
&= \textstyle\sum_{n=1}^{\infty}|x_n - y_n| \\
&= \textstyle\sum_{n\in A}|x_n - y_n| + \sum_{n\in B}|x_n - y_n| \\
&= \textstyle\sum_{n\in A}(x_n - y_n) + \sum_{n\in B}(y_n - x_n) \\
&= \textstyle\sum_{n\in A} x_n + \sum_{n\in B}(-x_n) - \sum_{n\in A} y_n - \sum_{n\in B}(-y_n) \\
&= \textstyle\sum_{n=1}^{\infty}|x_n| - \sum_{n=1}^{\infty}|y_n| \\
&= \|x\| - \|y\|.
\end{aligned}$$

This implies that $\|y\| \ge 1$. Since $y \in \mathbb{B}$, then $\|y\| \le 1$, which proves (3.11). By (3.10) and (3.11), we obtain that, if $y \in P_{\mathbb{B}}(x)$ then

$$\|x - y\| = \|x\| - 1, \ \|y\| = 1, \ \text{and} \ \begin{matrix} 0 \le y_n \le x_n, & \text{if } x_n \ge 0, \\ x_n \le y_n \le 0, & \text{if } x_n < 0, \end{matrix} \quad \text{for } n = 1, 2, \dots . \tag{3.12}$$

Conversely, we prove that, for $y \in \mathbb{B}$, if $y$ satisfies (3.12), then $y \in P_{\mathbb{B}}(x)$. By (3.2), we already proved

$$\|x - v\| \geq \|x\| - 1, \text{ for any } v \in \mathbb{B}. \tag{3.13}$$

Let $y \in \mathbb{B}$. Suppose that $y$ satisfies (3.12). By $\mathbb{N} = A \cup B$, we calculate

$$\begin{aligned}
&\|x - y\| \\
&= \sum_{n=1}^{\infty} |x_n - y_n| \\
&= \sum_{n \in A} |x_n - y_n| + \sum_{n \in B} |x_n - y_n| \\
&= \sum_{n \in A} (x_n - y_n) + \sum_{n \in B} (y_n - x_n) \\
&= \sum_{n \in A} x_n + \sum_{n \in B} (-x_n) - \sum_{n \in A} y_n - \sum_{n \in B} (-y_n) \\
&= \sum_{n=1}^{\infty} |x_n| - \sum_{n=1}^{\infty} |y_n| \\
&= \|x\| - \|y\| \\
&= \|x\| - 1.
\end{aligned}$$

By (3.13), this proves

$$y \in \mathbb{B} \text{ and } y \text{ satisfies (3.12)} \quad \Longrightarrow \quad y \in P_{\mathbb{B}}(x). \tag{3.14}$$

By (3.12) and (3.14), for any $y \in \mathbb{B}$,

$$y \in P_{\mathbb{B}}(x) \iff \|x - y\| = \|x\| - 1, \ \|y\| = 1 \text{ and } \begin{matrix} 0 \leq y_n \leq x_n, & \text{if } x_n \geq 0, \\ x_n \leq y_n \leq 0, & \text{if } x_n < 0, \end{matrix} \text{ for } n = 1, 2, \ldots \tag{3.15}$$

By the condition $\|y\| = 1$ in (3.15), it can be rewritten as, for any $y \in \mathbb{B}$,

$$y \in P_{\mathbb{B}}(x) \iff \|x - y\| = \|x\| - 1, \|y\| = 1 \text{ and } \begin{matrix} 0 \leq y_n \leq x_n \wedge 1, & \text{if } x_n \geq 0, \\ x_n \vee (-1) \leq y_n \leq 0, & \text{if } x_n < 0, \end{matrix} \text{ for } n = 1, 2, \ldots \tag{3.16}$$

Part (I) is proved by (3.1) and (3.16). Parts (II) and (III) are indeed special cases of (3.16), which proves this theorem. □

Now, we use the representation of the set-valued metric projection operator $P_{\mathbb{B}}$ to study its Mordukhovich differentiability. Since the complexity of the Mordukhovich derivative of $P_{\mathbb{B}}$, we do not find the explicit solutions of $\widehat{D}^*P_{\mathbb{B}}(\bar{x}, \bar{y})$ in general. However, in next section, we investigate the explicit solutions of $\widehat{D}^*P_{\mathbb{B}}(\bar{x}, \bar{y})$ in a 2-diensiaonal Banach space, from which, we see that is very complicated to find the explicit solutions of $\widehat{D}^*P_{\mathbb{B}}(\bar{x}, \bar{y})$.

**Theorem 3.2**. *Let $\bar{x} = (\bar{x}_1, \bar{x}_2, \ldots) \in l_1$ with $\|\bar{x}\| \geq 1$ and $\bar{y} = (\bar{y}_1, \bar{y}_2, \ldots) \in P_{\mathbb{B}}(\bar{x})$. Then for a given $y^* = (y_1^*, y_2^*, \ldots) \in l_\infty$, the Mordukhovich derivative $\widehat{D}^*P_{\mathbb{B}}(\bar{x}, \bar{y})(y^*)$ has the following properties.*

(A) *If $\bar{x} \in K$, then,*

$$\widehat{D}^*P_{\mathbb{B}}(\bar{x}, \bar{y})(y^*) \subseteq (-K).$$

(B) *If $\bar{x} \in -K$, then,*

$$\widehat{D}^*P_{\mathbb{B}}(\bar{x}, \bar{y})(y^*) \subseteq K.$$

(C) *If* $y_1^* = y_2^* = \cdots$, *then*

(I) $\sum_{n=1}^{\infty} \bar{y}_n = 1$ and $y_1^* \le 0 \quad \Longrightarrow \quad \theta \in \widehat{D}^* P_{\mathbb{B}}(\bar{x}, \bar{y})(y^*)$.

(II) $\sum_{n=1}^{\infty} \bar{y}_n = -1$ and $y_1^* \ge 0 \quad \Longrightarrow \quad \theta \in \widehat{D}^* P_{\mathbb{B}}(\bar{x}, \bar{y})(y^*)$.

(III) $\bar{x} \in K$ *and* $y_1^* \le 0 \quad \Longrightarrow \quad \theta \in \widehat{D}^* P_{\mathbb{B}}(\bar{x}, \bar{y})(y^*)$.

(IV) $\bar{x} \in -K$ *and* $y_1^* \ge 0 \quad \Longrightarrow \quad \theta \in \widehat{D}^* P_{\mathbb{B}}(\bar{x}, \bar{y})(y^*)$.

(V) $\bar{x} \in K \quad \Longrightarrow \quad \theta \in \widehat{D}^* P_{\mathbb{B}}(\bar{x}, \bar{y})(-\beta)$.

(VI) $\bar{x} \in -K \quad \Longrightarrow \quad \theta \in \widehat{D}^* P_{\mathbb{B}}(\bar{x}, \bar{y})(\beta)$.

(D) *If* $y_1^* = y_2^* = \cdots$, *with* $y_1^* \neq 0$, *then*

$$-1 < \sum_{n=1}^{\infty} \bar{y}_n < 1 \quad \Longrightarrow \quad \theta \notin \widehat{D}^* P_{\mathbb{B}}(\bar{x}, \bar{y})(y^*).$$

(E) *For any* $\bar{x} = (\bar{x}_1, \bar{x}_2, \ldots) \in \mathbb{B}^0$, *we have* $P_{\mathbb{B}}(\bar{x}) = \bar{x}$ *and*

$$\widehat{D}^* P_{\mathbb{B}_l}(\bar{x}, \bar{x})(y^*) = y^*, \text{ for any } y^* = (y_1^*, y_2^*, \ldots) \in l_\infty.$$

*Proof.* Let $\bar{x} = (\bar{x}_1, \bar{x}_2, \ldots) \in l_1$ with $\|\bar{x}\| \ge 1$. Let $\bar{y} = (\bar{y}_1, \bar{y}_2, \ldots) \in P_{\mathbb{B}}(\bar{x})$. Let $y^* = (y_1^*, y_2^*, \ldots) \in l_\infty$ and $x^* = (x_1^*, x_2^*, \ldots) \in l_\infty$. We calculate

$$\limsup_{\substack{(u,v)\to(\bar{x},\bar{y}) \\ v=(v_1,v_2,\ldots)\in P_{\mathbb{B}}(u)}} \frac{\langle x^*, u-\bar{x}\rangle - \langle y^*, v-\bar{y}\rangle}{\|u-\bar{x}\| + \|v-\bar{y}\|_\infty}$$

$$= \limsup_{\substack{(u,v)\to(\bar{x},\bar{y}) \\ v=(v_1,v_2,\ldots)\in P_{\mathbb{B}}(u)}} \frac{\langle (x_1^*, x_2^*, \ldots),\ (u_1, u_2, \ldots) - (\bar{x}_1, \bar{x}_2, \ldots)\rangle - \langle (y_1^*, y_2^*, \ldots),\ (v_1, v_2, \ldots) - (\bar{y}_1, \bar{y}_2, \ldots)\rangle}{\|u-\bar{x}\| + \|v-\bar{y}\|_\infty}$$

$$= \limsup_{\substack{(u,v)\to(\bar{x},\bar{y}) \\ v=(v_1,v_2,\ldots)\in P_{\mathbb{B}}(u)}} \frac{\sum_{n=1}^{\infty} x_n^*(u_n - \bar{x}_n) - \sum_{n=1}^{\infty} y_n^*(v_n - \bar{y}_n)}{\|u-\bar{x}\| + \|v-\bar{y}\|_\infty}. \tag{3.17}$$

Proof of (A). Let $\bar{x} = (\bar{x}_1, \bar{x}_2, \ldots) \in K$ and $\bar{y} = (\bar{y}_1, \bar{y}_2, \ldots) \in P_{\mathbb{B}}(\bar{x})$. Then, $\bar{x}_n \ge 0$, for $n = 1, 2, \ldots$ . By Theorem 3.1, $0 \le \bar{y}_n \le \bar{x}_n \wedge 1$, for $n = 1, 2, \ldots$ .

Let $x^* = (x_1^*, x_2^*, \ldots) \notin -K$. Then there is a positive number $m$ such that $x_m^* > 0$. For $t > 0$, define $u(t) = (u_1, u_2, \ldots) \in l_1$, for $n = 1, 2, \ldots$ , by

$$u_n = \begin{cases} \bar{x}_n, & \text{if } n \neq m, \\ \bar{x}_m + t, & \text{if } n = m. \end{cases} \tag{3.18}$$

Let $v = (v_1, v_2, \ldots) = \bar{y} = (\bar{y}_1, \bar{y}_2, \ldots)$. By $\bar{x} = (\bar{x}_1, \bar{x}_2, \ldots) \in K$ with $\bar{x}_n \ge 0$, for $n = 1, 2, \ldots$ , we have

$$\|u(t) - v\| = \|u(t) - \bar{y}\| = \|\bar{x} - \bar{y}\| + t = \|\bar{x}\| - 1 + t = \|u(t)\| - 1.$$

We see that $\|v\| = \|\bar{y}\| = 1$ and $0 \le v_n = \bar{y}_n \le \bar{x}_n \wedge 1 \le u_n \wedge 1$, for $n = 1, 2, \ldots$ . By Theorem 3.1, this proves that

$$v = \bar{y} \in P_{\mathbb{B}}(u(t)), \text{ for } t > 0. \tag{3.19}$$

Then, substituting (3.18) and $v = \bar{y}$ to (3.17), we obtain

$$\limsup_{\substack{(u,v)\to(\bar{x},\bar{y})\\ v=(v_1,v_2,\dots)\in P_{\mathbb{B}}(u)}} \frac{\langle x^*, u-\bar{x}\rangle - \langle y^*, v-\bar{y}\rangle}{\|u-\bar{x}\| + \|v-\bar{y}\|_\infty}$$

$$\geq \limsup_{\substack{(u(t),v)\to(\bar{x},\bar{y})\\ v=\bar{y}\in P_{\mathbb{B}}(u(t)),t\downarrow 0}} \frac{x_m^* t}{t}$$

$$= x_m^* > 0.$$

This implies that under the conditions that $\bar{x} \in K$ and $\bar{y} \in P_{\mathbb{B}}(\bar{x})$, for any given $y^* = (y_1^*, y_2^*, \dots) \in l_\infty$, for $x^* = (x_1^*, x_2^*, \dots)$, if there is a positive number $m$ such that $x_m^* > 0$, then

$$x^* \notin \widehat{D}^* P_{\mathbb{B}}(\bar{x}, \bar{y})(y^*), \text{ for any } x^* \notin -K.$$

This proves part (A).

Proof of (B). Let $\bar{x} = (\bar{x}_1, \bar{x}_2, \dots) \in -K$ and $\bar{y} = (\bar{y}_1, \bar{y}_2, \dots) \in P_{\mathbb{B}}(\bar{x})$. Then, $\bar{x}_n \leq 0$, for $n$ = 1, 2, … . By Theorem 3.1, $\bar{x}_n \vee (-1) \leq \bar{y}_n \leq 0$, for $n$ = 1, 2, … .

Let $x^* = (x_1^*, x_2^*, \dots) \notin K$. Then there is a positive number $m$ such that $x_m^* < 0$. For $t > 0$, define $u(t) = (u_1, u_2, \dots) \in l_1$, for $n$ = 1, 2, … , by

$$u_n = \begin{cases} \bar{x}_n, & \text{if } n \neq m, \\ \bar{x}_m - t, & \text{if } n = m. \end{cases} \tag{3.20}$$

Let $v = (v_1, v_2, \dots) = \bar{y} = (\bar{y}_1, \bar{y}_2, \dots)$. By the condition that $\bar{x}_m \leq 0$, similarly to the proof of (3.19), we can show that

$$v = \bar{y} \in P_{\mathbb{B}}(u(t)), \text{ for } t > 0.$$

Then, substituting (3.20) and $v = \bar{y}$ to (3.17), we obtain

$$\limsup_{\substack{(u,v)\to(\bar{x},\bar{y})\\ v=(v_1,v_2,\dots)\in P_{\mathbb{B}}(u)}} \frac{\langle x^*, u-\bar{x}\rangle - \langle y^*, v-\bar{y}\rangle}{\|u-\bar{x}\| + \|v-\bar{y}\|_\infty}$$

$$\geq \limsup_{\substack{(u(t),v)\to(\bar{x},\bar{y})\\ v=\bar{y}\in P_{\mathbb{B}}(u(t)),t\downarrow 0}} \frac{x_m^*(-t)}{|t|}$$

$$= -x_m^* > 0.$$

This implies that under the conditions that $\bar{x} \in -K$ and $\bar{y} \in P_{\mathbb{B}}(\bar{x})$, for any given $y^* = (y_1^*, y_2^*, \dots) \in l_\infty$, for $x^* = (x_1^*, x_2^*, \dots)$, if there is a positive number $m$ such that $x_m^* < 0$, then

$$x^* \notin \widehat{D}^* P_{\mathbb{B}}(\bar{x}, \bar{y})(y^*), \text{ for any } x^* \notin K.$$

This proves part (B).

Proof (I) of (C). Let $\bar{x} \in l_1$ with $\|\bar{x}\| \geq 1$. Let $\bar{y} \in P_{\mathbb{B}}(\bar{x})$ such that $\sum_{n=1}^{\infty} \bar{y}_n = 1$. Suppose $y_1^* \leq 0$. Let $u = (u_1, u_2, \ldots) \in l_1$ with $\|u\| \geq 1$ and $v = \in P_{\mathbb{B}}(u)$. By Theorem 3.1, $\|v\| = 1$. Substituting $x^*$ by $\theta$ in (3.17), we have

$$\limsup_{\substack{(u,v)\to(\bar{x},\bar{y}) \\ v=(v_1,v_2,\ldots)\in P_{\mathbb{B}}(u)}} \frac{\sum_{n=1}^{\infty} x_n^*(u_n-\bar{x}_n)-\sum_{n=1}^{\infty} y_n^*(v_n-\bar{y}_n)}{\|u-\bar{x}\|+\|v-\bar{y}\|_\infty}$$

$$= \limsup_{\substack{(u,v)\to(\bar{x},\bar{y}) \\ v=(v_1,v_2,\ldots)\in P_{\mathbb{B}}(u)}} \frac{-y_1^*\sum_{n=1}^{\infty}(v_n-\bar{y}_n)}{\|u-\bar{x}\|+\|v-\bar{y}\|_\infty}$$

$$= \limsup_{\substack{(u,v)\to(\bar{x},\bar{y}) \\ v=(v_1,v_2,\ldots)\in P_{\mathbb{B}}(u)}} \frac{-y_1^*(\sum_{n=1}^{\infty} v_n-\sum_{n=1}^{\infty}\bar{y}_n)}{\|u-\bar{x}\|+\|v-\bar{y}\|_\infty}$$

$$= \limsup_{\substack{(u,v)\to(\bar{x},\bar{y}) \\ v=(v_1,v_2,\ldots)\in P_{\mathbb{B}}(u)}} \frac{-y_1^*(\sum_{n=1}^{\infty} v_n-\sum_{n=1}^{\infty}\bar{y}_n)}{\|u-\bar{x}\|+\|v-\bar{y}\|_\infty}$$

$$= \limsup_{\substack{(u,v)\to(\bar{x},\bar{y}) \\ v=(v_1,v_2,\ldots)\in P_{\mathbb{B}}(u)}} \frac{-y_1^*(\sum_{n=1}^{\infty} v_n-1)}{\|u-\bar{x}\|+\|v-\bar{y}\|_\infty}$$

$$\leq \limsup_{\substack{(u,v)\to(\bar{x},\bar{y}) \\ v=(v_1,v_2,\ldots)\in P_{\mathbb{B}}(u)}} \frac{-y_1^*(\sum_{n=1}^{\infty}|v_n|-1)}{\|u-\bar{x}\|+\|v-\bar{y}\|_\infty} \qquad (y_1^* \leq 0)$$

$$= \limsup_{\substack{(u,v)\to(\bar{x},\bar{y}) \\ v\in P_{\mathbb{B}}(u)}} \frac{-y_1^*(1-1)}{\|u-\bar{x}\|+\|v-\bar{y}\|_M} = 0.$$

This proves (I) of (C).

Proof (II) of (C). Let $\bar{x} \in l_1$ with $\|\bar{x}\| \geq 1$. Let $\bar{y} \in P_{\mathbb{B}}(\bar{x})$ such that $\sum_{n=1}^{\infty} \bar{y}_n = -1$. Suppose $y_1^* \geq 0$. Substituting $x^*$ by $\theta$ in (3.17), we have

$$\limsup_{\substack{(u,v)\to(\bar{x},\bar{y}) \\ v=(v_1,v_2,\ldots)\in P_{\mathbb{B}}(u)}} \frac{\sum_{n=1}^{\infty} x_n^*(u_n-\bar{x}_n)-\sum_{n=1}^{\infty} y_n^*(v_n-\bar{y}_n)}{\|u-\bar{x}\|+\|v-\bar{y}\|_\infty}$$

$$= \limsup_{\substack{(u,v)\to(\bar{x},\bar{y}) \\ v=(v_1,v_2,\ldots)\in P_{\mathbb{B}}(u)}} \frac{-y_1^*\sum_{n=1}^{\infty}(v_n-\bar{y}_n)}{\|u-\bar{x}\|+\|v-\bar{y}\|_\infty}$$

$$= \limsup_{\substack{(u,v)\to(\bar{x},\bar{y}) \\ v=(v_1,v_2,\ldots)\in P_{\mathbb{B}}(u)}} \frac{-y_1^*(\sum_{n=1}^{\infty} v_n-\sum_{n=1}^{\infty}\bar{y}_n)}{\|u-\bar{x}\|+\|v-\bar{y}\|_\infty}$$

$$= \limsup_{\substack{(u,v)\to(\bar{x},\bar{y}) \\ v=(v_1,v_2,\ldots)\in P_{\mathbb{B}}(u)}} \frac{-y_1^*(\sum_{n=1}^{\infty} v_n+1)}{\|u-\bar{x}\|+\|v-\bar{y}\|_\infty}$$

$$\leq \limsup_{\substack{(u,v)\to(\bar{x},\bar{y}) \\ v=(v_1,v_2,\dots)\in P_{\mathbb{B}}(u)}} \frac{-y_1^*(-\sum_{n=1}^{\infty}|v_n|+1)}{\|u-\bar{x}\|+\|v-\bar{y}\|_\infty} \qquad (y_1^* \geq 0)$$

$$= \limsup_{\substack{(u,v)\to(\bar{x},\bar{y}) \\ v\in P_{\mathbb{B}}(u)}} \frac{-y_1^*(-1+1)}{\|u-\bar{x}\|+\|v-\bar{y}\|_\infty} = 0.$$

This proves (II) of (C). All other parts (III) and (IV) of (C) follow from (I) and (II) immediately.

Proof of (D). Let $\bar{x} = (\bar{x}_1, \bar{x}_2, \dots) \in l_1$ with $\|\bar{x}\| \geq 1$. Let $\bar{y} = (\bar{y}_1, \bar{y}_2, \dots) \in P_{\mathbb{B}}(\bar{x})$. Suppose that $-1 < \sum_{n=1}^{\infty} \bar{y}_n < 1$. By Theorem 3.1, $\sum_{n=1}^{\infty}|\bar{y}_n| = 1$, this implies that there are two positive numbers *m* and *k* such that $\bar{y}_m > 0$ and $\bar{y}_k < 0$. By Theorem 3.1 again, this implies that $\bar{x}_m > 0$ and $\bar{x}_k < 0$. By the conditions $-1 < \sum_{n=1}^{\infty} \bar{y}_n < 1$ and $\sum_{n=1}^{\infty}|\bar{y}_n| = 1$, it is clear that $\bar{y}_m \neq 1$ and $\bar{y}_k \neq -1$.

Let $y^* = (y_1^*, y_2^*, \dots) \in l_\infty$. Suppose $y_1^* = y_2^* = \cdots$. We consider two cases with respect to $y_1^* \neq 0$.

Case 1. Suppose $y_1^* < 0$. Let $\delta \coloneqq \min\{\bar{y}_m, 1-\bar{y}_m, -\bar{y}_k, 1+\bar{y}_k\}$.

For real number *t* with $0 < t < \delta$, let $u(t) = (u_1, u_2, \dots) \in l_1$ be defined, for *n* = 1, 2, … , by

$$u_n = \begin{cases} \bar{x}_n, & \text{if } n \neq m, \\ \bar{x}_m + t, & \text{if } n = m, \\ \bar{x}_k + t, & \text{if } n = k. \end{cases} \tag{3.21}$$

Let $v = (v_1, v_2, \dots)$ be defined, for *n* = 1, 2, … , by

$$v_n = \begin{cases} \bar{y}_n, & \text{if } n \neq m, \\ \bar{y}_m + t, & \text{if } n = m, \\ \bar{y}_k + t, & \text{if } n = k. \end{cases} \tag{3.22}$$

By $\bar{y}_m > 0$, $\bar{y}_k < 0$, $\bar{x}_m > 0$, $\bar{x}_k < 0$, $\bar{y}_m \neq 1$ and $\bar{y}_k \neq -1$, with $0 < t < \delta$, we have

$$0 < \bar{y}_m + t < 1 \text{ and } -1 < \bar{y}_k + t < 0.$$

Then, for *n* = *m* or *k*, we have

$$0 \leq v_m = \bar{y}_m + t = (\bar{y}_m + t) \wedge 1 = u_m \wedge 1$$

$$\text{and} \quad 0 > v_k = \bar{y}_k + t = (\bar{y}_k + t) \vee (-1) = u_k \vee (-1). \tag{3.23}$$

By (3.22) and by $0 < t < \delta$, we have

$$\begin{aligned} \|v\| &= \textstyle\sum_{n=1}^{\infty}|v_n| \\ &= \textstyle\sum_{n\notin\{m,k\}}|\bar{y}_n| + |\bar{y}_m + t| + |\bar{y}_k + t| \\ &= \textstyle\sum_{n\notin\{m,k\}}|\bar{y}_n| + |\bar{y}_m| + t + |\bar{y}_k| - t \\ &= \textstyle\sum_{n=1}^{\infty}|\bar{y}_n| \\ &= 1. \end{aligned} \tag{3.24}$$

By (3.21), we have

$$\|u(t)\| = \sum_{n=1}^{\infty}|u_n|$$
$$= \sum_{n\notin\{m,k\}}|\bar{x}_n| + |\bar{x}_m + t| + |\bar{x}_k + t|$$
$$= \sum_{n\notin\{m,k\}}|\bar{x}_n| + |\bar{x}_m| + t + |\bar{x}_k| - t$$
$$= \sum_{n=1}^{\infty}|\bar{x}_n|$$
$$= \|\bar{x}\|. \tag{3.25}$$

By $\bar{y} \in P_{\mathbb{B}}(\bar{x})$, Theorem 3.1 and (3.25), we calculate

$$\|u(t) - v\|$$
$$= \sum_{n=1}^{\infty}|u_n - v_n|$$
$$= \sum_{n\notin\{m,k\}}|u_n - v_n| + |u_m - v_m| + |u_k - v_k|$$
$$= \sum_{n\notin\{m,k\}}|\bar{x}_n - \bar{y}_n| + |\bar{x}_m + t - (\bar{y}_m + t)| + |\bar{x}_k + t - (\bar{y}_k + t)|$$
$$= \sum_{n=1}^{\infty}|\bar{x}_n - \bar{y}_n|$$
$$= \|\bar{x} - \bar{y}\|$$
$$= \|\bar{x}\| - 1$$
$$= \|u(t)\| - 1. \tag{3.26}$$

By (3.23), (3.24), (3.25) and (3.26), we obtain that, $v \in P_{\mathbb{B}}(u(t))$, for $0 < t < \delta$. Substituting $x^*$ by $\theta$ in (3.17), we have

$$\limsup_{\substack{(u,v)\to(\bar{x},\bar{y})\\ v=(v_1,v_2,\dots)\in P_{\mathbb{B}}(u)}} \frac{\sum_{n=1}^{\infty} x_n^*(u_n-\bar{x}_n) - \sum_{n=1}^{\infty} y_n^*(v_n-\bar{y}_n)}{\|u-\bar{x}\| + \|v-\bar{y}\|_\infty}$$
$$= \limsup_{\substack{(u,v)\to(\bar{x},\bar{y})\\ v=(v_1,v_2,\dots)\in P_{\mathbb{B}}(u)}} \frac{y_1^* \sum_{n=1}^{\infty}(\bar{y}_n - v_n)}{\|u-\bar{x}\| + \|v-\bar{y}\|_\infty}$$
$$\geq \limsup_{\substack{(u(t),v)\to(\bar{x},\bar{y}),0<t<\delta\\ v\in P_{\mathbb{B}}(u),t\downarrow 0}} \frac{y_1^*\big((\bar{y}_m-(\bar{y}_m+t))+(\bar{y}_k-(\bar{y}_k+t))\big)}{2t+t}$$
$$= \limsup_{\substack{(u(t),v)\to(\bar{x},\bar{y}),0<t<\delta\\ v\in P_{\mathbb{B}}(u),t\downarrow 0}} \frac{-2ty_1^*}{2t+t}$$
$$= \frac{-2y_1^*}{3} > 0.$$

This implies that, for $\bar{x} \in l_1$ with $\|\bar{x}\| \geq 1$ and $\bar{y} = (\bar{y}_1, \bar{y}_2, \dots) \in P_{\mathbb{B}}(\bar{x})$ with $-1 < \sum_{n=1}^{\infty} \bar{y}_n < 1$, and

$y_1^* = y_2^* = \cdots$ ,

$$\theta \notin \widehat{D}^* P_{\mathbb{B}}(\bar{x}, \bar{y})(y^*), \text{ for } y_1^* < 0. \tag{3.27}$$

Case 2. Suppose $y_1^* > 0$. Recall that $\delta \coloneqq \min\{\bar{y}_m, 1 - \bar{y}_m, -\bar{y}_k, 1 + \bar{y}_k\}$. For number $t$ with $0 < t < \delta$, let $u(t) = (u_1, u_2, \ldots) \in l_1$ be defined, for $n$ = 1, 2, … , by

$$u_n = \begin{cases} \bar{x}_n, & \text{if } n \neq m, \\ \bar{x}_m - t, & \text{if } n = m, \\ \bar{x}_k - t, & \text{if } n = k. \end{cases} \tag{3.28}$$

Let $v = (v_1, v_2, \ldots)$ be defined, for $n$ = 1, 2, … , by

$$v_n = \begin{cases} \bar{y}_n, & \text{if } n \neq m, \\ \bar{y}_m - t, & \text{if } n = m, \\ \bar{y}_k - t, & \text{if } n = k. \end{cases} \tag{3.29}$$

By $\bar{y}_m > 0$, $\bar{y}_k < 0$, $\bar{x}_m > 0$, $\bar{x}_k < 0$, $\bar{y}_m \neq 1$ and $\bar{y}_k \neq -1$, with $0 < t < \delta$, we have

$$0 < \bar{y}_m - t < 1 \text{ and } -1 < \bar{y}_k - t < 0.$$

Then, for $n = m$ or $k$, we have

$$0 \leq v_m = \bar{y}_m - t = (\bar{y}_m - t) \wedge 1 = u_m \wedge 1$$

$$\text{and} \quad 0 > v_k = \bar{y}_k - t = (\bar{y}_k - t) \vee (-1) = u_k \vee (-1). \tag{3.30}$$

By (3.29), similarly to the proof of (3.24), we can show that

$$\|v\| = \sum_{n=1}^{\infty} |v_n| = 1. \tag{3.31}$$

By (3.28), similarly to the proof of (3.25), we can show that

$$\|u(t)\| = \sum_{n=1}^{\infty} |u_n| = \|\bar{x}\|. \tag{3.32}$$

By $\bar{y} \in P_{\mathbb{B}}(\bar{x})$, Theorem 3.1, (3.31) and (3.32), similarly to (3.26), we calculate

$$\|u(t) - v\| = \|u(t)\| - 1. \tag{3.33}$$

By (3.30), (3.32) and (3.33), we obtain $v \in P_{\mathbb{B}}(u(t))$, for $0 < t < \delta$. Substituting $x^*$ by $\theta$ in (3.17), we have

$$\begin{aligned} &\limsup_{\substack{(u,v) \to (\bar{x},\bar{y}) \\ v = (v_1, v_2, \ldots) \in P_{\mathbb{B}}(u)}} \frac{\sum_{n=1}^{\infty} x_n^*(u_n - \bar{x}_n) - \sum_{n=1}^{\infty} y_n^*(v_n - \bar{y}_n)}{\|u - \bar{x}\| + \|v - \bar{y}\|_\infty} \\ &= \limsup_{\substack{(u,v) \to (\bar{x},\bar{y}) \\ v = (v_1, v_2, \ldots) \in P_{\mathbb{B}}(u)}} \frac{y_1^* \sum_{n=1}^{\infty} (\bar{y}_n - v_n)}{\|u - \bar{x}\| + \|v - \bar{y}\|_\infty} \\ &= \limsup_{\substack{(u(t),v) \to (\bar{x},\bar{y}), 0 < t < \delta \\ v \in P_{\mathbb{B}}(u), t \downarrow 0}} \frac{y_1^*\big((\bar{y}_m - (\bar{y}_m - t)) + (\bar{y}_k - (\bar{y}_k - t))\big)}{2t + t} \end{aligned}$$

$$= \limsup_{\substack{(u(t),v)\to(\bar{x},\bar{y}),0<t<\delta \\ v\in P_{\mathbb{B}}(u),t\downarrow 0}} \frac{2ty_1^*}{2t+t}$$

$$= \frac{2y_1^*}{3} > 0.$$

This implies that, for $\bar{x} \in l_1$ with $\|\bar{x}\| \geq 1$ and $\bar{y} = (\bar{y}_1, \bar{y}_2, \dots) \in P_{\mathbb{B}}(\bar{x})$ with $-1 < \sum_{n=1}^{\infty} \bar{y}_n < 1$, and $y_1^* = y_2^* = \cdots$

$$\theta \notin \widehat{D}^* P_{\mathbb{B}}(\bar{x}, \bar{y})(y^*), \text{ for } y_1^* > 0. \tag{3.34}$$

Part (D) is proved by (3.27), (3.28) and (3.34).

Proof of (E). We use two ways to prove part (E).

Way 1. Recall that $\mathbb{B}^o$ denotes the topological interior of $\mathbb{B}$ in $l_1$. Then $P_{\mathbb{B}}|_{\mathbb{B}^o} = I_{\mathbb{B}^o}$, which is the identity mapping on $\mathbb{B}^o$. It implies that the Fréchet derivative of $I_{\mathbb{B}^o}$ on $\mathbb{B}^o$ is $I_{\mathbb{B}^o}$. By the connection between Fréchet and Mordukhovich derivatives for single-valued mappings in Banach spaces reviewed in (2.2), for given $\bar{x} = (\bar{x}_1, \bar{x}_2, \dots) \in \mathbb{B}^o$, we obtain that, for any $y^* = (y_1^*, y_2^*) \in \mathbb{R}_M^2$,

$$\widehat{D}^* P_{\mathbb{B}}(\bar{x}, \bar{x})(y^*) = (\nabla P_{\mathbb{B}}(\bar{x}))^*(y^*) = I_{\mathbb{B}^o}^*(y^*) = y^*, \text{ for any } y^* = (y_1^*, y_2^*, \dots) \in l_\infty.$$

Way 2. Here, we provide a directly proof. For given $\bar{x} = (\bar{x}_1, \bar{x}_2, \dots) \in \mathbb{B}^o$ with $\|\bar{x}\| < 1$, there is a positive number $\delta$ such that

$$\|u - \bar{x}\| < \delta \implies \|u\| < 1, \text{ for } u \in l_1.$$

Then, for any $y^* = (y_1^*, y_2^*, \dots) \in l_\infty$, we calculate

$$\limsup_{\substack{(u,v)\to(\bar{x},\bar{x}) \\ v\in P_{\mathbb{B}}(u)}} \frac{\langle y^*, u-\bar{x}\rangle - \langle y^*, v-\bar{x}\rangle}{\|u-\bar{x}\| + \|v-\bar{x}\|_\infty}$$

$$= \limsup_{\substack{(u,v)\to(\bar{x},\bar{x}),\|u-\bar{x}\|<\delta \\ v\in P_{\mathbb{B}}(u)=u}} \frac{\langle y^*, u-\bar{x}\rangle - \langle y^*, u-\bar{x}\rangle}{\|u-\bar{x}\| + \|u-\bar{x}\|_\infty}$$

$$= 0.$$

This implies that

$$y^* \in \widehat{D}^* P_{\mathbb{B}}(\bar{x}, \bar{x})(y^*), \text{ for any } y^* \in l_\infty.$$

Next, let $x^* = (x_1^*, x_2^*, \dots) \in l_\infty$. Suppose that $x^* \neq y^*$. Then, there is a positive number *m* such that $x_m^* \neq y_m^*$. For real number *t* with $|t| < \delta$, we define $u(t) = (u_1, u_2, \dots) \in l_1$ as follows

$$u_n = \begin{cases} \bar{x}_n, & \text{if } n \neq m, \\ \bar{x}_m + t, & \text{if } n = m, \end{cases} \text{ for } n = 1, 2, \dots.$$

By the condition $|t| < \delta$, we have $u(t) = (u_1, u_2, \dots) \in \mathbb{B}^o$. Then, we calculate

$$\limsup_{\substack{(u,v)\to(\bar{x},\bar{x})\\ v=(v_1,v_2,\dots)\in P_{\mathbb{B}}(u)}} \frac{\langle x^*,u-\bar{x}\rangle-\langle y^*,v-\bar{x}\rangle}{\|u-\bar{x}\|+\|v-\bar{x}\|_\infty}$$

$$\geq \limsup_{\substack{(u(t),u(t))\to(\bar{x},\bar{x}),t\to0\\ v(t)\in P_{\mathbb{B}}(u(t))=u(t)}} \frac{\langle x^*,u(t)-\bar{x}\rangle-\langle y^*,u(t)-\bar{x}\rangle}{\|u(t)-\bar{x}\|+\|u(t)-\bar{x}\|_\infty}$$

$$= \limsup_{\substack{(u(t),u(t))\to(\bar{x},\bar{x}),t\to0\\ v(t)\in P_{\mathbb{B}}(u(t))=u(t)}} \frac{\langle x^*-y^*,\ u(t)-\bar{x}\rangle}{\|u(t)-\bar{x}\|+\|u(t)-\bar{x}\|_\infty}$$

$$= \limsup_{t\to0}\frac{(x_m^*-y_m^*)t}{2|t|}$$

$$= \frac{|x_m^*-y_m^*|}{2} > 0.$$

This implies that, for a given $y^* \in l_\infty$,

$$x^* \neq y^* \implies x^* \notin \widehat{D}^*P_{\mathbb{B}}(\bar{x},\bar{x})(y^*), \text{ for } x^* \in l_\infty. \qquad \square$$

**Lemma 3.3**. *Let* $\bar{x} = (\bar{x}_1, \bar{x}_2, \dots) \in l_1$. *If* $\|\bar{x}\| \geq 2$, *Then there is* $\bar{y} = (\bar{y}_1, \bar{y}_2, \dots) \in P_{\mathbb{B}}(\bar{x})$ *such that*

$$\sum_{n=1}^{\infty} \bar{y}_n = \pm 1. \tag{3.35}$$

*Proof*. Let $\bar{x} = (\bar{x}_1, \bar{x}_2, \dots) \in l_1$ with $\|\bar{x}\| \geq 2$. (It is clear that (3.35) holds for $\bar{x} \in (K \cup (-K))$). Let

$$M = \{n \in \mathbb{N}: \bar{x}_n > 0\} \quad \text{and} \quad N = \{n \in \mathbb{N}: \bar{x}_n < 0\}.$$

Then

$$2 \leq \|\bar{x}\| = \sum_{n\in M} \bar{x}_n + \sum_{n\in N}(-\bar{x}_n).$$

This implies that

$$\max\{\sum_{n\in M} \bar{x}_n, \sum_{n\in N}(-\bar{x}_n)\} \geq 1.$$

Case 1. $\sum_{n\in M} \bar{x}_n \geq 1$. In this case, let $\bar{y} = (\bar{y}_1, \bar{y}_2, \dots) \in l_1$ be defined, for $n$ = 1, 2, … by

$$\bar{y}_n = \begin{cases} \frac{\bar{x}_n}{\sum_{i\in M} \bar{x}_i}, & \text{if } n \in M, \\ 0, & \text{if } n \notin M. \end{cases} \tag{3.36}$$

It is clear that

$$\sum_{n=1}^{\infty} |\bar{y}_n| = \sum_{n=1}^{\infty} \bar{y}_n = \sum_{n\in M} \frac{\bar{x}_n}{\sum_{i\in M} \bar{x}_i} = 1. \tag{3.37}$$

By $\sum_{n\in M} \bar{x}_n \geq 1$, we calculate

$$\sum_{n=1}^{\infty} |\bar{y}_n - \bar{x}_n|$$

$$= \sum_{n\in M} |\bar{y}_n - \bar{x}_n| + \sum_{n\in N} |\bar{y}_n - \bar{x}_n| + \sum_{n\notin M\cup N} |\bar{y}_n - \bar{x}_n|$$

$$= \sum_{n\in M}\left|\frac{\bar{x}_n}{\sum_{i\in M}\bar{x}_i} - \bar{x}_n\right| + \sum_{n\in N}|\bar{x}_n| + 0$$

$$= \frac{\sum_{i\in M}\bar{x}_i - 1}{\sum_{i\in M}\bar{x}_i}\sum_{n\in M}|\bar{x}_n| + \sum_{n\in N}|\bar{x}_n|$$

$$= \sum_{n\in M}|\bar{x}_n| - \frac{1}{\sum_{i\in M}\bar{x}_i}\sum_{n\in M}|\bar{x}_n| + \sum_{n\in N}|\bar{x}_n|$$

$$= \sum_{n=1}^{\infty}|\bar{x}_n| - \frac{1}{\sum_{i\in M}\bar{x}_i}\sum_{n\in M}\bar{x}_n$$

$$= \|\bar{x}\| - 1. \tag{3.38}$$

By $\sum_{n\in M}\bar{x}_n \geq 1$, we can check

$$0 \leq \bar{y}_n = \frac{\bar{x}_n}{\sum_{i\in M}\bar{x}_i} \leq \bar{x}_n \wedge 1, \quad \text{for } \bar{x}_n > 0,$$

$$0 = \bar{y}_n \geq \bar{x}_n \vee (-1), \quad \text{for } \bar{x}_n \leq 0. \tag{3.39}$$

By (3.36), (3.37), (3.38) and (3.39), it shows that $\bar{y} = (\bar{y}_1, \bar{y}_2, \ldots) \in P_{\mathbb{B}}(\bar{x})$, which is defined by (3.37). By (3.37), $\bar{y}$ satisfies that $\sum_{n=1}^{\infty}\bar{y}_n = 1$. Then, (3.35) is proved.

Case 2. $\sum_{n\in N}\bar{x}_n \leq -1$. In this case, let $\bar{y} = (\bar{y}_1, \bar{y}_2, \ldots) \in l_1$ be defined, for $n = 1, 2, \ldots$ by

$$\bar{y}_n = \begin{cases} \frac{\bar{x}_n}{\sum_{i\in N}(-\bar{x}_i)}, & \text{if } n \in N, \\ 0, & \text{if } n \notin N. \end{cases}$$

It is clear that

$$\sum_{n=1}^{\infty}|\bar{y}_n| = -\sum_{n=1}^{\infty}\bar{y}_n = \sum_{n\in N}\frac{\bar{x}_n}{\sum_{i\in N}(-\bar{x}_i)} = -1.$$

By $\sum_{n\in N}\bar{x}_n \leq -1$, we calculate

$$\sum_{n=1}^{\infty}|\bar{y}_n - \bar{x}_n|$$

$$= \sum_{n\in M}|\bar{y}_n - \bar{x}_n| + \sum_{n\in N}|\bar{y}_n - \bar{x}_n|$$

$$= +\sum_{n\in M}|\bar{x}_n| + \sum_{n\in N}\left|\frac{\bar{x}_n}{\sum_{i\in N}(-\bar{x}_i)} - \bar{x}_n\right|$$

$$= +\sum_{n\in M}|\bar{x}_n| + \frac{\sum_{i\in N}(-\bar{x}_i) - 1}{\sum_{i\in N}(-\bar{x}_i)}\sum_{n\in N}|\bar{x}_n|$$

$$= \sum_{n\in M}|\bar{x}_n| + \sum_{n\in N}|\bar{x}_n| - \frac{1}{\sum_{i\in N}(-\bar{x}_i)}\sum_{n\in N}|\bar{x}_n|$$

$$= \sum_{n=1}^{\infty}|\bar{x}_n| - \frac{1}{\sum_{i\in N}|\bar{x}_i|}\sum_{n\in N}|\bar{x}_n|$$

$$= \|\bar{x}\| - 1.$$

By $\sum_{n\in N}\bar{x}_n \leq -1$, we can check

$$0 \le \bar{y}_n = 0 \le \bar{x}_n \wedge 1, \qquad \text{for } \bar{x}_n \ge 0,$$

$$0 \ge \bar{y}_n = \frac{\bar{x}_n}{\sum_{n \in N}(-\bar{x}_n)} \ge \bar{x}_n \vee (-1), \qquad \text{for } \bar{x}_n < 0.$$

By Theorem 3.1, they show that $\bar{y} = (\bar{y}_1, \bar{y}_2, \dots) \in P_{\mathbb{B}}(\bar{x})$ and $\bar{y}$ satisfies that $\sum_{n=1}^{\infty} \bar{y}_n = -1$. Then, (3.35) is proved.

**Theorem 3.4**. *Let $\bar{x} = (\bar{x}_1, \bar{x}_2, \dots) \in l_1$ and $\bar{y} = (\bar{y}_1, \bar{y}_2, \dots) \in P_{\mathbb{B}}(\bar{x})$. Then, the covering constant $\hat{\alpha}(P_{\mathbb{B}}, \bar{x}, \bar{y})$ has the following properties*.

(i) *If $\|\bar{x}\| \ge 1$ and $\sum_{n=1}^{\infty} \bar{y}_n = \pm 1$, then*

$$\hat{\alpha}(P_{\mathbb{B}}, \bar{x}, \bar{y}) = 0. \tag{3.40}$$

(ii) *If $\|\bar{x}\| \ge 1$ and $\bar{x} \in K \cup (-K)$, then,* (3.40) *holds with respect to each $\bar{y} \in P_{\mathbb{B}}(\bar{x})$.*

(iii) *If $\|\bar{x}\| \ge 2$, then, there is $\bar{y} = (\bar{y}_1, \bar{y}_2, \dots) \in P_{\mathbb{B}}(\bar{x})$ such that $\sum_{n=1}^{\infty} \bar{y}_n = \pm 1$ and* (3.40) *holds.*

(iv) *For any $\bar{x} \in \mathbb{B}^o$ with $P_{\mathbb{B}}(\bar{x}) = \bar{x}$, we have*

$$\hat{\alpha}(P_{\mathbb{B}}, \bar{x}, \bar{x}) = 1.$$

*Proof*. Proof of (i). Let $\bar{x} = (\bar{x}_1, \bar{x}_2, \dots) \in l_1$ with $\|\bar{x}\| \ge 1$ and $\bar{y} = (\bar{y}_1, \bar{y}_2, \dots) \in P_{\mathbb{B}}(\bar{x})$ with $\sum_{n=1}^{\infty} \bar{y}_n = \pm 1$. Let $\eta > 0$ be arbitrarily given. Let $\mathbb{B}(\bar{x}, \eta)$ and $\mathbb{B}(\bar{y}, \eta)$ be the closed balls in $l_1$ with radius $\eta$ and centered at $\bar{x}$ and $\bar{y}$, respectively. Let $y^* = (y_1^*, y_2^*, \dots) \in l_\infty$ with $y_1^* = y_2^* = \cdots$. By Theorem 3.2, more precisely, if $\sum_{n=1}^{\infty} \bar{y}_n = 1$, we take $y_1^* = -1$; if $\sum_{n=1}^{\infty} \bar{y}_n = -1$, we take $y_1^* = 1$ such that

$$\theta \in \widehat{D}^* P_{\mathbb{B}}(\bar{x}, \bar{y})(y^*).$$

This implies that, for this arbitrarily given $\eta > 0$, we have

$$\begin{aligned}
&\hat{\alpha}(P_{\mathbb{B}}, \bar{x}, \bar{y}) \\
&= \inf\left\{\|z^*\|_\infty : z^* \in \widehat{D}^* P_{\mathbb{B}}(x, y)(w^*), x \in \mathbb{B}(\bar{x}, \eta), y \in P_{\mathbb{B}}(x) \in \mathbb{B}(\bar{y}, \eta), \|w^*\|_\infty = 1\right\} \\
&\le \inf\left\{\|\theta\|_\infty : \theta \in \widehat{D}^* P_{\mathbb{B}}(\bar{x}, \bar{y})(y^*), \|y^*\|_\infty = 1\right\} \\
&= 0.
\end{aligned}$$

Proof of (ii). If $\|\bar{x}\| \ge 1$ and $\bar{x} \in K \cup (-K)$, by Theorem 3.1, for any $\bar{y} = (\bar{y}_1, \bar{y}_2, \dots) \in P_{\mathbb{B}}(\bar{x})$, we have

$$\sum_{n=1}^{\infty} \bar{y}_n = \begin{cases} 1, & \text{if } \bar{x} \in K, \\ -1, & \text{if } \bar{x} \in -K. \end{cases}$$

By part (i), this implies that $\hat{\alpha}(P_{\mathbb{B}}, \bar{x}, \bar{y})$ satisfies (3.40).

Proof of (iii). If $\|\bar{x}\| \ge 2$, by Lemma 3.3, there is $\bar{y} = (\bar{y}_1, \bar{y}_2, \dots) \in P_{\mathbb{B}}(\bar{x})$ such that $\sum_{n=1}^{\infty} \bar{y}_n = \pm 1$, at which, by part (i), (3.40) *holds*.

Proof of (iv). Let $\bar{x} = (\bar{x}_1, \bar{x}_2, \dots) \in \mathbb{B}^o$. There is $\delta > 0$ such that

$$\mathbb{B}(\bar{x},\eta) \subseteq \mathbb{B}^o, \text{ for } \eta < \delta \quad \text{and} \quad \mathbb{B}(\bar{x},\eta)\backslash\mathbb{B}^o \neq \emptyset, \text{ for } \eta \geq \delta.$$

For $\bar{x} = (\bar{x}_1, \bar{x}_2, \dots) \in \mathbb{B}^o$, $P_{\mathbb{B}}(\bar{x}) = \bar{x}$, we have

$$\sup_{0<\eta<\delta} \inf\left\{\|z^*\|_\infty : z^* \in \widehat{D}^* P_{\mathbb{B}}(x,y)(w^*), x \in \mathbb{B}(\bar{x},\eta), y \in P_{\mathbb{B}}(x) \in \mathbb{B}(\bar{y},\eta), \|w^*\|_\infty = 1\right\}$$

$$= \sup_{0<\eta<\delta} \inf\left\{\|z^*\|_\infty : z^* \in \widehat{D}^* P_{\mathbb{B}}(x,y)(w^*), x \in \mathbb{B}(\bar{x},\eta), x = P_{\mathbb{B}}(x) \in \mathbb{B}(\bar{x},\eta), \|w^*\|_\infty = 1\right\}$$

$$= \sup_{0<\eta<\delta} \inf\left\{\|w^*\|_\infty : w^* \in \widehat{D}^* P_{\mathbb{B}}(x,y)(w^*), x \in \mathbb{B}(\bar{x},\eta), x = P_{\mathbb{B}}(x) \in \mathbb{B}(\bar{x},\eta), \|w^*\|_\infty = 1\right\}$$

$$= 1.$$

This proves part (iv). □

## 4. Mordukhovich Differentiability of the Metric Projection Operator in the 2-d Banach Space

### 4.1. The Representation of the Set-Valued Metric Projection Operator in the 2-d Banach Space

In the previous section, we find the explicit solutions of the set-valued metric projection $P_{\mathbb{B}}: l_1 \rightrightarrows \mathbb{B}$. However, since the solutions are complicated, it is extremely difficult to precisely calculate the Mordukhovich derivatives and the covering constants for $P_{\mathbb{B}}$. This is why in Theorems 6.2 and 6.4, we only consider some special cases, which are lack some general details. To put flesh on the details, in this section, we find Mordukhovich derivatives and the covering constants for the set-valued metric projection in the 2-dimensional Banach space.

Let $\mathbb{R}^2$ be the ordinary 2-dimansional real vector space. When $\mathbb{R}^2$ is equipped with the $l_1$-norm, denoted by $\|\cdot\|_l$, then it becomes the 2-dimansional $l_1$-Banach space $\left(\mathbb{R}_l^2, \|\cdot\|_l\right)$, in which the $l_1$-norm $\|\cdot\|_l$ on $\mathbb{R}_l^2$ is defined by

$$\|(x_1, x_2)\|_l = |x_1| + |x_2|, \text{ for any } (x_1, x_2) \in \mathbb{R}_l^2. \tag{4.1}$$

The topological dual space of $\left(\mathbb{R}_l^2, \|\cdot\|_l\right)$ is the Banach space $(\mathbb{R}_M^2, \|\cdot\|_M)$, in which the maximum-norm $\|\cdot\|_M$ on $\mathbb{R}_M^2$ is defined by

$$\|(y_1, y_2)\|_M = \max\{|y_1|, |y_2|\}, \text{ for any } (y_1, y_2) \in \mathbb{R}_M^2. \tag{4.2}$$

As usual, the real pairing between $\mathbb{R}_M^2$ and $\mathbb{R}_l^2$ is written by $\langle\cdot,\cdot\rangle$ and

$$\langle(y_1, y_2), (x_1, x_2)\rangle = x_1 y_1 + x_2 y_2, \text{ for any } (x_1, x_2) \in \mathbb{R}_l^2 \text{ and } (y_1, y_2) \in \mathbb{R}_M^2. \tag{4.3}$$

$\mathbb{R}_l^2$ and $\mathbb{R}_M^2$ have the same set of elements. In particular, their null element is denoted by $\theta = (0, 0)$. Let $\mathbb{B}_l$ and $\mathbb{B}_M$ denote the unit balls in $\mathbb{R}_l^2$ and $\mathbb{R}_M^2$, respectively. Actually, $\mathbb{B}_l$ is the closed square in the plane $\mathbb{R}^2$ with vertexes (1, 0), (0, 1), $(-1, 0)$ and $(0, -1)$. Let $P_{\mathbb{B}_l}: \mathbb{R}_l^2 \rightrightarrows \mathbb{B}_l$ denote the set-valued metric projection operator, which is defined, for any $x = (x_1, x_2) \in \mathbb{R}_l^2$, by

$$P_{\mathbb{B}_l}(x) = \left\{y \in \mathbb{B}_l : \|x - y\|_l = \min_{v \in \mathbb{B}_l} \|x - v\|_l\right\}. \tag{4.4}$$

In particular, $P_{\mathbb{B}_l}\big((x_1, x_2)\big) = \{(x_1, x_2)\}$, for any $(x_1, x_2) \in \mathbb{B}_l$. The square $\mathbb{B}_l$ in the plane $\mathbb{R}^2$ is enclosed by the following four lines:

$$x_2 - x_1 - 1 = 0, \quad x_2 - x_1 + 1 = 0, \quad x_2 + x_1 - 1 = 0 \quad \text{and} \quad x_2 + x_1 + 1 = 0.$$

Let $a$ and $b$ be real numbers, as usual, we denote $a \wedge b = \min\{a, b\}$, $a \vee b = \max\{a, b\}$ and $[a, a] = \{a\}$. Let $\overline{a, b} = \overline{b, a}$ denote the closed segment in $\mathbb{R}^2$ with ending points $a$ and $b$. In particular, $\overline{a, a} = \{a\}$.

In the following proposition, we find the explicit solutions of the set-valued metric projection operator $P_{\mathbb{B}_l}: \mathbb{R}_l^2 \rightrightarrows \mathbb{B}_l$, which can be considered as special cases of the results of Theorem 3.1.

**Proposition 4.1**. *Let* $x = (x_1, x_2) \in \mathbb{R}_l^2$ *with* $\|x\|_l \geq 1$. $P_{\mathbb{B}_l}: \mathbb{R}_l^2 \rightrightarrows \mathbb{B}_l$ *has the following properties.*

$$P_{\mathbb{B}_l}(x) = \left\{ y \in \mathbb{B}_l : \|y\|_l = 1, \|x - y\|_l = \|x\|_l - 1, \begin{matrix} 0 \leq y_n \leq x_n \wedge 1, & \text{if } x_n \geq 0, \\ x_n \vee (-1) \leq y_n \leq 0, & \text{if } x_n < 0, \end{matrix} \; for\ n = 1, 2 \right\}. \quad (4.5)$$

*More precisely,* $P_{\mathbb{B}_l}$ *satisfies the following equations*:

(I) *If* $x_1 \geq 0, x_2 \geq 0$ *and* $x_1 + x_2 > 1$, *then*

$$P_{\mathbb{B}_l}(x) = \overline{(x_1 \wedge 1, 1 - x_1 \wedge 1), (1 - x_2 \wedge 1, x_2 \wedge 1)}.$$

*In particular,*

(a) $x_1 \geq 1$ *and* $x_2 \geq 1 \quad \Longrightarrow \quad P_{\mathbb{B}_l}(x_1, x_2) = \overline{A_l B_l} = \overline{(1, 0), (0, 1)}$;
(b) $0 \leq x_1 \leq 1, 0 \leq x_2 \leq 1$ *and* $x_2 + x_1 \geq 1 \quad \Longrightarrow \quad P_{\mathbb{B}_l}(x_1, x_2) = \overline{(x_1, 1 - x_1), (1 - x_2, x_2)}$;
(c) $x_1 \geq 1$ *and* $0 \leq x_2 \leq 1 \quad \Longrightarrow \quad P_{\mathbb{B}_l}(x_1, x_2) = \overline{(1, 0), \ (1 - x_2, x_2)}$;
(d) $0 \leq x_1 \leq 1$ *and* $x_2 \geq 1 \quad \Longrightarrow \quad P_{\mathbb{B}_l}(x_1, x_2) = \overline{(x_1, 1 - x_1), (0, 1)}$.

(II) *If* $x_1 \leq 0, x_2 \geq 0$ *and* $x_2 - x_1 \geq 1$, *then*

$$P_{\mathbb{B}_l}(x_1, x_2) = \overline{\big(x_1 \vee (-1), \ 1 + x_1 \vee (-1)\big), (x_2 \wedge 1 - 1, \ x_2 \wedge 1)}.$$

(III) *If* $x_1 \leq 0, x_2 \leq 0$ *and* $x_2 + x_1 \leq -1$, *then*

$$P_{\mathbb{B}_l}(x_1, x_2) = \overline{\big(x_1 \vee (-1) + 1, \ x_1 \vee (-1)\big), (x_2 \vee (-1), \ x_2 \vee (-1) + 1)}.$$

(IV) *If* $x_1 \geq 0, x_2 \leq 0$ *and* $x_2 - x_1 \leq -1$, *then*

$$P_{\mathbb{B}_l}(x_1, x_2) = \overline{(x_1 \wedge 1, \ x_1 \wedge 1 - 1), \ \big(1 + x_2 \vee (-1), \ x_2 \vee (-1)\big)}.$$

*In particular,*

(i) $P_{\mathbb{B}_l}(x_1, 0) = (1, 0)$, *for any* $x_1 \geq 1$;
(ii) $P_{\mathbb{B}_l}(0, \ x_2) = (0, 1)$, *for any* $x_2 \geq 1$;
(iii) $P_{\mathbb{B}_l}(x_1, 0) = (-1, 0)$, *for any* $x_1 \leq -1$;
(iv) $P_{\mathbb{B}_l}(0, \ x_2) = (0, -1)$, *for any* $x_2 \leq -1$.

*Proof*. The solutions of $P_{\mathbb{B}_l}: \mathbb{R}_l^2 \rightrightarrows \mathbb{B}_l$ represented by (4.5) follows from Theorem 3.1 immediately. We prove part (I) by using (4.5). Parts (II−IV) can be similarly proved. A directly proof of part (I) of this proposition is provided in the Appendix.

For $x = (x_1, x_2) \in$ Quadrant I with $x_1 \geq 0, x_2 \geq 0$ and $x_1 + x_2 > 1$, by (4.5), we have

$$P_{\mathbb{B}_l}(x) = \{y = (y_1, y_2) \in \mathbb{B}: \|y\|_l = 1, \|x - y\|_l = \|x\|_l - 1 \text{ and } 0 \leq y_i \leq x_i \wedge 1, \text{for } i = 1, 2\}.$$

This implies that

$$y = (y_1, y_2) \in P_{\mathbb{B}_l}(x) \quad \Leftrightarrow \quad y_1 + y_2 = 1 \text{ and } 0 \leq y_i \leq x_i \wedge 1, \text{for } i = 1, 2.$$

We get

$$0 \leq 1 - y_2 \leq x_1 \wedge 1 \ \text{ and } \ 0 \leq 1 - y_1 \leq x_2 \wedge 1.$$

Together, we obtain

$$1 - x_2 \wedge 1 \leq y_1 \leq x_1 \wedge 1 \ \text{ and } \ 1 - x_1 \wedge 1 \leq y_2 \leq x_2 \wedge 1.$$

This proves that

$$P_{\mathbb{B}_l}(x) = \overline{(x_1 \wedge 1, 1 - x_1 \wedge 1), (1 - x_2 \wedge 1, x_2 \wedge 1)}.$$

This completes the proof of (I) of Proposition 4.1. □

### 4.2. The Mordukhovich derivative of the Metric Projection Operator in a 2-d Banach Space

In this subsection, by using the explicit representation of the metric projection operator $P_{\mathbb{B}_l}$ proved in the previous subsection, we investigate the Mordukhovich differentiability of $P_{\mathbb{B}_l}$ with more details than the results of Theorem 3.2.

For *j* = *I*, *II*, *III*, *IV*, let $Q_j$ denote the Quadrat *j* in the plane $\mathbb{R}^2$. In particular, $Q_I$ is the positive cone in $\mathbb{R}_l^2$. Let $\preccurlyeq_2$ denote the partial order on $\mathbb{R}_l^2$ induced by $Q_I$, which is defined, for any $x = (x_1, x_2)$ and $u = (u_1, u_2) \in \mathbb{R}_l^2$, by

$$x \preccurlyeq_2 u \quad \Leftrightarrow \quad x_i \leq u_i \ \text{ and } \ x \prec_2 u \quad \Leftrightarrow \quad x_i < u_i, \text{ for } i = 1, 2.$$

In the following theorem, we mainly investigate the solutions of the Mordukhovich derivative of $P_{\mathbb{B}_l}$ on $Q_I$. The similar details of Mordukhovich derivative of $P_{\mathbb{B}_l}$ on the Quadrats II, III, and IV can be similarly studied.

**Theorem 4.2**. *Let $\bar{x} = (\bar{x}_1, \bar{x}_2) \in \mathbb{R}_l^2$ with $\|\bar{x}\|_l \geq 1$. Let $\bar{y} = (\bar{y}_1, \bar{y}_2) \in P_{\mathbb{B}_l}(\bar{x})$. Then, for any $y^* = (y_1^*, y_2^*) \in \mathbb{R}_M^2$, the Mordukhovich derivative $\widehat{D}^* P_{\mathbb{B}_l}(\bar{x}, \bar{y})(y^*)$ has the following properties.*

(A) *If $\bar{x} \in Q_I$, then,*

$$\widehat{D}^* P_{\mathbb{B}_l}(\bar{x}, \bar{y})(y^*) \subseteq Q_{III} \ (= -Q_I).$$

(B) *If $\bar{x} \in Q_{III}$, then,*

$$\widehat{D}^* P_{\mathbb{B}_l}(\bar{x}, \bar{y})(y^*) \subseteq Q_I \ (= -Q_{III}).$$

(C) *Suppose* $\bar{x} \in Q_I^o \cup Q_{III}^o$. *If* $y_1^* = y_2^*$, *then*

$$\theta \in \widehat{D}^* P_{\mathbb{B}_l}(\bar{x}, \bar{y})(y^*).$$

*More precisely*, $\widehat{D}^* P_{\mathbb{B}_l}(\bar{x}, \bar{y})$ *has the following explicit properties*.

(i) *Suppose* $\bar{x}_1 > 1$ *and* $\bar{x}_2 > 1$.

(a) *If* $(\bar{y}_1, \bar{y}_2) \notin \{(1,0), (0,1)\}$, *then*,

(I) $$y_1^* \neq y_2^* \implies \widehat{D}^* P_{\mathbb{B}_l}(\bar{x}, \bar{y})(y^*) = \emptyset,$$

(II) $$y_1^* = y_2^* \implies \widehat{D}^* P_{\mathbb{B}_l}(\bar{x}, \bar{y})(y^*) = \{\theta\}.$$

(b) *Let* $(\bar{y}_1, \bar{y}_2) = (1, 0)$.

(I) *If* $y_1^* \leq y_2^*$, *then*

$$\widehat{D}^* P_{\mathbb{B}_l}\big(\bar{x}, (1,0)\big)(y^*) = \{\theta\}.$$

(II) *If* $y_1^* > y_2^*$, *then*

$$\widehat{D}^* P_{\mathbb{B}_l}(\bar{x}, (1,0))(y^*) = \emptyset.$$

(c) *Let* $(\bar{y}_1, \bar{y}_2) = (0, 1)$.

(I) *If* $y_1^* \geq y_2^*$, *then*

$$\widehat{D}^* P_{\mathbb{B}_l}\big(\bar{x}, (0,1)\big)(y^*) = \{\theta\}.$$

(II) *If* $y_1^* < y_2^*$, *then*

$$\widehat{D}^* P_{\mathbb{B}_l}(\bar{x}, (0,1))(y^*) = \emptyset.$$

(ii) *Let* $0 < \bar{x}_1 < 1$ *and* $\bar{x}_2 > 1$ *with* $P_{\mathbb{B}_l}(\bar{x}) = \overline{(\bar{x}_1, 1 - \bar{x}_1), (0,1)}$ *and* $\bar{y} \in P_{\mathbb{B}_l}(\bar{x})$.

(a) *Suppose* $0 < \bar{y}_1 < \bar{x}_1$.

(I) *If* $y_1^* \neq y_2^*$, *then*

$$\widehat{D}^* P_{\mathbb{B}_l}(\bar{x}, \bar{y})(y^*) = \emptyset.$$

(II) *If* $y_1^* = y_2^*$, *then*

$$\widehat{D}^* P_{\mathbb{B}_l}(\bar{x}, \bar{y})(y^*) = \{\theta\}.$$

(b) *Let* $(\bar{y}_1, \bar{y}_2) = (0, 1)$.

(I) *If* $y_1^* > y_2^*$, *then*

$$\widehat{D}^* P_{\mathbb{B}_l}(\bar{x}, (0,1))(y^*) = \{\theta\}.$$

(II) *If* $y_1^* < y_2^*$, *then*

$$\widehat{D}^* P_{\mathbb{B}_l}(\bar{x}, (0,1))(y^*) = \emptyset.$$

(c) *Let* $(\bar{y}_1, \bar{y}_2) = (\bar{x}_1, 1 - \bar{x}_1)$.

(I) *If* $y_1^* = y_2^*$, *then*

$$\widehat{D}^* P_{\mathbb{B}_l}\big((\bar{x}_1, \bar{x}_2), (\bar{x}_1, 1 - \bar{x}_1)\big)(y^*) = \{\theta\}.$$

(II) *If* $y_1^* \neq y_2^*$, *then, for any* $x^* = (x_1^*, x_2^*) \in \mathbb{R}_M^2$,

$$\|x^*\|_M < \frac{|y_1^* - y_2^*|}{2} \implies (x_1^*, x_2^*) \notin \widehat{D}^* P_{\mathbb{B}_l}\big((\bar{x}_1, \bar{x}_2), (\bar{x}_1, 1 - \bar{x}_1)\big)(y^*).$$

(iii) *Let* $\bar{x}_1 = 0$ *and* $\bar{x}_2 > 1$ *with* $P_{\mathbb{B}_l}(0, \bar{x}_2) = (0,1)$. *Let* $(\bar{y}_1, \bar{y}_2) = (0,1) = P_{\mathbb{B}_l}(0, \bar{x}_2)$.

(a) *Suppose* $y_1^* = y_2^*$.

(I) *If* $y_1^* = y_2^* < 0$, *then for* $x^* = (x_1^*, 0) \in \mathbb{R}_M^2$,

$$x_1^* \leq 0 \implies (x_1^*, 0) \in \widehat{D}^* P_{\mathbb{B}_l}\big((0, \bar{x}_2), (0,1)\big)(y^*).$$

(II) *If* $y_1^* = y_2^* > 0$, *then for* $x^* = (x_1^*, x_2^*) \in \mathbb{R}_M^2$,

$$\|(x_1^*, x_2^*)\|_M < y_1^* \implies (x_1^*, x_2^*) \notin \widehat{D}^* P_{\mathbb{B}_l}\big((0, \bar{x}_2), (0,1)\big)(y^*).$$

(b) *Suppose* $y_1^* + y_2^* = 0$. *If* $y_1^* > y_2^*$, *then for* $x^* = (x_1^*, 0) \in \mathbb{R}_M^2$,

$$x_1^* \leq y_2^* - y_1^* \implies (x_1^*, 0) \notin \widehat{D}^* P_{\mathbb{B}_l}\big((0, \bar{x}_2), (0,1)\big)(y^*).$$

(iv) *Let* $0 < \bar{x}_1 < 1$ *and* $0 < \bar{x}_2 < 1$ *with* $\bar{x}_1 + \bar{x}_2 > 1$ *and* $P_{\mathbb{B}_l}(\bar{x}) = \overline{(\bar{x}_1, 1 - \bar{x}_1), (1 - \bar{x}_2, \bar{x}_2)}$. *Let* $\bar{y} = (\bar{y}_1, \bar{y}_2) \in P_{\mathbb{B}_l}(\bar{x})$. *For* $y^* = (y_1^*, y_2^*) \in \mathbb{R}_M^2 \backslash \{\theta\}$, *we have*

$$y_1^* \neq y_2^* \implies \widehat{D}^* P_{\mathbb{B}_l}(\bar{x}, \bar{y})(y^*) = \emptyset,$$

$$y_1^* = y_2^* \implies \widehat{D}^* P_{\mathbb{B}_l}(\bar{x}, \bar{y})(y^*) = \{\theta\}.$$

(v) *Let* $\bar{x} = (\bar{x}_1, \bar{x}_2) \in \mathbb{R}_l^2$ *with* $0 \leq \bar{x}_1 \leq 1, 0 \leq \bar{x}_2 \leq 1$ *and* $\bar{x}_1 + \bar{x}_2 = 1$. *Then* $P_{\mathbb{B}_l}(\bar{x}) = \bar{x}$. *Let* $y^* = (y_1^*, y_2^*) \in \mathbb{R}_M^2 \backslash \{\theta\}$.

(a) *If* $y_1^* = y_2^*$, *then,*

(I) $$y_1^* = y_2^* > 0 \implies \theta \notin \widehat{D}^* P_{\mathbb{B}_l}(\bar{x}, \bar{x})(y^*),$$

(II) $$y_1^* = y_2^* < 0 \implies \theta \in \widehat{D}^* P_{\mathbb{B}_l}(\bar{x}, \bar{x})(y^*).$$

(b) *If* $y_1^* \neq y_2^*$, *then,*

$$\theta \notin \widehat{D}^* P_{\mathbb{B}_l}(\bar{x}, \bar{x})(y^*).$$

(vi) *Let* $\bar{x} = (\bar{x}_1, \bar{x}_2) \in \mathbb{R}_l^2$ *with* $0 \leq \bar{x}_1 < 1, 0 \leq \bar{x}_2 < 1$ *and* $\bar{x}_1 + \bar{x}_2 < 1$. *Then* $P_{\mathbb{B}_l}(\bar{x}) = \bar{x}$ *and,*

$$\widehat{D}^* P_{\mathbb{B}_l}(\bar{x}, \bar{x})(y^*) = y^*, \text{ for any } y^* \in \mathbb{R}_M^2.$$

*Proof*. Notice that $Q_I$ is the positive cone in $\mathbb{R}_l^2$ and $-Q_I = Q_{III}$. Then, parts (A) and (B) of this theorem follow from parts (A) and (B) of Theorem 3.2, respectively. We only prove (C). Let $\bar{x} \in Q_I^o \cup Q_{III}^o$ and let $\bar{y} = (\bar{y}_1, \bar{y}_2) \in P_{\mathbb{B}_l}(\bar{x})$. There is $\delta > 0$ such that, for $u = (u_1, u_2) \in \mathbb{R}_l^2$ and $v = (v_1, v_2) \in P_{\mathbb{B}_l}(u)$,

$$\|u - \bar{x}\|_l < \delta \implies u \in Q_I^o \cup Q_{III}^o.$$

Let $v = (v_1, v_2) \in P_{\mathbb{B}_l}(u)$. By Proposition 4.1, this implies that

$$\|u - \bar{x}\|_l < \delta \implies v_1 + v_2 = \bar{y}_1 + \bar{y}_2 = \begin{cases} 1, & \text{if } \bar{x} \in Q_I^o, \\ -1, & \text{if } \bar{x} \in Q_{III}^o. \end{cases}$$

Let $y^* = (y_1^*, y_2^*, \dots) \in \mathbb{R}_M^2$ with $y_1^* = y_2^*$ and let $x^* = \theta \in \mathbb{R}_M^2$. We calculate

$$\limsup_{\substack{(u,v)\to(\bar{x},\bar{y}) \\ v=(v_1,v_2)\in P_{\mathbb{B}_l}(u)}} \frac{\langle x^*, u-\bar{x}\rangle - \langle y^*, v-\bar{y}\rangle}{\|u-\bar{x}\|_l + \|v-\bar{y}\|_M}$$

$$= \limsup_{\substack{(u,v)\to(\bar{x},\bar{y}), \|u-\bar{x}\|_l<\delta \\ v=(v_1,v_2)\in P_{\mathbb{B}_l}(u)}} \frac{-\langle (y_1^*, y_2^*), (v_1,v_2)-(\bar{y}_1,\bar{y}_2)\rangle}{\|u-\bar{x}\|_l + \|v-\bar{y}\|_M}$$

$$= \limsup_{\substack{(u,v)\to(\bar{x},\bar{y}), \|u-\bar{x}\|_l<\delta \\ v=(v_1,v_2)\in P_{\mathbb{B}_l}(u)}} \frac{-y_1^*(v_1+v_2-(\bar{y}_1+\bar{y}_2))}{\|u-\bar{x}\|_l + \|v-\bar{y}\|_M}$$

$$= \limsup_{\substack{(u,v)\to(\bar{x},\bar{y}), \|u-\bar{x}\|_l<\delta \\ v=(v_1,v_2)\in P_{\mathbb{B}_l}(u)}} \frac{0}{\|u-\bar{x}\|_l + \|v-\bar{y}\|_M}$$

$$= 0.$$

This proves that $\theta \in \widehat{D}^* P_{\mathbb{B}_l}(\bar{x}, \bar{y})(y^*)$.

Proof of (I) of (a) of (i). $(\bar{y}_1, \bar{y}_2) \neq (1, 0)$ and $(\bar{y}_1, \bar{y}_2) \neq (0, 1)$. In this case, we have $0 < \bar{y}_1 < 1$ and $0 < \bar{y}_2 < 1$. By $(\bar{x}_1, \bar{x}_2) \in \mathbb{R}_l^2$ with $\bar{x}_1 > 1$ and $\bar{x}_2 > 1$, there is $\delta > 0$ such that, for $u = (u_1, u_2) \in \mathbb{R}_l^2$

$$\|(u_1, u_2) - (\bar{x}_1, \bar{x}_2)\|_l < \delta \implies u_1 > 1 \text{ and } u_2 > 1.$$

By (4.5) or part (I) in Proposition 4.1, this implies that, for $u = (u_1, u_2) \in \mathbb{R}_l^2$,

$$\|u - \bar{x}\|_l < \delta \implies P_{\mathbb{B}_l}(u) = \overline{(1,0), (0,1)}. \tag{4.6}$$

Let $y^* = (y_1^*, y_2^*) \in \mathbb{R}_M^2$. We prove that

$$y_1^* \neq y_2^* \implies \theta \notin \widehat{D}^* P_{\mathbb{B}_l}(\bar{x}, \bar{y})(y^*). \tag{4.7}$$

By (4.6), we calculate

$$\limsup_{\substack{(u,v)\to(\bar{x},\bar{y})\\ u\in\mathbb{R}_l^2 \text{ and } v\in P_{\mathbb{B}_l}(u)}} \frac{\langle\theta,u-\bar{x}\rangle-\langle y^*,v-\bar{y}\rangle}{\|u-\bar{x}\|_l+\|v-\bar{y}\|_M}$$

$$= \limsup_{\substack{(u,v)\to(\bar{x},\bar{y}),\|u-\bar{x}\|_l<\delta\\ v=(v_1,v_2)\in\overline{(1,0),(0,1)}}} \frac{-\langle(y_1^*,y_2^*),\ (v_1,v_2)-(\bar{y}_1,\bar{y}_2)\rangle}{\|(u_1,u_2)-(\bar{x}_1,\bar{x}_2)\|_l+\|(v_1,v_2)-(\bar{y}_1,\bar{y}_2)\|_M}$$

$$= \limsup_{\substack{(u,v)\to(\bar{x},\bar{y}),\|u-\bar{x}\|_l<\delta\\ v=(v_1,v_2)\in\overline{(1,0),(0,1)}}} \frac{-\langle(y_1^*,y_2^*),(v_1,1-v_1)-(\bar{y}_1,1-\bar{y}_1)\rangle}{\|(u_1,u_2)-(\bar{x}_1,\bar{x}_2)\|_l+\|(v_1,1-v_1)-(\bar{y}_1,1-\bar{y}_1)\|_M}$$

$$= \limsup_{\substack{(u,v)\to(\bar{x},\bar{y}),\|u-\bar{x}\|_l<\delta\\ v=(v_1,v_2)\in\overline{(1,0),(0,1)}}} \frac{-\langle(y_1^*,y_2^*),(v_1-\bar{y}_1,\bar{y}_1-v_1)\rangle}{\|(u_1,u_2)-(\bar{x}_1,\bar{x}_2)\|_l+\|(v_1-\bar{y}_1,\bar{y}_1-v_1)\|_M}. \tag{4.8}$$

Let $t$ be a real number. Let $u(t) = (\bar{x}_1+t,\bar{x}_2)$ and $v(t) = (\bar{y}_1+t,1-(\bar{y}_1+t))$. By the conditions that $0<\bar{y}_1<1$ and $0<\bar{y}_2<1$, $\bar{x}_1>1$ and $\bar{x}_2>1$, if $|t|$ is small enough, then $u(t)$ satisfies (4.6) and

$$v(t)\in P_{\mathbb{B}_l}(u(t))=\overline{(1,0),(0,1)}.$$

Substituting $u(t)$ and $v(t)$ into (4.8), we obtain that

$$\limsup_{\substack{(u,v)\to(\bar{x},\bar{y})\\ u\in\mathbb{R}_l^2 \text{ and } v\in P_{\mathbb{B}_l}(u)}} \frac{\langle\theta,u-\bar{x}\rangle-\langle y^*,v-\bar{y}\rangle}{\|u-\bar{x}\|_l+\|v-\bar{y}\|_M}$$

$$\geq \limsup_{\substack{(u(t),v(t))\to(\bar{x},\bar{y}),\|u(t)-\bar{x}\|_l<\delta\\ v(t)\in P_{\mathbb{B}_l}(u(t)),t\to 0}} \frac{-\langle(y_1^*,y_2^*),(t,-t)\rangle}{\|(t,0)\|_l+\|(t,-t)\|_M}$$

$$= \limsup_{\substack{(u(t),v(t))\to(\bar{x},\bar{y}),\|u(t)-\bar{x}\|_l<\delta\\ v(t)\in P_{\mathbb{B}_l}(u(t)),t\to 0}} \frac{-t(y_1^*-y_2^*)}{|t|+|t|}$$

$$= \frac{|y_1^*-y_2^*|}{2} > 0.$$

This proves that $\theta\notin\widehat{D}^*P_{\mathbb{B}_l}(\bar{x},\bar{y})(y^*)$, which proves (4.7).

Let $y^*=(y_1^*,y_2^*)\in\mathbb{R}_M^2$ and $x^*=(x_1^*,x_2^*)\in\mathbb{R}_M^2$. For $v=(v_1,v_2)\in P_{\mathbb{B}_l}(u)$, $\bar{y}=(\bar{y}_1,\bar{y}_2)\in P_{\mathbb{B}_l}(\bar{x})$, we have $v_2=1-v_1$ and $\bar{y}_2=1-\bar{y}_1$.

$$\limsup_{\substack{(u,v)\to(\bar{x},\bar{y})\\ u\in\mathbb{R}_l^2 \text{ and } v\in P_{\mathbb{B}_l}(u)}} \frac{\langle x^*,u-\bar{x}\rangle-\langle y^*,v-\bar{y}\rangle}{\|u-\bar{x}\|_l+\|v-\bar{y}\|_M}$$

$$= \limsup_{\substack{(u,v)\to(\bar{x},\bar{y})\\ v=(v_1,v_2)\in P_{\mathbb{B}_l}(u)}} \frac{x_1^*(u_1-\bar{x}_1)+x_2^*(u_2-\bar{x}_2)-\big(y_1^*(v_1-\bar{y}_1)+y_2^*(v_2-\bar{y}_2)\big)}{\|u-\bar{x}\|_l+\|v-\bar{y}\|_M}$$

$$= \limsup_{\substack{(u,v)\to(\bar{x},\bar{y}) \\ v=(v_1,v_2)\in P_{\mathbb{B}_l}(u)}} \frac{x_1^*(u_1-\bar{x}_1)+x_2^*(u_2-\bar{x}_2)-\left(y_1^*(v_1-\bar{y}_1)+y_2^*\left((1-v_1)-(1-\bar{y}_1)\right)\right)}{\|u-\bar{x}\|_l+\|v-\bar{y}\|_M}$$

$$= \limsup_{\substack{(u,v)\to(\bar{x},\bar{y}) \\ v=(v_1,v_2)\in P_{\mathbb{B}_l}(u)}} \frac{x_1^*(u_1-\bar{x}_1)+x_2^*(u_2-\bar{x}_2)+(y_2^*-y_1^*)(v_1-\bar{y}_1)}{\|u-\bar{x}\|_l+\|v-\bar{y}\|_M}. \quad (4.9)$$

Let $t$ be a real number. Let $u(t) = (\bar{x}_1+t, \bar{x}_2+t)$ and $v = \bar{y} = (\bar{y}_1, \bar{y}_2)$. By the conditions that $0 < \bar{y}_1 < 1$ and $0 < \bar{y}_2 < 1$, if $|t|$ is small enough, then $u(t)$ satisfies (4.6) and

$$v = \bar{y} \in P_{\mathbb{B}_l}(u(t)) = \overline{(1,0),(0,1)}.$$

Then, substituting $u(t)$ and $v = \bar{y}$ into (4.9), we obtain that

$$\limsup_{\substack{(u,v)\to(\bar{x},\bar{y}) \\ u\in\mathbb{R}_l^2 \text{ and } v\in P_{\mathbb{B}_l}(u)}} \frac{\langle x^*, u-\bar{x}\rangle - \langle y^*, v-\bar{y}\rangle}{\|u-\bar{x}\|_l+\|v-\bar{y}\|_M}$$

$$\geq \limsup_{\substack{(u(t),v)\to(\bar{x},\bar{y}), \|u(t)-\bar{x}\|_l<\delta \\ v=\bar{y}, t\downarrow 0}} \frac{t(x_1^*+x_2^*) - \langle (y_1^*,y_2^*),(0,0)\rangle}{\|(t,t)\|_l+\|(0,0)\|_M}$$

$$= \limsup_{\substack{(u,v)\to(\bar{x},\bar{y}), \|u-\bar{x}\|_l<\delta \\ v=\bar{y}, t\downarrow 0}} \frac{t(x_1^*+x_2^*)}{2|t|}$$

$$= \frac{|x_1^*+x_2^*|}{2}. \quad (4.10)$$

On another hand, let $t$ be a real number. Let $w(t) = (\bar{x}_1+t, \bar{x}_2-t)$ and $v = \bar{y} = (\bar{y}_1, \bar{y}_2)$. If $|t|$ is small enough, then $w(t)$ satisfies (4.6) and

$$v = \bar{y} = (\bar{y}_1, \bar{y}_2) \in P_{\mathbb{B}_l}(u(t)) = \overline{(1,0),(0,1)}.$$

Then, substituting $w(t)$ and $v$ into (4.9), we obtain that

$$\limsup_{\substack{(u,v)\to(\bar{x},\bar{y}) \\ u\in\mathbb{R}_l^2 \text{ and } v\in P_{\mathbb{B}_l}(u)}} \frac{\langle x^*, u-\bar{x}\rangle - \langle y^*, v-\bar{y}\rangle}{\|u-\bar{x}\|_l+\|v-\bar{y}\|_M}$$

$$\geq \limsup_{\substack{(u(t),v)\to(\bar{x},\bar{y}), \|u(t)-\bar{x}\|_l<\delta \\ v=\bar{y}\in\overline{(1,0),(0,1)}}} \frac{\langle (x_1^*,x_2^*),(t,-t)\rangle - \langle (y_1^*,y_2^*),(0,0)\rangle}{\|(t,-t)\|_l+\|(0,0)\|_M}$$

$$= \limsup_{\substack{(u(t),v)\to(\bar{x},\bar{y}), \|u(t)-\bar{x}\|_l<\delta \\ (v_1,v_2)=\bar{y}\in\overline{(1,0),(0,1)}}} \frac{t(x_1^*-x_2^*)}{2|t|}$$

$$= \frac{|x_1^*-x_2^*|}{2}. \quad (4.11)$$

By the condition that $(x_1^*, x_2^*) \neq (0,0)$, from (4.10) and (4.11), we obtain that

$$\limsup_{\substack{(u,v)\to(\bar{x},\bar{y}) \text{ in } \mathbb{R}_l^2 \\ u\in\mathbb{R}_l^2 \text{ and } v\in P_{\mathbb{B}_l}(u)}} \frac{\langle x^*, u-\bar{x}\rangle - \langle y^*, v-\bar{y}\rangle}{\|u-\bar{x}\|_l + \|v-\bar{y}\|_M}$$

$$\geq \max\left\{\frac{|x_1^*+x_2^*|}{2}, \frac{|x_1^*-x_2^*|}{2}\right\} > 0.$$

We obtain that

$$x^* \notin \widehat{D}^* P_{\mathbb{B}_l}(\bar{x},\bar{y})(y^*), \text{ for any } x^* \neq \theta. \tag{4.12}$$

This proves the part (I) of (a) in (i).

Proof of (II) of (a) of (i). Let $y^* = (y_1^*, y_2^*) \in \mathbb{R}_M^2 \backslash \{\theta\}$ with $y_1^* = y_2^*$. By part (B) and (4.12), we have

$$\widehat{D}^* P_{\mathbb{B}_l}(\bar{x},\bar{y})(y^*) = \{\theta\}.$$

Proof (I) of (b) of (i). Let $(\bar{y}_1, \bar{y}_2) = (1, 0)$. Let $y^* = (y_1^*, y_2^*) \in \mathbb{R}_M^2 \backslash \{\theta\}$. By (4.6), for $u \in \mathbb{R}_M^2$ , if $\|u-\bar{x}\|_l < \delta$ and $v = (v_1, 1-v_1) \in P_{\mathbb{B}_l}(u) = \overline{(1,0),(0,1)}$, then $0 \leq v_1 \leq 1$. Let $x^* = \theta$ in (4.8).

$$\limsup_{\substack{(u,v)\to(\bar{x},(1,0)) \\ v\in P_{\mathbb{B}_l}(u)}} \frac{\langle \theta, u-\bar{x}\rangle - \langle y^*, v-\bar{y}\rangle}{\|u-\bar{x}\|_l + \|v-\bar{y}\|_M}$$

$$= \limsup_{\substack{(u,v)\to(\bar{x},(1,0)), \|u-\bar{x}\|_l<\delta \\ v=(v_1,v_2)\in P_{\mathbb{B}_l}(u)}} \frac{(y_2^*-y_1^*)(v_1-1)}{\|u-\bar{x}\|_l + \|v-\bar{y}\|_M}$$

$$= \limsup_{\substack{(u,v)\to(\bar{x},(1,0)) \\ (v_1,1-v_1)\in P_{\mathbb{B}_l}(u)}} \frac{-(y_1^*-y_2^*)(v_1-1)}{\|(u_1,u_2)-(\bar{x}_1,\bar{x}_2)\|_l + \|(v_1-1, 1-v_1)\|_M}$$

$$= \limsup_{\substack{(u,v)\to(\bar{x},(1,0)), \|u-\bar{x}\|_l<\delta \\ (v_1,1-v_1)\to(1,0),\ v_1\leq 1}} \frac{(y_1^*-y_2^*)(1-v_1)}{\|(u_1,u_2)-(\bar{x}_1,\bar{x}_2)\|_l + 1 - v_1}$$

$$= \limsup_{\substack{(u,v)\to(\bar{x},(1,0)), v_1\uparrow 1 \\ \|u-\bar{x}\|_l = 1-v_1<\delta}} \frac{(y_1^*-y_2^*)(1-v_1)}{1-v_1+1-v_1}$$

$$= \frac{y_1^*-y_2^*}{2}. \tag{4.13}$$

This implies that

$$y_1^* \leq y_2^* \implies \theta \in \widehat{D}^* P_{\mathbb{B}_l}\big((\bar{x}_1,\bar{x}_2),(1,0)\big)(y_1^*, y_2^*) \tag{4.14}$$

Let $x^* = (x_1^*, x_2^*) \in \mathbb{R}_M^2 \backslash \{\theta\}$. Without loose of the generality, suppose $x_1^* \neq 0$. Similarly to the above proof, we calculate

$$\limsup_{\substack{(u,v)\to(\bar{x},(1,0)) \\ u\in\mathbb{R}_l^2 \text{ and } v\in P_{\mathbb{B}_l}(u)}} \frac{\langle x^*, u-\bar{x}\rangle - \langle y^*, v-\bar{y}\rangle}{\|u-\bar{x}\|_l + \|v-\bar{y}\|_M}$$

$$= \limsup_{\substack{(u,v)\to(\bar{x},(1,0)),\|u-\bar{x}\|_l<\delta \\ v=(v_1,v_2)\in\overline{(1,0),(0,1)}}} \frac{x_1^*(u_1-\bar{x}_1)+x_2^*(u_2-\bar{x}_2)+(y_2^*-y_1^*)(v_1-1)}{\|u-\bar{x}\|_l+\|v-\bar{y}\|_M}$$

$$\geq \limsup_{\substack{u=(\bar{x}_1+t,\bar{x}_2),|t|<\delta \\ v=(1,0)\in\overline{(1,0),(0,1)}}} \frac{x_1^*(u_1-\bar{x}_1)+x_2^*(u_2-\bar{x}_2)+(y_2^*-y_1^*)(v_1-1)}{\|u-\bar{x}\|_l+\|v-\bar{y}\|_M}$$

$$\geq \limsup_{\substack{u=(\bar{x}_1+t,\bar{x}_2),|t|<\delta \\ v=(1,0)\in\overline{(1,0),(0,1)}}} \frac{x_1^*t+0+0}{|t|+0}$$

$$= |x_1^*|$$

$$> 0. \tag{4.15}$$

(4.15) implies that, for any $x^* = (x_1^*, x_2^*) \in \mathbb{R}_M^2\backslash\{\theta\}$

$$x^* \notin \widehat{D}^* P_{\mathbb{B}_l}\big((\bar{x}_1, \bar{x}_2), (1,0)\big)(y^*). \tag{4.16}$$

By (4.14) and (4.16), we have that

$$\widehat{D}^* P_{\mathbb{B}_l}\big((\bar{x}_1, \bar{x}_2), (1,0)\big)(y^*) = \{\theta\}.$$

Proof (II) of (b) of (i). From (4.13), we immediately obtain that

$$y_1^* > y_2^* \implies \theta \notin \widehat{D}^* P_{\mathbb{B}_l}\big((\bar{x}_1, \bar{x}_2), (1,0)\big)(y^*). \tag{4.17}$$

Let $x^* = (x_1^*, x_2^*) \in \mathbb{R}_M^2\backslash\{(0,0)\}$. Suppose $y_1^* > y_2^*$. Similarly to the above proof, we calculate

$$\limsup_{\substack{(u,v)\to(\bar{x},(1,0)) \\ u\in\mathbb{R}_l^2 \text{ and } v\in P_{\mathbb{B}_l}(u)}} \frac{\langle x^*, u-\bar{x}\rangle - \langle y^*, v-\bar{y}\rangle}{\|u-\bar{x}\|_l+\|v-\bar{y}\|_M}$$

$$= \limsup_{\substack{(u,v)\to(\bar{x},(1,0)),\|u-\bar{x}\|_l<\delta \\ (v_1,v_2)\in\overline{(1,0),(0,1)}}} \frac{\langle x^*, u-\bar{x}\rangle - \langle (y_1^*,y_2^*),(v_1,v_2)-(1,0)\rangle}{\|(u_1,u_2)-(\bar{x}_1,\bar{x}_2)\|_l+\|(v_1,v_2)-(1,0)\|_M}$$

$$= \limsup_{\substack{(u,v)\to(\bar{x},(1,0)) \\ (v_1,v_2)\in\overline{(1,0),(0,1)}}} \frac{\langle x^*, u-\bar{x}\rangle - \langle (y_1^*,y_2^*),(v_1,1-v_1)-(1,0)\rangle}{\|(u_1,u_2)-(\bar{x}_1,\bar{x}_2)\|_l+\|(v_1,1-v_1)-(1,0)\|_M}$$

$$= \limsup_{\substack{(u,v)\to(\bar{x},(1,0)) \\ (v_1,v_2)\in\overline{(1,0),(0,1)}}} \frac{\langle x^*, u-\bar{x}\rangle - \langle (y_1^*,y_2^*),(v_1-1,1-v_1)\rangle}{\|(u_1,u_2)-(\bar{x}_1,\bar{x}_2)\|_l+\|(v_1,1-v_1)-(1,0)\|_M}$$

$$= \limsup_{\substack{(u,v)\to(\bar{x},(1,0)) \\ (v_1,1-v_1)\in\overline{(1,0),(0,1)}}} \frac{\langle x^*, u-\bar{x}\rangle - (y_1^*-y_2^*)(v_1-1)}{\|(u_1,u_2)-(\bar{x}_1,\bar{x}_2)\|_l+\|(v_1-1,1-v_1)\|_M}$$

$$= \limsup_{\substack{(u,v)\to(\bar{x},(1,0)),\|u-\bar{x}\|_l<\frac{y_1^*-y_2^*}{2\|x^*\|_M}(1-v_1)\\ (v_1,1-v_1)\in\overline{(1,0),(0,1)},\ v_1\uparrow 1}} \frac{\langle x^*,u-\bar{x}\rangle+(y_1^*-y_2^*)(1-v_1)}{\|(u_1,u_2)-(\bar{x}_1,\bar{x}_2)\|_l+1-v_1}$$

$$\geq \limsup_{\substack{(u,v)\to(\bar{x},(1,0)),\|u-\bar{x}\|_l<\frac{y_1^*-y_2^*}{2\|x^*\|_M}(1-v_1)\\ (v_1,1-v_1)\in\overline{(1,0),(0,1)},\ v_1\uparrow 1}} \frac{-\|x^*\|_M\frac{y_1^*-y_2^*}{2\|x^*\|_M}(1-v_1)+(y_1^*-y_2^*)(1-v_1)}{\|(u_1,u_2)-(\bar{x}_1,\bar{x}_2)\|_l+1-v_1}$$

$$= \limsup_{\substack{(u,v)\to(\bar{x},(1,0)),\|u-\bar{x}\|_l<\frac{y_1^*-y_2^*}{2\|x^*\|_M}(1-v_1)\\ (v_1,1-v_1)\in\overline{(1,0),(0,1)},\ v_1\uparrow 1}} \frac{-\frac{y_2^*-y_1^*}{2}(1-v_1)+(y_1^*-y_2^*)(1-v_1)}{\|(u_1,u_2)-(\bar{x}_1,\bar{x}_2)\|_l+1-v_1}$$

$$\geq \limsup_{\substack{(u,v)\to(\bar{x},(1,0)),\|u-\bar{x}\|_l<\frac{y_1^*-y_2^*}{2\|x^*\|_M}(1-v_1)\\ (v_1,1-v_1)\in\overline{(1,0),(0,1)},\ v_1\uparrow 1}} \frac{\frac{1}{2}(y_1^*-y_2^*)(1-v_1)}{\frac{y_1^*-y_2^*}{2\|x^*\|_M}(1-v_1)+1-v_1}$$

$$= \frac{\frac{1}{2}(y_1^*-y_2^*)}{\frac{y_1^*-y_2^*}{2\|x^*\|_M}+1} > 0.$$

This implies that for $y^*=(y_1^*,y_2^*)\in\mathbb{R}_M^2\backslash\{\theta\}$, if $y_1^*>y_2^*$, then for any $x^*=(x_1^*,x_2^*)\in\mathbb{R}_M^2$,

$$x^*\neq\theta \Longrightarrow x^*\notin\widehat{D}^*P_{\mathbb{B}_l}\big((\bar{x}_1,\bar{x}_2),(1,0)\big)(y^*). \tag{4.18}$$

Then, (II) of (b) of (i) is proved by (4.17) and (4.18).

Proof of (I) of (c) of (i). $(\bar{y}_1,\bar{y}_2)=(0,1)$. Let $y^*=(y_1^*,y_2^*)\in\mathbb{R}_M^2\backslash\{\theta\}$. We have

$$\limsup_{\substack{(u,v)\to(\bar{x},(0,1))\\ u\in\mathbb{R}_l^2 \text{ and } v\in P_{\mathbb{B}_l}(u)}} \frac{\langle\theta,u-\bar{x}\rangle-\langle y^*,v-\bar{y}\rangle}{\|u-\bar{x}\|_l+\|v-\bar{y}\|_M}$$

$$= \limsup_{\substack{(u,v)\to(\bar{x},(0,1)),\|u-\bar{x}\|_l<\delta\\ (v_1,v_2)\in\overline{(1,0),(0,1)}}} \frac{-\langle(y_1^*,y_2^*),(v_1,v_2)-(0,1)\rangle}{\|(u_1,u_2)-(\bar{x}_1,\bar{x}_2)\|_l+\|(v_1,v_2)-(0,1)\|_M}$$

$$= \limsup_{\substack{(u,v)\to(\bar{x},(1,0))\\ (v_1,v_2)\in\overline{(1,0),(0,1)}}} \frac{-\langle(y_1^*,y_2^*),(v_1,1-v_1)-(0,1)\rangle}{\|(u_1,u_2)-(\bar{x}_1,\bar{x}_2)\|_l+\|(v_1,1-v_1)-(0,1)\|_M}$$

$$= \limsup_{\substack{(u,v)\to(\bar{x},(0,1))\\ (v_1,v_2)\in\overline{(1,0),(0,1)}}} \frac{-\langle(y_1^*,y_2^*),(v_1,-v_1)\rangle}{\|(u_1,u_2)-(\bar{x}_1,\bar{x}_2)\|_l+\|(v_1,1-v_1)-(0,1)\|_M}$$

$$= \limsup_{\substack{(u,(v_1,1-v_1))\to(\bar{x},(0,1))\\ (v_1,1-v_1)\in\overline{(1,0),(0,1)}}} \frac{-(y_1^*-y_2^*)v_1}{\|(u_1,u_2)-(\bar{x}_1,\bar{x}_2)\|_l+\|(v_1,-v_1)\|_M}$$

$$= \limsup_{\substack{(u,(v_1,1-v_1))\to(\bar{x},(0,1)),\|u-\bar{x}\|_l<\delta \\ (v_1,1-v_1)\in\overline{(1,0),(0,1)},\ v_1\downarrow 0}} \frac{-(y_1^*-y_2^*)v_1}{\|(u_1,u_2)-(\bar{x}_1,\bar{x}_2)\|_l+v_1}$$

$$= -\frac{y_1^*-y_2^*}{2}. \tag{4.19}$$

This implies that, for $y^* = (y_1^*, y_2^*) \in \mathbb{R}_M^2\backslash\{\theta\}$,

$$y_1^* \geq y_2^* \implies \theta \in \widehat{D}^*P_{\mathbb{B}_l}\big((\bar{x}_1,\bar{x}_2),(0,1)\big)(y^*). \tag{4.20}$$

Similarly to the proof of (4.15) of part (I) in (b), for any $x^* = (x_1^*, x_2^*) \in \mathbb{R}_M^2\backslash\{\theta\}$, we can prove

$$x^* \notin \widehat{D}^*P_{\mathbb{B}_l}\big((\bar{x}_1,\bar{x}_2),(0,1)\big)(y^*), \text{ for any } x^* \in \mathbb{R}_M^2\backslash\{\theta\}. \tag{4.21}$$

By (4.20) and (4.21), we have

$$\widehat{D}^*P_{\mathbb{B}_l}\big((\bar{x}_1,\bar{x}_2),(0,1)\big)(y^*) = \{\theta\}.$$

Proof of (II) of (c) of (i). $(\bar{y}_1,\bar{y}_2) = (0,1)$. By the calculation of (4.19), we have

$$y_1^* < y_2^* \implies \theta \notin \widehat{D}^*P_{\mathbb{B}_l}\big((\bar{x}_1,\bar{x}_2),(0,1)\big)(y^*). \tag{4.22}$$

Similarly to the proof of (4.18), of part (I) in (c), we can prove

$$x^* \notin \widehat{D}^*P_{\mathbb{B}_l}\big((\bar{x}_1,\bar{x}_2),(0,1)\big)(y^*), \text{ for any } x^* \in \mathbb{R}_M^2\backslash\{\theta\}. \tag{4.23}$$

By (4.22) and (4.23), we have that if $y_1^* < y_2^*$, then

$$\widehat{D}^*P_{\mathbb{B}_l}\big((\bar{x}_1,\bar{x}_2),(0,1)\big)(y^*) = \emptyset.$$

(ii) Let $0 < \bar{x}_1 < 1$ and $\bar{x}_2 > 1$. Let $\bar{y} = (\bar{y}_1,\bar{y}_2) \in P_{\mathbb{B}_l}(\bar{x}_1,\bar{x}_2) = \overline{(\bar{x}_1, 1-\bar{x}_1),(0,1)}$.

Proof (I) of (a) in (ii). Let $0 < \bar{y}_1 < \bar{x}_1$. By $\bar{y}_1 + \bar{y}_2 = 1$, this implies that $1 - \bar{x}_1 < \bar{y}_2 < 1$. There is $\lambda > 0$ such that, for $u = (u_1,u_2) \in \mathbb{R}_l^2$

$$\|(u_1,u_2)-(\bar{x}_1,\bar{x}_2)\|_l < \lambda \implies \bar{y}_1 < u_1 < 1 \text{ and } u_2 > 1.$$

This implies that

$$\|(u_1,u_2)-(\bar{x}_1,\bar{x}_2)\|_l < \lambda \implies \bar{y}_1 < u_1 < 1, u_2 > 1 \text{ and } 1-u_1 < 1-\bar{y}_1 = \bar{y}_2.$$

Then, we obtain that, for $u = (u_1,u_2) \in \mathbb{R}_l^2$,

$$\|(u_1,u_2)-(\bar{x}_1,\bar{x}_2)\|_l < \lambda \implies (\bar{y}_1,\bar{y}_2) \in P_{\mathbb{B}_l}(u_1,u_2) = \overline{(u_1, 1-u_1),(0,1)}. \tag{4.24}$$

Let $y^* = (y_1^*, y_2^*) \in \mathbb{R}_M^2\backslash\{(0,0)\}$. Suppose $y_1^* \neq y_2^*$. For any $x^* = (x_1^*, x_2^*) \in \mathbb{R}_M^2$, by (4.24), we calculate

$$\limsup_{\substack{(u,v)\to(\bar{x},\bar{y}) \\ u\in\mathbb{R}_l^2 \text{ and } v\in P_{\mathbb{B}_l}(u)}} \frac{\langle x^*, u-\bar{x}\rangle - \langle y^*, v-\bar{y}\rangle}{\|u-\bar{x}\|_l + \|v-\bar{y}\|_M}$$

$$= \limsup_{\substack{(u,v)\to((\bar{x}_1,\bar{x}_2),(\bar{y}_1,\bar{y}_2)),\|u-\bar{x}\|_l<\lambda \\ v=(v_1,v_2)\in\overline{(u_1,1-u_1),(0,1)}}} \frac{\langle(x_1^*,x_2^*),(u_1,u_2)-(\bar{x}_1,\bar{x}_2)\rangle - \langle(y_1^*,y_2^*),(v_1,v_2)-(\bar{y}_1,\bar{y}_2)\rangle}{\|(u_1,u_2)-(\bar{x}_1,\bar{x}_2)\|_l + \|(v_1,v_2)-(\bar{y}_1,\bar{y}_2)\|_M}$$

$$= \limsup_{\substack{(u,v)\to((\bar{x}_1,\bar{x}_2),(\bar{y}_1,1-\bar{y}_1)),\|u-\bar{x}\|_l<\lambda \\ v=(v_1,v_2)\in\overline{(u_1,1-u_1),(0,1)}}} \frac{\langle(x_1^*,x_2^*),(u_1,u_2)-(\bar{x}_1,\bar{x}_2)\rangle - \langle(y_1^*,y_2^*),(v_1,1-v_1)-(\bar{y}_1,1-\bar{y}_1)\rangle}{\|(u_1,u_2)-(\bar{x}_1,\bar{x}_2)\|_l + \|(v_1,1-v_1)-(\bar{y}_1,1-\bar{y}_1)\|_M}$$

$$= \limsup_{\substack{(u,v)\to(\bar{x},(\bar{y}_1,1-\bar{y}_1)),\|u-\bar{x}\|_l<\lambda \\ v=(v_1,v_2)\in\overline{(u_1,1-u_1),(0,1)}}} \frac{x_1^*(u_1-\bar{x}_1) + x_2^*(u_2-\bar{x}_2)-(y_1^*-y_2^*)(v_1-\bar{y}_1)}{\|(u_1,u_2)-(\bar{x}_1,\bar{x}_2)\|_l + \|(v_1-\bar{y}_1,-(v_1-\bar{y}_1))\|_M} \,. \tag{4.25}$$

Case 1. Let $x^* = \theta$ in (4.25). We have

$$\limsup_{\substack{(u,v)\to(\bar{x},\bar{y}) \\ u\in\mathbb{R}_l^2 \text{ and } v\in P_{\mathbb{B}_l}(u)}} \frac{\langle\theta,u-\bar{x}\rangle - \langle y^*,v-\bar{y}\rangle}{\|u-\bar{x}\|_l + \|v-\bar{y}\|_M}$$

$$= \limsup_{\substack{(u,v)\to(\bar{x},(\bar{y}_1,1-\bar{y}_1)),\|u-\bar{x}\|_l<\lambda \\ (v_1,v_2)\in\overline{(u_1,1-u_1),(0,1)}}} \frac{-(y_1^*-y_2^*)(v_1-\bar{y}_1)}{\|(u_1,u_2)-(\bar{x}_1,\bar{x}_2)\|_l + \|(v_1-\bar{y}_1,-(v_1-\bar{y}_1))\|_M}. \tag{4.26}$$

Let $t$ be a real number satisfying $|t| < \min\{\lambda, \bar{x}_1 - \bar{y}_1\}$ (Recall that in this part (a), $0 < \bar{y}_1 < \bar{x}_1$ and $1 - \bar{x}_1 < 1 - \bar{y}_1 < 1$). We write $u(t) = (1+t)\bar{x}$ and $v_1(t) = \bar{y}_1 + t$. Substituting these to (4.26), we obtain

$$\limsup_{\substack{(u,v)\to(\bar{x},\bar{y}) \\ u\in\mathbb{R}_l^2 \text{ and } v\in P_{\mathbb{B}_l}(u)}} \frac{\langle\theta,u-\bar{x}\rangle - \langle y^*,v-\bar{y}\rangle}{\|u-\bar{x}\|_l + \|v-\bar{y}\|_M}$$

$$\geq \limsup_{t\to 0} \frac{-(y_1^*-y_2^*)t}{|t|\|\bar{x}\|_l + |t|}$$

$$= \limsup_{t\to 0} \frac{(y_2^*-y_1^*)t}{|t|\|\bar{x}\|_l + |t|}$$

$$= \frac{|y_1^*-y_2^*|}{\|\bar{x}\|_l + 1} > 0, \text{ for } y_1^* \neq y_2^*.$$

For $\bar{y} = (\bar{y}_1, \bar{y}_2)$ with $0 < \bar{y}_1 < \bar{x}_1$, this implies that

$$y_1^* \neq y_2^* \implies \theta \notin \widehat{D}^* P_{\mathbb{B}_l}(\bar{x}, \bar{y})(y^*). \tag{4.27}$$

Case 2. If $x^* = (x_1^*, x_2^*) \in \mathbb{R}_M^2 \backslash \{\theta\}$. By (4.25), we have

$$\limsup_{\substack{(u,v)\to(\bar{x},\bar{y}) \\ u\in\mathbb{R}_l^2 \text{ and } v\in P_{\mathbb{B}_l}(u)}} \frac{\langle x^*,u-\bar{x}\rangle - \langle y^*,v-\bar{y}\rangle}{\|u-\bar{x}\|_l + \|v-\bar{y}\|_M}$$

$$= \limsup_{\substack{(u,v)\to(\bar{x},(\bar{y}_1,1-\bar{y}_1)),\|u-\bar{x}\|_l<\lambda \\ v=(v_1,v_2)\in\overline{(u_1,1-u_1),(0,1)}}} \frac{x_1^*(u_1-\bar{x}_1) + x_2^*(u_2-\bar{x}_2)-(y_1^*-y_2^*)(v_1-\bar{y}_1)}{\|(u_1,u_2)-(\bar{x}_1,\bar{x}_2)\|_l + \|(v_1-\bar{y}_1,-(v_1-\bar{y}_1))\|_M}$$

$$= \limsup_{\substack{(u,v)\to(\bar{x},(\bar{y}_1,1-\bar{y}_1)),\|u-\bar{x}\|_l<\lambda \\ v=(v_1,v_2)\in\overline{(u_1,1-u_1),(0,1)}}} \frac{x_1^*(u_1-\bar{x}_1)+x_2^*(u_2-\bar{x}_2)+(y_2^*-y_1^*)(v_1-\bar{y}_1)}{\|(u_1,u_2)-(\bar{x}_1,\bar{x}_2)\|_l+\|(v_1-\bar{y}_1,-(v_1-\bar{y}_1))\|_M}$$

$$= \limsup_{\substack{(u,v)\to(\bar{x},(\bar{y}_1,1-\bar{y}_1)),\|u-\bar{x}\|_l<\lambda \\ v=(v_1,v_2)\in\overline{(u_1,1-u_1),(0,1)}}} \frac{\langle x^*,u-\bar{x}\rangle+(y_2^*-y_1^*)(v_1-\bar{y}_1)}{\|u-\bar{x}\|_l+|v_1-\bar{y}_1|}. \tag{4.28}$$

If $y_2^*-y_1^*>0$, then we take $u_1-\bar{x}_1>0$ and $v_1-\bar{y}_1>0$. Otherwise, if $y_2^*-y_1^*<0$, then we take $u_1-\bar{x}_1<0$ and $v_1-\bar{y}_1<0$. Substituting these to (4.28), we obtain

$$\limsup_{\substack{(u,v)\to(\bar{x},\bar{y}) \\ u\in\mathbb{R}_l^2 \text{ and } v\in P_{\mathbb{B}_l}(u)}} \frac{\langle x^*,u-\bar{x}\rangle-\langle y^*,v-\bar{y}\rangle}{\|u-\bar{x}\|_l+\|v-\bar{y}\|_M}$$

$$\geq \limsup_{\substack{(u,v)\to(\bar{x},(\bar{y}_1,1-\bar{y}_1)),\|u-\bar{x}\|_l<\lambda \\ (v_1,v_2)\in\overline{(u_1,1-u_1),(0,1)}}} \frac{\langle x^*,u-\bar{x}\rangle+|y_2^*-y_1^*||v_1-\bar{y}_1|}{\|u-\bar{x}\|_l+|v_1-\bar{y}_1|}. \tag{4.29}$$

For $\|u-\bar{x}\|_l<\lambda$, by (4.24), $(\bar{y}_1,\bar{y}_2)\in P_{\mathbb{B}_l}(u_1,u_2)=\overline{(u_1,1-u_1),(0,1)}$. As $u$ approaches to $\bar{x}$, $u_1$ will approach to $\bar{x}_1$, which satisfies that $\bar{y}_1<\bar{x}_1$. Then, we can always take $(v_1,v_2)\in\overline{(u_1,1-u_1),(0,1)}$ such that $v_1$ is arbitrarily close to $\bar{y}_1$ and

$$\|u-\bar{x}\|_l<\frac{|y_2^*-y_1^*|}{2\|x^*\|_M}|v_1-\bar{y}_1|<\lambda. \tag{4.30}$$

Substituting (4.30) to (4.29), we have

$$\limsup_{\substack{(u,v)\to(\bar{x},\bar{y}) \\ u\in\mathbb{R}_l^2 \text{ and } v\in P_{\mathbb{B}_l}(u)}} \frac{\langle x^*,u-\bar{x}\rangle-\langle y^*,v-\bar{y}\rangle}{\|u-\bar{x}\|_l+\|v-\bar{y}\|_M}$$

$$\geq \limsup_{\substack{(u,v)\to(\bar{x},(\bar{y}_1,1-\bar{y}_1)),\|u-\bar{x}\|_l<\frac{|y_2^*-y_1^*|}{2\|x^*\|_M}|v_1-\bar{y}_1| \\ v=(v_1,v_2)\in\overline{(u_1,1-u_1),(0,1)}\ni(\bar{y}_1,1-\bar{y}_1)}} \frac{-\|x^*\|_M\frac{|y_2^*-y_1^*|}{2\|x^*\|_M}|v_1-\bar{y}_1|+|y_2^*-y_1^*||v_1-\bar{y}_1|}{\|u-\bar{x}\|_l+|v_1-\bar{y}_1|}$$

$$\geq \limsup_{\substack{(u,v)\to(\bar{x},(\bar{y}_1,1-\bar{y}_1)),\|u-\bar{x}\|_l<\frac{|y_2^*-y_1^*|}{2\|x^*\|_M}|v_1-\bar{y}_1| \\ v=(v_1,v_2)\in\overline{(u_1,1-u_1),(0,1)}\ni(\bar{y}_1,1-\bar{y}_1)}} \frac{-\|x^*\|_M\frac{|y_2^*-y_1^*|}{2\|x^*\|_M}|v_1-\bar{y}_1|+|y_2^*-y_1^*||v_1-\bar{y}_1|}{\frac{|y_2^*-y_1^*|}{2\|x^*\|_M}|v_1-\bar{y}_1|+|v_1-\bar{y}_1|}$$

$$= \frac{\frac{1}{2}|y_2^*-y_1^*|}{\frac{|y_2^*-y_1^*|}{2\|x^*\|_M}+1}>0.$$

This implies that, for $y^*=(y_1^*,y_2^*)\in\mathbb{R}_M^2\backslash\{\theta\}$, if $y_1^*\neq y_2^*$, then, for $(\bar{y}_1,\bar{y}_2)\in P_{\mathbb{B}_l}(\bar{x}_1,\bar{x}_2)=\overline{(\bar{x}_1,1-\bar{x}_1),(0,1)}$ with $0<\bar{y}_1<\bar{x}_1$,

$$x^*\notin\widehat{D}^*P_{\mathbb{B}_l}(\bar{x},\bar{y})(y^*), \text{ for any } x^*\in\mathbb{R}_M^2\backslash\{\theta\}. \tag{4.31}$$

By (4.27) and (4.31), we obtain that

$$\widehat{D}^* P_{\mathbb{B}_l}(\bar{x}, \bar{y})(y^*) = \emptyset, \text{ for any } y^* = (y_1^*, y_2^*) \in \mathbb{R}_M^2 \text{ with } y_1^* \neq y_2^*.$$

Proof (II) of (a) in (ii). This part is proved by (B) and the above (4.31).

Proof of (b) in (ii). The proofs of (I) and (II) of (b) in (ii) are similar to the proofs of (II) and (I) of (b) in (i), respectively, which are omitted here.

Proof of (I) of (c) in (ii). In this case, recall that $\bar{x} = (\bar{x}_1, \bar{x}_2)$ with $0 < \bar{x}_1 < 1, \bar{x}_2 > 1$ and $P_{\mathbb{B}_l}(\bar{x}_1, \bar{x}_2) = \overline{(\bar{x}_1, 1 - \bar{x}_1), (0,1)}$. In particular, let $\bar{y} = (\bar{y}_1, \bar{y}_2) = (\bar{x}_1, 1 - \bar{x}_1) \in P_{\mathbb{B}_l}(\bar{x}_1, \bar{x}_2)$. Let $y^* = (y_1^*, y_2^*) \in \mathbb{R}_M^2$ and $x^* = (x_1^*, x_2^*) \in \mathbb{R}_M^2$. We calculate

$$\limsup_{\substack{(u,v)\to(\bar{x},\bar{y}) \\ u\in\mathbb{R}_l^2 \text{ and } v \in P_{\mathbb{B}_l}(u)}} \frac{\langle x^*, u-\bar{x}\rangle - \langle y^*, v-\bar{y}\rangle}{\|u-\bar{x}\|_l + \|v-\bar{y}\|_M}$$

$$= \limsup_{\substack{(u,v)\to((\bar{x}_1,\bar{x}_2),(\bar{x}_1,1-\bar{x}_1)), \|u-\bar{x}\|_l<\lambda \\ v=(v_1,v_2) \in \overline{(u_1,1-u_1),(0,1)}}} \frac{\langle (x_1^*,x_2^*),(u_1,u_2)-(\bar{x}_1,\bar{x}_2)\rangle - \langle (y_1^*,y_2^*),(v_1,v_2)-(\bar{x}_1,1-\bar{x}_1)\rangle}{\|(u_1,u_2)-(\bar{x}_1,\bar{x}_2)\|_l + \|(v_1,v_2)-(\bar{x}_1,1-\bar{x}_1)\|_M}$$

$$= \limsup_{\substack{(u,v)\to((\bar{x}_1,\bar{x}_2),(\bar{x}_1,1-\bar{x}_1)), \|u-\bar{x}\|_l<\lambda \\ v=(v_1,v_2) \in \overline{(u_1,1-u_1),(0,1)}}} \frac{\langle (x_1^*,x_2^*),(u_1,u_2)-(\bar{x}_1,\bar{x}_2)\rangle - \langle (y_1^*,y_2^*),(v_1,1-v_1)-(\bar{x}_1,1-\bar{x}_1)\rangle}{\|(u_1,u_2)-(\bar{x}_1,\bar{x}_2)\|_l + \|(v_1,1-v_1)-(\bar{x}_1,1-\bar{x}_1)\|_M}$$

$$= \limsup_{\substack{(u,v)\to((\bar{x}_1,\bar{x}_2),(\bar{x}_1,1-\bar{x}_1)), \|u-\bar{x}\|_l<\lambda \\ v=(v_1,v_2) \in \overline{(u_1,1-u_1),(0,1)}}} \frac{\langle (x_1^*,x_2^*),(u_1-\bar{x}_1,u_2-\bar{x}_2)\rangle - \langle (y_1^*,y_2^*),(v_1-\bar{x}_1,-(v_1-\bar{x}_1))\rangle}{\|(u_1,u_2)-(\bar{x}_1,\bar{x}_2)\|_l + \|(v_1-\bar{x}_1,-(v_1-\bar{x}_1))\|_M}$$

$$= \limsup_{\substack{(u,v)\to((\bar{x}_1,\bar{x}_2),(\bar{x}_1,1-\bar{x}_1)), \|u-\bar{x}\|_l<\lambda \\ v=(v_1,v_2) \in \overline{(u_1,1-u_1),(0,1)}}} \frac{x_1^*(u_1-\bar{x}_1) + x_2^*(u_2-\bar{x}_2) - (y_1^*-y_2^*)(v_1-\bar{x}_1)}{\|(u_1,u_2)-(\bar{x}_1,\bar{x}_2)\|_l + \|(v_1-\bar{x}_1,-(v_1-\bar{x}_1))\|_M}$$

$$= \limsup_{\substack{(u,v)\to((\bar{x}_1,\bar{x}_2),(\bar{x}_1,1-\bar{x}_1)), \|u-\bar{x}\|_l<\lambda \\ v=(v_1,v_2) \in \overline{(u_1,1-u_1),(0,1)}}} \frac{x_1^*(u_1-\bar{x}_1) + x_2^*(u_2-\bar{x}_2) - (y_1^*-y_2^*)(v_1-\bar{x}_1)}{\|(u_1,u_2)-(\bar{x}_1,\bar{x}_2)\|_l + |v_1-\bar{x}_1|}. \tag{4.32}$$

In particular, for $y^* = (y_1^*, y_2^*) \in \mathbb{R}_M^2$ with $y_1^* = y_2^*$ in (4.32), we consider two cases.

Case 1. Let $x^* = (x_1^*, x_2^*) = \theta$. By (4.32), we have

$$\limsup_{\substack{(u,v)\to(\bar{x},\bar{y}) \\ u\in\mathbb{R}_l^2 \text{ and } v \in P_{\mathbb{B}_l}(u)}} \frac{\langle x^*, u-\bar{x}\rangle - \langle y^*, v-\bar{y}\rangle}{\|u-\bar{x}\|_l + \|v-\bar{y}\|_M} = 0.$$

This implies that, for any $y^* = (y_1^*, y_2^*) \in \mathbb{R}_M^2$, if $y_1^* = y_2^*$, then

$$\theta \in \widehat{D}^* P_{\mathbb{B}_l}\big((\bar{x}_1, \bar{x}_2), (\bar{x}_1, 1 - \bar{x}_1)\big)(y^*). \tag{4.33}$$

Case 2. Let $x^* = (x_1^*, x_2^*) \neq \theta$. By (4.32), we have

$$\limsup_{\substack{(u,v)\to(\bar{x},\bar{y}) \\ u\in\mathbb{R}_l^2 \text{ and } v\in P_{\mathbb{B}_l}(u)}} \frac{\langle x^*, u-\bar{x}\rangle - \langle y^*, v-\bar{y}\rangle}{\|u-\bar{x}\|_l + \|v-\bar{y}\|_M}$$

$$= \limsup_{\substack{(u,v)\to((\bar{x}_1,\bar{x}_2),(\bar{x}_1,1-\bar{x}_1)), \|u-\bar{x}\|_l<\lambda \\ (v_1,v_2)\in\overline{(u_1,1-u_1),(0,1)}}} \frac{x_1^*(u_1-\bar{x}_1) + x_2^*(u_2-\bar{x}_2)}{\|(u_1,u_2)-(\bar{x}_1,\bar{x}_2)\|_l + |v_1-\bar{x}_1|}.$$

By the condition that $(x_1^*, x_2^*) \neq \theta$, without loose of the generality, we suppose $x_1^* \neq 0$. Let $t$ be real. Let $u(t) = (\bar{x}_1 - t, \bar{x}_2)$ and $v(t) = (\bar{x}_1 - t, 1 - (\bar{x}_1 - t))$. Here we assume that $|t|$ is small enough such that $0 < \bar{x}_1 - t < 1$. Substituting these to (4.32), we have

$$\limsup_{\substack{(u,v)\to(\bar{x},\bar{y}) \\ u\in\mathbb{R}_l^2 \text{ and } v\in P_{\mathbb{B}_l}(u)}} \frac{\langle x^*, u-\bar{x}\rangle - \langle y^*, v-\bar{y}\rangle}{\|u-\bar{x}\|_l + \|v-\bar{y}\|_M}$$

$$\geq \limsup_{\substack{(u(t),v(t))\to((\bar{x}_1,\bar{x}_2),(\bar{x}_1,1-\bar{x}_1)), \|u(t)-\bar{x}\|_l<\lambda \\ v(t) = (\bar{x}_1-t, 1-(\bar{x}_1-t))}} \frac{x_1^*(\bar{x}_1-t-\bar{x}_1) + x_2^*(\bar{x}_2-\bar{x}_2)}{\|(\bar{x}_1+t,\bar{x}_2)-(\bar{x}_1,\bar{x}_2)\|_l + |\bar{x}_1+t-\bar{x}_1|}$$

$$\geq \limsup_{\substack{(u(t),v(t))\to((\bar{x}_1,\bar{x}_2),(\bar{x}_1,1-\bar{x}_1)), \|u(t)-\bar{x}\|_l<\lambda \\ v(t) = (\bar{x}_1-t, 1-(\bar{x}_1-t))}} \frac{-x_1^* t}{2|t|}$$

$$= \frac{|x_1^*|}{2} > 0 \quad \text{(under the condition } x_1^* \neq 0\text{)}.$$

This implies that for $y^* = (y_1^*, y_2^*)$, if $y_1^* = y_2^*$, then

$$x^* \notin \widehat{D}^* F\big((\bar{x}_1, \bar{x}_2), (\bar{x}_1, 1 - \bar{x}_1)\big)(y^*), \text{ for any } x^* \neq \theta. \tag{4.34}$$

By (4.33) and (4.34), for any $y^* = (y_1^*, y_2^*)$, if $y_1^* = y_2^*$, then

$$\widehat{D}^* P_{\mathbb{B}_l}\big((\bar{x}_1, \bar{x}_2), (\bar{x}_1, 1 - \bar{x}_1)\big)(y^*) = \{\theta\}.$$

Proof of (II) of (c) in (ii). Let $(\bar{y}_1, \bar{y}_2) = (\bar{x}_1, 1 - \bar{x}_1)$. Let $y^* = (y_1^*, y_2^*) \in \mathbb{R}_M^2$ and $x^* = (x_1^*, x_2^*) \in \mathbb{R}_M^2$. There are two cases for consideration.

Case 1. Suppose $y_1^* < y_2^*$. Let $x^* = (x_1^*, x_2^*) \in \mathbb{R}_M^2$. If $\|x^*\|_M < \frac{y_2^* - y_1^*}{2}$, then by (4.32), we have

$$\limsup_{\substack{(u,v)\to(\bar{x},\bar{y}) \\ u\in\mathbb{R}_l^2 \text{ and } v\in P_{\mathbb{B}_l}(u)}} \frac{\langle x^*, u-\bar{x}\rangle - \langle y^*, v-\bar{y}\rangle}{\|u-\bar{x}\|_l + \|v-\bar{y}\|_M}$$

$$= \limsup_{\substack{(u,v)\to((\bar{x}_1,\bar{x}_2),(\bar{x}_1,1-\bar{x}_1)), \|u-\bar{x}\|_l<\lambda \\ (v_1,v_2)\in\overline{(u_1,1-u_1),(0,1)}}} \frac{x_1^*(u_1-\bar{x}_1) + x_2^*(u_2-\bar{x}_2) - (y_1^*-y_2^*)(v_1-\bar{x}_1)}{\|(u_1,u_2)-(\bar{x}_1,\bar{x}_2)\|_l + |v_1-\bar{x}_1|}$$

$$= \limsup_{\substack{(u,v)\to((\bar{x}_1,\bar{x}_2),(\bar{x}_1,1-\bar{x}_1)), \|u-\bar{x}\|_l<\lambda \\ (v_1,v_2)\in\overline{(u_1,1-u_1),(0,1)}}} \frac{\langle x^*, u-\bar{x}\rangle + (y_2^*-y_1^*)(v_1-\bar{x}_1)}{\|u-\bar{x}\|_l + |v_1-\bar{x}_1|}$$

$$\geq \limsup_{\substack{((u_1,\bar{x}_2),(u_1,1-u_1))\to((\bar{x}_1,\bar{x}_2),(\bar{x}_1,1-\bar{x}_1)),\|u-\bar{x}\|_l<\lambda \\ u_1\downarrow\bar{x}_1}} \frac{\langle x^*,(u_1,\bar{x}_2)-(\bar{x}_1,\bar{x}_2)\rangle+(y_2^*-y_1^*)(u_1-\bar{x}_1)}{\|(u_1,\bar{x}_2)-(\bar{x}_1,\bar{x}_2)\|_l+|u_1-\bar{x}_1|}$$

$$\geq \limsup_{\substack{((u_1,\bar{x}_2),(u_1,1-u_1))\to((\bar{x}_1,\bar{x}_2),(\bar{x}_1,1-\bar{x}_1)),\|u-\bar{x}\|_l<\lambda \\ u_1\downarrow\bar{x}_1}} \frac{-\|x^*\|_M\|(u_1,\bar{x}_2)-(\bar{x}_1,\bar{x}_2)\|_l+(y_2^*-y_1^*)(u_1-\bar{x}_1)}{|u_1-\bar{x}_1|+|u_1-\bar{x}_1|}$$

$$\geq \limsup_{\substack{((u_1,\bar{x}_2),(u_1,1-u_1))\to((\bar{x}_1,\bar{x}_2),(\bar{x}_1,1-\bar{x}_1)),\|u-\bar{x}\|_l<\lambda \\ u_1\downarrow\bar{x}_1}} \frac{-\|x^*\|_M|u_1-\bar{x}_1|+(y_2^*-y_1^*)(u_1-\bar{x}_1)}{|u_1-\bar{x}_1|+|u_1-\bar{x}_1|}$$

$$> \limsup_{\substack{((u_1,\bar{x}_2),(u_1,1-u_1))\to((\bar{x}_1,\bar{x}_2),(\bar{x}_1,1-\bar{x}_1)),\|u-\bar{x}\|_l<\lambda \\ u_1\downarrow\bar{x}_1}} \frac{-\frac{y_2^*-y_1^*}{2}|u_1-\bar{x}_1|+(y_2^*-y_1^*)(u_1-\bar{x}_1)}{|u_1-\bar{x}_1|+|u_1-\bar{x}_1|}$$

$$= \limsup_{\substack{((u_1,\bar{x}_2),(u_1,1-u_1))\to((\bar{x}_1,\bar{x}_2),(\bar{x}_1,1-\bar{x}_1)),\|u-\bar{x}\|_l<\lambda \\ u_1\downarrow\bar{x}_1}} \frac{\frac{y_2^*-y_1^*}{2}|u_1-\bar{x}_1|}{|u_1-\bar{x}_1|+|u_1-\bar{x}_1|}$$

$$= \frac{y_2^*-y_1^*}{4} > 0 \quad \text{(under the condition } y_1^* < y_2^*\text{).}$$

This implies that for $y^* = (y_1^*, y_2^*)$, with $y_1^* < y_2^*$, for any $x^* = (x_1^*, x_2^*) \neq \theta$, if $\|x^*\|_M < \frac{y_2^*-y_1^*}{2}$, then

$$x^* \notin \widehat{D}^* P_{\mathbb{B}_l}\big((\bar{x}_1,\bar{x}_2),(\bar{x}_1,1-\bar{x}_1)\big)(y^*). \tag{4.35}$$

Case 2. Suppose $y_1^* > y_2^*$. Let $x^* = (x_1^*, x_2^*) \in \mathbb{R}_M^2$. If $\|x^*\|_M < \frac{y_1^*-y_2^*}{2}$, then, by (4.32) and similarly to the proof of (4.35), we have

$$\limsup_{\substack{(u,v)\to(\bar{x},\bar{y}) \\ u\in\mathbb{R}_l^2 \text{ and } v\in P_{\mathbb{B}_l}(u)}} \frac{\langle x^*,u-\bar{x}\rangle-\langle y^*,v-\bar{y}\rangle}{\|u-\bar{x}\|_l+\|v-\bar{y}\|_M}$$

$$\geq \limsup_{\substack{((u_1,\bar{x}_2),(u_1,1-u_1))\to((\bar{x}_1,\bar{x}_2),(\bar{x}_1,1-\bar{x}_1)),\|u-\bar{x}\|_l<\lambda \\ u_1\uparrow\bar{x}_1}} \frac{\langle x^*,(u_1,\bar{x}_2)-(\bar{x}_1,\bar{x}_2)\rangle_X+(y_1^*-y_2^*)(\bar{x}_1-u_1)}{\|(u_1,\bar{x}_2)-(\bar{x}_1,\bar{x}_2)\|_l+|u_1-\bar{x}_1|}$$

$$> \limsup_{\substack{((u_1,\bar{x}_2),(u_1,1-u_1))\to((\bar{x}_1,\bar{x}_2),(\bar{x}_1,1-\bar{x}_1)),\|u-\bar{x}\|_l<\lambda \\ u_1\uparrow\bar{x}_1}} \frac{-\frac{y_1^*-y_2^*}{2}|u_1-\bar{x}_1|+(y_1^*-y_2^*)(\bar{x}_1-u_1)}{|u_1-\bar{x}_1|+|u_1-\bar{x}_1|}$$

$$= \frac{y_1^*-y_2^*}{4} > 0 \quad \text{(under the condition } y_1^* > y_2^*\text{).}$$

This implies that for $y^* = (y_1^*, y_2^*)$, with $y_1^* > y_2^*$, for any $x^* \neq \theta$, if $\|x^*\|_M < \frac{y_1^*-y_2^*}{2}$, then

$$x^* \notin \widehat{D}^* P_{\mathbb{B}_l}\big((\bar{x}_1,\bar{x}_2),(\bar{x}_1,1-\bar{x}_1)\big)(y^*).$$

This completes the proof of part (ii).

Proof of (I) of (a) in (iii). Let $\bar{x}_1 = 0$ and $\bar{x}_2 > 1$. In this case, since $P_{\mathbb{B}_l}(0, \bar{x}_2) = (0,1)$, that is a singleton, if $\bar{y} = (\bar{y}_1, \bar{y}_2) \in P_{\mathbb{B}_l}(0, \bar{x}_2)$, we must have $(\bar{y}_1, \bar{y}_2) = (0,1)$. So, for $y^* = (y_1^*, y_2^*) \in \mathbb{R}_M^2$, we consider

$$\widehat{D}^* P_{\mathbb{B}_l}\big((0, \bar{x}_2), (0,1)\big)(y^*).$$

Let $x^* = (x_1^*, x_2^*) \in \mathbb{R}_M^2$, we calculate

$$\limsup_{\substack{(u,v)\to((0,\bar{x}_2),(0,1)) \\ u\in\mathbb{R}_l^2 \text{ and } v\in P_{\mathbb{B}_l}(u)}} \frac{\langle x^*, u-\bar{x}\rangle - \langle y^*, v-\bar{y}\rangle}{\|u-\bar{x}\|_l + \|v-\bar{y}\|_M}$$

$$= \limsup_{\substack{((u_1,u_2),(v_1,v_2))\to((0,\bar{x}_2),(0,1)) \\ v=(v_1,v_2)\in P_{\mathbb{B}_l}(u_1,u_2)}} \frac{x_1^*(u_1-0) + x_2^*(u_2-\bar{x}_2) - (y_1^*(v_1-0) + y_2^*(v_2-1))}{\|(u_1,u_2)-(0,\bar{x}_2)\|_l + \|(v_1,v_2)-(0,1)\|_M}. \tag{4.36}$$

For the above limit, we consider two cases with respect to $u_1$:

Case 1. $u_1 \downarrow 0$ and $u_2 \to \bar{x}_2$. In this case, for $v = (v_1, v_2) \in P_{\mathbb{B}_l}(u_1, u_2) = \overline{(u_1, 1-u_1), (0,1)}$, it implies that $(v_1, v_2) = (v_1, 1 - v_1)$. Substituting these to (4.36), we have

$$\limsup_{\substack{(u,v)\to((0,\bar{x}_2),(0,1)) \\ u\in\mathbb{R}_l^2 \text{ and } v\in P_{\mathbb{B}_l}(u)}} \frac{\langle x^*, u-\bar{x}\rangle - \langle y^*, v-\bar{y}\rangle}{\|u-\bar{x}\|_l + \|v-\bar{y}\|_M}$$

$$= \limsup_{\substack{u_1\downarrow 0 \text{ and } u_2\to\bar{x}_2 \\ (v_1, 1-v_2)\to(0,1),\ v_1\downarrow 0}} \frac{x_1^* u_1 + x_2^*(u_2-\bar{x}_2) - (y_1^* - y_2^*)v_1}{\|(u_1,u_2)-(0,\bar{x}_2)\|_l + \|(v_1,-v_1)\|_M}. \tag{4.37}$$

Case 2. $u_1 \uparrow 0$ with $u_1 > -1$ and $u_2 \to \bar{x}_2$. Then, by Proposition 4.1, for $v = (v_1, v_2) \in P_{\mathbb{B}_l}(u_1, u_2) = \overline{(u_1, 1+u_1), (0,1)}$, it implies that $(v_1, v_2) = (v_1, 1 + v_1)$ with $-1 < v_1 < 0$ and $v_2 - v_1 = 1$. Substituting these to (4.36), we have

$$\limsup_{\substack{(u,v)\to((0,\bar{x}_2),(0,1)) \\ u\in\mathbb{R}_l^2 \text{ and } v\in P_{\mathbb{B}_l}(u)}} \frac{\langle x^*, u-\bar{x}\rangle - \langle y^*, v-\bar{y}\rangle}{\|u-\bar{x}\|_l + \|v-\bar{y}\|_M}$$

$$= \limsup_{\substack{u_1\uparrow 0 \text{ and } u_2\to\bar{x}_2 \\ (v_1, 1+v_2)\to(0,1),\ v_1\uparrow 0}} \frac{x_1^*(u_1-0) + x_2^*(u_2-\bar{x}_2) - (y_1^*(v_1) + y_2^* v_1)}{\|(u_1,u_2)-(0,\bar{x}_2)\|_l + \|(v_1,-v_1)\|_M}$$

$$= \limsup_{\substack{u_1\uparrow 0 \text{ and } u_2\to\bar{x}_2 \\ (v_1, 1+v_2)\to(0,1),\ v_1\uparrow 0}} \frac{x_1^* u_1 + x_2^*(u_2-\bar{x}_2) - (y_1^* + y_2^*)v_1}{\|(u_1,u_2)-(0,\bar{x}_2)\|_l + \|(v_1,v_1)\|_M}. \tag{4.38}$$

Let $y_1^* = y_2^* < 0$. For any $x^* = (x_1^*, 0) \in \mathbb{R}_M^2$, by (A) in this theorem, we only prove

$$x_1^* \le 0 \quad \Longrightarrow \quad (x_1^*, 0) \in \widehat{D}^* P_{\mathbb{B}_l}\big((0, \bar{x}_2), (0,1)\big)(y^*). \tag{4.39}$$

Proof of (4.39). By (4.37) and (4.38), we have

$$\limsup_{\substack{(u,v)\to((0,\bar{x}_2),(0,1))\\ u\in\mathbb{R}^2_l \text{ and } v\in P_{\mathbb{B}_l}(u)}} \frac{\langle x^*,u-\bar{x}\rangle-\langle y^*,v-\bar{y}\rangle}{\|u-\bar{x}\|_l+\|v-\bar{y}\|_M}$$

$$= \max\left\{\limsup_{\substack{u_1\downarrow 0 \text{ and } u_2\to\bar{x}_2\\ (v_1,1-v_2)\to(0,1), v_1\downarrow 0}} \frac{\langle x^*,u-\bar{x}\rangle_X-(y_1^*-y_2^*)v_1}{\|(u_1,u_2)-(0,\bar{x}_2)\|_l+\|(v_1,-v_1)\|_M}, \limsup_{\substack{u_1\uparrow 0 \text{ and } u_2\to\bar{x}_2\\ (v_1,1+v_2)\to(0,1), v_1\uparrow 0}} \frac{-(y_1^*+y_2^*)v_1}{\|(u_1,u_2)-(0,\bar{x}_2)\|_l+\|(v_1,v_1)\|_M}\right\}$$

$$= \max\left\{\limsup_{\substack{u_1\downarrow 0 \text{ and } u_2\to\bar{x}_2\\ (v_1,1-v_2)\to(0,1),\ v_1\downarrow 0}} \frac{x_1^*(u_1-0)+x_2^*(u_2-\bar{x}_2)}{\|(u_1,u_2)-(0,\bar{x}_2)\|_l+\|(v_1,-v_1)\|_M}, \limsup_{\substack{u_1\uparrow 0 \text{ and } u_2\to\bar{x}_2\\ (v_1,1+v_2)\to(0,1), v_1\uparrow 0}} \frac{-(y_1^*+y_2^*)v_1}{\|(u_1,u_2)-(0,\bar{x}_2)\|_l+\|(v_1,v_1)\|_M}\right\}$$

$$= \max\left\{\limsup_{\substack{u_1\downarrow 0 \text{ and } u_2\to\bar{x}_2\\ (v_1,1-v_2)\to(0,1),\ v_1\downarrow 0}} \frac{x_1^*u_1}{\|(u_1,u_2)-(0,\bar{x}_2)\|_l+\|(v_1,-v_1)\|_M}, \limsup_{\substack{u_1\uparrow 0 \text{ and } u_2\to\bar{x}_2\\ (v_1,1+v_2)\to(0,1), v_1\uparrow 0}} \frac{-(y_1^*+y_2^*)v_1}{\|(u_1,u_2)-(0,\bar{x}_2)\|_l+\|(v_1,v_1)\|_M}\right\}$$

$\leq \max\{0,0\}$ (by $x_1^*\leq 0$ and $y_1^*=y_2^*<0$)

$\leq 0.$

This proves (4.39)

Proof of (II) of (a) in (iii). If $y_1^*=y_2^*>0$, then, for $x^*=(x_1^*,x_2^*)$, we prove

$$\|x^*\|_M<y_1^* \implies x^*\notin \widehat{D}^*P_{\mathbb{B}_l}\big((0,\bar{x}_2),(0,1)\big)(y^*). \tag{4.40}$$

Proof of (4.40). Similarly to the proof of (4.39), by (4.37) and (4.38), we have

$$\limsup_{\substack{(u,v)\to((0,\bar{x}_2),(0,1))\\ u\in\mathbb{R}^2_l \text{ and } v\in P_{\mathbb{B}_l}(u)}} \frac{\langle x^*,u-\bar{x}\rangle-\langle y^*,v-\bar{y}\rangle}{\|u-\bar{x}\|_l+\|v-\bar{y}\|_M}$$

$$= \max\left\{\limsup_{\substack{u_1\downarrow 0 \text{ and } u_2\to\bar{x}_2\\ (v_1,1-v_2)\to(0,1), v_1\downarrow 0}} \frac{\langle x^*,u-\bar{x}\rangle-(y_1^*-y_2^*)v_1}{\|(u_1,u_2)-(0,\bar{x}_2)\|_l+\|(v_1,-v_1)\|_M}, \limsup_{\substack{u_1\uparrow 0 \text{ and } u_2\to\bar{x}_2\\ (v_1,1+v_2)\to(0,1), v_1\uparrow 0}} \frac{\langle x^*,u-\bar{x}\rangle-(y_1^*+y_2^*)v_1}{\|(u_1,u_2)-(0,\bar{x}_2)\|_l+\|(v_1,v_1)\|_M}\right\}$$

$$\geq \max\left\{\limsup_{\substack{u_1\downarrow 0 \text{ and } u_2\to\bar{x}_2\\ (v_1,1-v_2)\to(0,1), v_1\downarrow 0}} \frac{\langle x^*,u-\bar{x}\rangle}{\|(u_1,u_2)-(0,\bar{x}_2)\|_l+\|(v_1,-v_1)\|_M}, \limsup_{\substack{u_1\uparrow 0 \text{ and } u_2=\bar{x}_2\\ (v_1,1+v_2)\to(0,1), v_1=u_1}} \frac{-\|x^*\|_M|u_1|-2y_1^*u_1}{|u_1|+|u_1|}\right\}$$

$$= \max\left\{\limsup_{\substack{u_1\downarrow 0 \text{ and } u_2\to\bar{x}_2\\ (v_1,1-v_2)\to(0,1), v_1\downarrow 0}} \frac{\langle x^*,u-\bar{x}\rangle}{\|(u_1,u_2)-(0,\bar{x}_2)\|_l+\|(v_1,-v_1)\|_M}, \limsup_{\substack{u_1\uparrow 0 \text{ and } u_2=\bar{x}_2\\ (v_1,1+v_2)\to(0,1), v_1=u_1}} \frac{-\|x^*\|_M|u_1|+2y_1^*|u_1|}{|u_1|+|u_1|}\right\}$$

$\geq \frac{y_1^*}{2}>0$ (by $\|x^*\|_M<y_1^*$)

This proves (4.40).

Proof of (b) in (iii). Let $y_1^* + y_2^* = 0$. If $y_1^* > y_2^*$, then for any $x^* = (x_1^*, 0) \in \mathbb{R}_M^2$, we prove

$$x_1^* \le y_2^* - y_1^* \implies (x_1^*, 0) \notin \widehat{D}^* P_{\mathbb{B}_l}\big((0, \bar{x}_2), (0, 1)\big)(y^*). \tag{4.41}$$

Similarly to the proof of (4.40), by (4.37) and (4.38), we have

$$\limsup_{\substack{(u,v)\to((0,\bar{x}_2),(0,1))\\ u\in\mathbb{R}_l^2 \text{ and } v\in P_{\mathbb{B}_l}(u)}} \frac{\langle x^*, u-\bar{x}\rangle - \langle y^*, v-\bar{y}\rangle}{\|u-\bar{x}\|_l + \|v-\bar{y}\|_M}$$

$$= \max\left\{ \limsup_{\substack{u_1\downarrow 0 \text{ and } u_2\to\bar{x}_2\\ (v_1,1-v_2)\to(0,1), v_1\downarrow 0}} \frac{\langle x^*, u-\bar{x}\rangle - (y_1^*-y_2^*)v_1}{\|(u_1,u_2)-(0,\bar{x}_2)\|_l + \|(v_1,-v_1)\|_M}, \limsup_{\substack{u_1\uparrow 0 \text{ and } u_2\to\bar{x}_2\\ (v_1,1+v_2)\to(0,1), v_1\uparrow 0}} \frac{\langle x^*, u-\bar{x}\rangle - (y_1^*+y_2^*)v_1}{\|(u_1,u_2)-(0,\bar{x}_2)\|_l + \|(v_1,v_1)\|_M} \right\}$$

$$= \max\left\{ \limsup_{\substack{u_1\downarrow 0 \text{ and } u_2\to\bar{x}_2\\ (v_1,1-v_2)\to(0,1), v_1\downarrow 0}} \frac{x_1^*(u_1-0) - (y_1^*-y_2^*)v_1}{\|(u_1,u_2)-(0,\bar{x}_2)\|_l + \|(v_1,-v_1)\|_M}, \limsup_{\substack{u_1\uparrow 0 \text{ and } u_2\to\bar{x}_2\\ (v_1,1+v_2)\to(0,1), v_1\uparrow 0}} \frac{x_1^*(u_1-0)}{\|(u_1,u_2)-(0,\bar{x}_2)\|_l + \|(v_1,v_1)\|_M} \right\}$$

$$= \max\left\{ \limsup_{\substack{u_1\downarrow 0 \text{ and } u_2\to\bar{x}_2\\ 0<v_1\le u_1}} \frac{x_1^*u_1 - (y_1^*-y_2^*)v_1}{\|(u_1,u_2)-(0,\bar{x}_2)\|_l + \|(v_1,-v_1)\|_M}, \limsup_{\substack{u_1\uparrow 0 \text{ and } u_2\to\bar{x}_2\\ (v_1,1+v_2)\to(0,1), v_1\uparrow 0}} \frac{x_1^*u_1}{\|(u_1,u_2)-(0,\bar{x}_2)\|_l + \|(v_1,v_1)\|_M} \right\}$$

$$\ge \max\left\{ \limsup_{\substack{u_1\downarrow 0 \text{ and } u_2\to\bar{x}_2\\ 0<v_1\le u_1}} \frac{(y_2^*-y_1^*)u_1 - (y_1^*-y_2^*)v_1}{\|(u_1,u_2)-(0,\bar{x}_2)\|_l + \|(v_1,-v_1)\|_M}, \limsup_{\substack{u_1\uparrow 0 \text{ and } u_2=\bar{x}_2\\ (v_1,1+v_2)\to(0,1), v_1\uparrow 0}} \frac{x_1^*u_1}{\|(u_1,u_2)-(0,\bar{x}_2)\|_l + \|(v_1,v_1)\|_M} \right\}$$

$$\ge \max\left\{ \limsup_{\substack{u_1\downarrow 0 \text{ and } u_2=\bar{x}_2\\ 0<v_1\le u_1}} \frac{(y_2^*-y_1^*)(u_1+v_1)}{\|(u_1,u_2)-(0,\bar{x}_2)\|_l + \|(v_1,-v_1)\|_M}, \limsup_{\substack{u_1\uparrow 0 \text{ and } u_2=\bar{x}_2\\ v_1=u_1,(v_1,1+v_2)\to(0,1)}} \frac{|x_1^*u_1|}{\|(u_1,u_2)-(0,\bar{x}_2)\|_l + \|(v_1,v_1)\|_M} \right\}$$

$$\ge \max\left\{ 0, \limsup_{\substack{u_1\uparrow 0 \text{ and } u_2=\bar{x}_2\\ v_1=u_1,(v_1,1+v_2)\to(0,1)}} \frac{|x_1^*u_1|}{|u_1|+|u_1|} \right\}$$

$$= \frac{|x_1^*|}{2}$$

$> 0$ (by $x_1^* \le y_2^* - y_1^* < 0$).

This proves (4.41), which completes the proof of part (iii).

Proof of (iv). For any $(\bar{x}_1, \bar{x}_2) \in \mathbb{R}_l^2$ with $0 < \bar{x}_1 < 1, 0 < \bar{x}_2 < 1$ and $\bar{x}_1 + \bar{x}_2 > 1$, by Proposition 4.1,

$$P_{\mathbb{B}_l}(\bar{x}_1, \bar{x}_2) = \overline{(\bar{x}_1, 1-\bar{x}_1), (1-\bar{x}_2, \bar{x}_2)}.$$

Then, the proof of part (iv) is similar to the proof of (i) and it is omitted here.

Proof (I) of (a) in (v). Let $\bar{x} = (\bar{x}_1, \bar{x}_2) \in \mathbb{R}_l^2$ with $0 \le \bar{x}_1 \le 1, 0 \le \bar{x}_2 \le 1$ and $\bar{x}_1 + \bar{x}_2 = 1$. Since $\bar{x} \in \mathbb{B}_l$, we have $P_{\mathbb{B}_l}(\bar{x}) = \bar{x}$, which is a singleton. Let $\bar{y} = (\bar{y}_1, \bar{y}_2) \in P_{\mathbb{B}_l}(\bar{x})$. we must have $\bar{y} = \bar{x}$. For $y^* = (y_1^*, y_2^*) \in \mathbb{R}_M^2$, we consider

$$\widehat{D}^* P_{\mathbb{B}_l}(\bar{x}, \bar{x})(y^*).$$

Let $x^* = (x_1^*, x_2^*) \in \mathbb{R}_M^2$, We calculate

$$\limsup_{\substack{(u,v)\to(\bar{x},\bar{x}) \\ u\in\mathbb{R}_l^2 \text{ and } v\in P_{\mathbb{B}_l}(u)}} \frac{\langle x^*, u-\bar{x}\rangle - \langle y^*, v-\bar{x}\rangle}{\|u-\bar{x}\|_l + \|v-\bar{y}\|_M}$$

$$= \limsup_{\substack{((u_1,u_2),(v_1,v_2))\to((\bar{x}_1,\bar{x}_2),(\bar{x}_1,\bar{x}_2)) \\ (v_1,v_2)\in P_{\mathbb{B}_l}(u_1,u_2)}} \frac{x_1^*(u_1-\bar{x}_1) + x_2^*(u_2-\bar{x}_2) - (y_1^*(v_1-\bar{x}_1) + y_2^*(v_2-\bar{x}_2))}{\|(u_1,u_2)-(\bar{x}_1,\bar{x}_2)\|_l + \|(v_1,v_2)-(\bar{x}_1,\bar{x}_2)\|_M}. \tag{4.42}$$

We consider two cases for the limit (4.42).

Case 1. $\|(u_1, u_2)\|_l > 1$. In this case, $P_{\mathbb{B}_l}(u_1, u_2) = \overline{(u_1, 1-u_1), (1-u_2, u_2)}$. Then, every $(v_1, v_2) \in P_{\mathbb{B}_l}(u_1, u_2)$ satisfies $(v_1, v_2) = (v_1, 1-v_1)$ and $u_1 \le v_1 \le 1-u_2$. By $\bar{x}_1 + \bar{x}_2 = 1$, substituting these to (4.42), we have

$$\limsup_{\substack{(u,v)\to((\bar{x}_1,\bar{x}_2),(\bar{x}_1,\bar{x}_2)) \\ u\in\mathbb{R}_l^2 \text{ and } v\in P_{\mathbb{B}_l}(u)}} \frac{\langle x^*, u-\bar{x}\rangle - \langle y^*, v-\bar{x}\rangle}{\|u-\bar{x}\|_l + \|v-\bar{y}\|_M}$$

$$= \limsup_{\substack{((u_1,u_2),(v_1,v_2))\to((\bar{x}_1,\bar{x}_2),(\bar{x}_1,\bar{x}_2)) \\ \|(u_1,u_2)\|_l>1,\ u_1\le v_1\le 1-u_2}} \frac{x_1^*(u_1-\bar{x}_1) + x_2^*(u_2-\bar{x}_2) - (y_1^*(v_1-\bar{x}_1) + y_2^*(v_2-\bar{x}_2))}{\|(u_1,u_2)-(\bar{x}_1,\bar{x}_2)\|_l + \|(v_1,v_2)-(\bar{x}_1,\bar{x}_2)\|_M}$$

$$= \limsup_{\substack{((u_1,u_2),(v_1,v_2))\to((\bar{x}_1,\bar{x}_2),(\bar{x}_1,\bar{x}_2)) \\ \|(u_1,u_2)\|_l>1,\ u_1\le v_1\le 1-u_2}} \frac{x_1^*(u_1-\bar{x}_1) + x_2^*(u_2-(1-\bar{x}_1)) - (y_1^*(v_1-\bar{x}_1) + y_2^*(\bar{x}_1-v_1))}{\|(u_1,u_2)-(\bar{x}_1,1-\bar{x}_1)\|_l + \|(v_1-\bar{x}_1,-(v_1-\bar{x}_1))\|_M}$$

$$= \limsup_{\substack{((u_1,u_2),(v_1,v_2))\to((\bar{x}_1,\bar{x}_2),(\bar{x}_1,\bar{x}_2)) \\ \|(u_1,u_2)\|_l>1,\ u_1\le v_1\le 1-u_2}} \frac{x_1^*(u_1-\bar{x}_1) + x_2^*(u_2-(1-\bar{x}_1)) - (y_1^*-y_2^*)(v_1-\bar{x}_1)}{\|(u_1,u_2)-(\bar{x}_1,1-\bar{x}_1)\|_l + \|(v_1-\bar{x}_1,-(v_1-\bar{x}_1))\|_M}$$

$$= \limsup_{\substack{((u_1,u_2),(v_1,v_2))\to((\bar{x}_1,\bar{x}_2),(\bar{x}_1,\bar{x}_2)) \\ \|(u_1,u_2)\|_l>1,\ u_1\le v_1\le 1-u_2}} \frac{x_1^*(u_1-\bar{x}_1) + x_2^*(u_2-\bar{x}_2) - (y_1^*-y_2^*)(v_1-\bar{x}_1)}{\|(u_1,u_2)-(\bar{x}_1,\bar{x}_2)\|_l + \|(v_1-\bar{x}_1,-(v_1-\bar{x}_1))\|_M}. \tag{4.43}$$

Case 2. $\|(u_1, u_2)\|_l \le 1$. In this case, $(u_1, u_2) \in \mathbb{B}_l$ and $P_{\mathbb{B}_l}(u_1, u_2) = (u_1, u_2)$. Then, every $(v_1, v_2) \in P_{\mathbb{B}_l}(u_1, u_2)$ satisfies $(v_1, v_2) = (u_1, u_2)$. Substituting these to (4.42), we have

$$\limsup_{\substack{(u,v)\to((\bar{x}_1,\bar{x}_2),(\bar{x}_1,\bar{x}_2)) \\ u\in\mathbb{R}_l^2 \text{ and } v\in P_{\mathbb{B}_l}(u)}} \frac{\langle x^*, u-\bar{x}\rangle - \langle y^*, v-\bar{x}\rangle}{\|u-\bar{x}\|_l + \|v-\bar{y}\|_M}$$

$$= \limsup_{\substack{((u_1,u_2),(v_1,v_2))\to((\bar{x}_1,\bar{x}_2),(\bar{x}_1,\bar{x}_2))\\ \|(u_1,u_2)\|_l\le 1,\ (v_1,v_2)=(u_1,u_2)}} \frac{x_1^*(u_1-\bar{x}_1)+x_2^*\big(u_2-(1-\bar{x}_1)\big)-\big(y_1^*(u_1-\bar{x}_1)+y_2^*\big(u_2-(1-\bar{x}_1)\big)\big)}{\|(u_1,u_2)-(\bar{x}_1,\bar{x}_2)\|_l+\|(u_1,u_2)-(\bar{x}_1,\bar{x}_2)\|_M}$$

$$= \limsup_{\substack{(u_1,u_2)\to(\bar{x}_1,\bar{x}_2)\\ \|(u_1,u_2)\|_l\le 1}} \frac{(x_1^*-y_1^*)(u_1-\bar{x}_1)+(x_2^*-y_2^*)\big(u_2-(1-\bar{x}_1)\big)}{\|(u_1,u_2)-(\bar{x}_1,\bar{x}_2)\|_l+\|(u_1,u_2)-(\bar{x}_1,\bar{x}_2)\|_M}$$

$$= \limsup_{\substack{(u_1,u_2)\to(\bar{x}_1,\bar{x}_2)\\ \|(u_1,u_2)\|_l\le 1}} \frac{(x_1^*-y_1^*)(u_1-\bar{x}_1)+(x_2^*-y_2^*)(u_2-\bar{x}_2)}{\|(u_1,u_2)-(\bar{x}_1,\bar{x}_2)\|_l+\|(u_1,u_2)-(\bar{x}_1,\bar{x}_2)\|_M}. \tag{4.44}$$

We prove

$$y_1^* = y_2^* > 0 \implies \theta \notin \widehat{D}^* P_{\mathbb{B}_l}(\bar{x},\bar{x})(y^*). \tag{4.45}$$

Let $(x_1^*, x_2^*) = \theta$ in (4.43) and (4.44), under the condition $y_1^* = y_2^*$, we have

$$\limsup_{\substack{(u,v)\to((\bar{x}_1,\bar{x}_2),(\bar{x}_1,\bar{x}_2))\\ u\in\mathbb{R}_l^2 \text{ and } v\in P_{\mathbb{B}_l}(u)}} \frac{\langle x^*,u-\bar{x}\rangle-\langle y^*,v-\bar{x}\rangle}{\|u-\bar{x}\|_l+\|v-\bar{y}\|_M}$$

$$= \max\left\{ \limsup_{\substack{((u_1,u_2),(v_1,v_2))\to((\bar{x}_1,\bar{x}_2),(\bar{x}_1,\bar{x}_2))\\ \|(u_1,u_2)\|_l>1,\ u_1\le v_1\le 1-u_2}} \frac{x_1^*(u_1-\bar{x}_1)+x_2^*(u_2-\bar{x}_2)-(y_1^*-y_2^*)(v_1-\bar{x}_1)}{\|(u_1,u_2)-(\bar{x}_1,\bar{x}_2)\|_l+\big\|\big(v_1-\bar{x}_1,-(v_1-\bar{x}_1)\big)\big\|_M}, \limsup_{\substack{(u_1,u_2)\to(\bar{x}_1,\bar{x}_2)\\ \|(u_1,u_2)\|_l\le 1}} \frac{(x_1^*-y_1^*)(u_1-\bar{x}_1)+(x_2^*-y_2^*)(u_2-\bar{x}_2)}{\|(u_1,u_2)-(\bar{x}_1,\bar{x}_2)\|_l+\|(u_1,u_2)-(\bar{x}_1,\bar{x}_2)\|_M} \right\}$$

$$= \max\left\{ \limsup_{\substack{((u_1,u_2),(v_1,v_2))\to((\bar{x}_1,\bar{x}_2),(\bar{x}_1,\bar{x}_2))\\ \|(u_1,u_2)\|_l>1,\ u_1\le v_1\le 1-u_2}} \frac{0}{\|(u_1,u_2)-(\bar{x}_1,\bar{x}_2)\|_l+\big\|\big(v_1-\bar{x}_1,-(v_1-\bar{x}_1)\big)\big\|_M}, \limsup_{\substack{(u_1,u_2)\to(\bar{x}_1,\bar{x}_2)\\ \|(u_1,u_2)\|_l\le 1}} \frac{-y_1^*(u_1-\bar{x}_1)-y_2^*(u_2-\bar{x}_2)}{\|(u_1,u_2)-(\bar{x}_1,\bar{x}_2)\|_l+\|(u_1,u_2)-(\bar{x}_1,\bar{x}_2)\|_M} \right\}$$

$$= \max\left\{ 0, \limsup_{\substack{(u_1,u_2)\to(\bar{x}_1,\bar{x}_2)\\ \|(u_1,u_2)\|_l\le 1}} \frac{-y_1^*\big((u_1-\bar{x}_1)+(u_2-\bar{x}_2)\big)}{\|(u_1,u_2)-(\bar{x}_1,\bar{x}_2)\|_l+\|(u_1,u_2)-(\bar{x}_1,\bar{x}_2)\|_M} \right\}$$

$$= \max\left\{ 0, \limsup_{\substack{(u_1,u_2)\to(\bar{x}_1,\bar{x}_2)\\ \|(u_1,u_2)\|_l\le 1}} \frac{-y_1^*\big((u_1+u_2)-1\big)}{\|(u_1,u_2)-(\bar{x}_1,\bar{x}_2)\|_l+\|(u_1,u_2)-(\bar{x}_1,\bar{x}_2)\|_M} \right\}. \tag{4.46}$$

By the conditions that $0 \le \bar{x}_1 \le 1, 0 \le \bar{x}_2 \le 1$ and $\bar{x}_1 + \bar{x}_2 = 1$, without loose of the generality,

we suppose that $\bar{x}_1 > 0$. For $t > 0$, let $u(t) = (\bar{x}_1 - t, \bar{x}_2)$. If $t$ is small enough, then $\|u\|_l < 1$. Notice that $y_1^* = y_2^* > 0$. Substituting it to (4.46), we have

$$\max\left\{0, \limsup_{\substack{(u_1,u_2)\to(\bar{x}_1,\bar{x}_2) \\ \|(u_1,u_2)\|_l\le 1}} \frac{-y_1^*((u_1+u_2)-1)}{\|(u_1,u_2)-(\bar{x}_1,\bar{x}_2)\|_l + \|(u_1,u_2)-(\bar{x}_1,\bar{x}_2)\|_M}\right\}$$

$$\ge \max\left\{0, \limsup_{\substack{(u_1,u_2)=(\bar{x}_1-t,\bar{x}_2) \\ \|u\|_l\le 1, t\downarrow 0}} \frac{-y_1^*(-t)}{\|(u_1,u_2)-(\bar{x}_1,\bar{x}_2)\|_l + \|(u_1,u_2)-(\bar{x}_1,\bar{x}_2)\|_M}\right\}$$

$$\ge \max\left\{0, \limsup_{t\downarrow 0} \frac{-y_1^*(-t)}{\|(-t,0)\|_l + \|(-t,0)\|_M}\right\}$$

$$= \frac{y_1^*}{2} > 0.$$

This proves (4.45).

Proof (II) of (a) in (v). Next, we prove

$$y_1^* = y_2^* < 0 \implies \theta \in \widehat{D}^* P_{\mathbb{B}_l}(\bar{x}, \bar{x})(y^*). \tag{4.47}$$

Substituting it to (4.46), we have

$$\limsup_{\substack{(u,v)\to((\bar{x}_1,\bar{x}_2),(\bar{x}_1,\bar{x}_2)) \\ u\in\mathbb{R}_l^2 \text{ and } v\in P_{\mathbb{B}_l}(u)}} \frac{\langle x^*, u-\bar{x}\rangle - \langle y^*, v-\bar{x}\rangle}{\|u-\bar{x}\|_l + \|v-\bar{y}\|_M}$$

$$= \max\left\{0, \limsup_{\substack{(u_1,u_2)\to(\bar{x}_1,\bar{x}_2) \\ \|(u_1,u_2)\|_l\le 1}} \frac{-y_1^*((u_1-\bar{x}_1)+(u_2-\bar{x}_2))}{\|(u_1,u_2)-(\bar{x}_1,\bar{x}_2)\|_l + \|(u_1,u_2)-(\bar{x}_1,\bar{x}_2)\|_M}\right\}$$

$$= \max\left\{0, \limsup_{\substack{(u_1,u_2)\to(\bar{x}_1,\bar{x}_2) \\ \|(u_1,u_2)\|_l\le 1}} \frac{-y_1^*((u_1+u_2)-(\bar{x}_1+\bar{x}_2))}{\|(u_1,u_2)-(\bar{x}_1,\bar{x}_2)\|_l + \|(u_1,u_2)-(\bar{x}_1,\bar{x}_2)\|_M}\right\}$$

$$= \max\left\{0, \limsup_{\substack{(u_1,u_2)\to(\bar{x}_1,\bar{x}_2) \\ u_1\ge 0, u_2\ge 0, u_1+u_2\le 1}} \frac{-y_1^*((u_1+u_2)-1)}{\|(u_1,u_2)-(\bar{x}_1,\bar{x}_2)\|_l + \|(u_1,u_2)-(\bar{x}_1,\bar{x}_2)\|_M}\right\}$$

$$\le \max\left\{0, \limsup_{\substack{(u_1,u_2)\to(\bar{x}_1,\bar{x}_2) \\ u_1\ge 0, u_2\ge 0, u_1+u_2\le 1}} \frac{0}{\|(u_1,u_2)-(\bar{x}_1,\bar{x}_2)\|_l + \|(u_1,u_2)-(\bar{x}_1,\bar{x}_2)\|_M}\right\} = 0.$$

This proves (4.47).

Proof of (b) in (v). Suppose $y_1^* - y_2^* \neq 0$. Then, we prove

$$y_1^* - y_2^* \neq 0 \implies \theta \notin \widehat{D}^* P_{\mathbb{B}_l}(\bar{x}, \bar{x})(y^*). \tag{4.48}$$

Case 1. $y_1^* - y_2^* > 0$. Let $x^* = (x_1^*, x_2^*) = \theta$ in (4.43) and (4.44), we have

$$\limsup_{\substack{(u,v)\to((\bar{x}_1,\bar{x}_2),(\bar{x}_1,\bar{x}_2)) \\ u\in\mathbb{R}_l^2 \text{ and } v\in P_{\mathbb{B}_l}(u)}} \frac{\langle x^*, u-\bar{x}\rangle - \langle y^*, v-\bar{x}\rangle}{\|u-\bar{x}\|_l + \|v-\bar{y}\|_M}$$

$$= \max\left\{\begin{array}{c} \limsup\limits_{\substack{((u_1,u_2),(v_1,v_2))\to((\bar{x}_1,\bar{x}_2),(\bar{x}_1,\bar{x}_2)) \\ \|(u_1,u_2)\|_l>1,\ u_1\le v_1\le 1-u_2}} \frac{-(y_1^*-y_2^*)(v_1-\bar{x}_1)}{\|(u_1,u_2)-(\bar{x}_1,\bar{x}_2)\|_l + \|(v_1-\bar{x}_1,-(v_1-\bar{x}_1))\|_M}, \\ \limsup\limits_{\substack{(u_1,u_2)\to(\bar{x}_1,\bar{x}_2) \\ \|(u_1,u_2)\|_l\le 1}} \frac{-y_1^*(u_1-\bar{x}_1)-y_2^*(u_2-\bar{x}_2))}{\|(u_1,u_2)-(\bar{x}_1,\bar{x}_2)\|_l + \|(u_1,u_2)-(\bar{x}_1,\bar{x}_2)\|_M} \end{array}\right\}.$$

For $t > 0$, let $u(t) = (\bar{x}_1 - t, \bar{x}_2 + 2t)$ with $\|u(t)\|_l = 1 + t > 1$ and let $v(t) = (\bar{x}_1 - t, 1 - \bar{x}_1 + t)$ with respect to small enough $t$. Substituting these to the first limit in the above equations, we have

$$\geq \max\left\{\limsup_{t\downarrow 0} \frac{-(y_1^*-y_2^*)(-t)}{\|(-t,2t)\|_l + \|(-t,-(-t))\|_M}, \limsup_{\substack{(u_1,u_2)\to(\bar{x}_1,\bar{x}_2) \\ \|(u_1,u_2)\|_l\le 1}} \frac{-y_1^*(u_1-\bar{x}_1)-y_2^*(u_2-\bar{x}_2))}{\|(u_1,u_2)-(\bar{x}_1,\bar{x}_2)\|_l + \|(u_1,u_2)-(\bar{x}_1,\bar{x}_2)\|_M}\right\}$$

$$= \max\left\{\frac{y_1^*-y_2^*}{4}, \limsup_{\substack{(u_1,u_2)\to(\bar{x}_1,\bar{x}_2) \\ \|(u_1,u_2)\|_l\le 1}} \frac{-y_1^*(u_1-\bar{x}_1)-y_2^*(u_2-\bar{x}_2))}{\|(u_1,u_2)-(\bar{x}_1,\bar{x}_2)\|_l + \|(u_1,u_2)-(\bar{x}_1,\bar{x}_2)\|_M}\right\}$$

$$\geq \frac{y_1^*-y_2^*}{4} > 0.$$

This proves that

$$y_1^* - y_2^* > 0 \implies \theta \notin \widehat{D}^* P_{\mathbb{B}_l}(\bar{x}, \bar{x})(y^*). \tag{4.49}$$

Case 2. $y_1^* - y_2^* < 0$. Let $x^* = (x_1^*, x_2^*) = \theta$ in (4.43) and (4.44), Similarly to the proof of (4.49),

$$\limsup_{\substack{(u,v)\to((\bar{x}_1,\bar{x}_2),(\bar{x}_1,\bar{x}_2)) \\ u\in\mathbb{R}_l^2 \text{ and } v\in P_{\mathbb{B}_l}(u)}} \frac{\langle x^*, u-\bar{x}\rangle - \langle y^*, v-\bar{x}\rangle}{\|u-\bar{x}\|_l + \|v-\bar{y}\|_M}$$

$$= \max\left\{ \limsup_{\substack{((u_1,u_2),(v_1,v_2))\to(\bar{x},\bar{x}) \\ \|(u_1,u_2)\|_l>1,\ u_1\le v_1\le 1-u_2}} \frac{x_1^*(u_1-\bar{x}_1) + x_2^*(u_2-\bar{x}_2)-(y_1^*-y_2^*)(v_1-\bar{x}_1)}{\|(u_1,u_2)-(\bar{x}_1,\bar{x}_2)\|_l + \|(v_1-\bar{x}_1,-(v_1-\bar{x}_1))\|_M},\right.$$

$$\left.\limsup_{\substack{(u_1,u_2)\to(\bar{x}_1,\bar{x}_2)\\ \|(u_1,u_2)\|_l\le 1}} \frac{(x_1^*-y_1^*)(u_1-\bar{x}_1)+(x_2^*-y_2^*)(u_2-\bar{x}_2)}{\|(u_1,u_2)-(\bar{x}_1,\bar{x}_2)\|_l+\|(u_1,u_2)-(\bar{x}_1,\bar{x}_2)\|_M}\right\}$$

$$= \max\left\{\limsup_{\substack{((u_1,u_2),(v_1,v_2))\to((\bar{x}_1,\bar{x}_2),(\bar{x}_1,\bar{x}_2))\\ \|(u_1,u_2)\|_l>1,\ u_1\le v_1\le 1-u_2}} \frac{-(y_1^*-y_2^*)(v_1-\bar{x}_1)}{\|(u_1,u_2)-(\bar{x}_1,\bar{x}_2)\|_l+\|(v_1-\bar{x}_1,-(v_1-\bar{x}_1))\|_M},\right.$$

$$\left.\limsup_{\substack{(u_1,u_2)\to(\bar{x}_1,\bar{x}_2)\\ \|(u_1,u_2)\|_l\le 1}} \frac{-y_1^*(u_1-\bar{x}_1)-y_2^*(u_2-\bar{x}_2))}{\|(u_1,u_2)-(\bar{x}_1,\bar{x}_2)\|_l+\|(u_1,u_2)-(\bar{x}_1,\bar{x}_2)\|_M}\right\}.$$

For $t > 0$, let $u(t) = (\bar{x}_1 + t, \bar{x}_2)$ and let $v(t) = (\bar{x}_1 + t, 1 - \bar{x}_1 - t)$ with respect to small enough $t$. Substituting these to the first limit in the above equations, we have

$$\ge \max\left\{\limsup_{t\downarrow 0} \frac{-(y_1^*-y_2^*)t}{\|(t,0)\|_l+\|(t,t)\|_M},\ \limsup_{\substack{(u_1,u_2)\to(\bar{x}_1,\bar{x}_2)\\ \|(u_1,u_2)\|_l\le 1}} \frac{-y_1^*(u_1-\bar{x}_1)-y_2^*(u_2-\bar{x}_2))}{\|(u_1,u_2)-(\bar{x}_1,\bar{x}_2)\|_l+\|(u_1,u_2)-(\bar{x}_1,\bar{x}_2)\|_M}\right\}$$

$$= \max\left\{\frac{-(y_1^*-y_2^*)}{2},\ \limsup_{\substack{(u_1,u_2)\to(\bar{x}_1,\bar{x}_2)\\ \|(u_1,u_2)\|_l\le 1}} \frac{-y_1^*(u_1-\bar{x}_1)-y_2^*(u_2-\bar{x}_2))}{\|(u_1,u_2)-(\bar{x}_1,\bar{x}_2)\|_l+\|(u_1,u_2)-(\bar{x}_1,\bar{x}_2)\|_M}\right\}$$

$$\ge \frac{-(y_1^*-y_2^*)}{2} > 0.$$

This proves that

$$y_1^* - y_2^* < 0 \Longrightarrow \ \theta \notin \widehat{D}^* P_{\mathbb{B}_l}\big((\bar{x}_1,\bar{x}_2),(\bar{x}_1,\bar{x}_2)\big)(y^*). \tag{4.50}$$

(4.48) is proved by (4.49) and (4.50).

Proof of (vi). This can be considered as a special case of part (E) in Theorem 3.2 and the proof of part (vi) is omitted here. □

### 4.3. The covering constant for the metric projection operator on $\mathbb{R}_l^2$.

As an application of the results of Theorem 4.2, in this subsection, we will calculate the covering constant for the set-valued metric projection operator $P_{\mathbb{B}_l}$ on $\mathbb{R}_l^2$, by using its Mordukhovich derivatives proved in Theorem 4.2. For this purpose, we give the following corollary of Theorem 4.2, which will be used in the proof of Theorem 4.4.

**Corollary 4.3.** *Let* $\bar{x} = (\bar{x}_1, \bar{x}_2) \in \mathbb{R}_l^2$ *with* $\bar{x}_1 > 0, \bar{x}_2 > 0$ *and* $\bar{x}_1 + \bar{x}_2 > 1$. *Let* $\bar{y} = (\bar{y}_1, \bar{y}_2) \in P_{\mathbb{B}_l}(\bar{x})$. *There is* $y^* = (y_1^*, y_2^*) \in \mathbb{B}_M$ *such that* $\|y^*\|_M = 1$ *and*

$$\theta \in \widehat{D}^* P_{\mathbb{B}_l}(\bar{x}, \bar{y})(y^*).$$

*Proof*. In Theorem 4.2, we can take $(y_1^*, y_2^*) = (1,1)$ or $(-1,-1)$. □

In the next theorem, we investigate the solutions of the covering constant for the metric projection operator $P_{\mathbb{B}_l}$ on $\mathbb{R}^2_l$. Similarly to Theorem 4.2, we concentrate our study in Quadrant I and $\mathbb{B}^o_l$.

**Theorem 4.4**. *Let $\bar{x} = (\bar{x}_1, \bar{x}_2) \in \mathbb{R}^2_l$ and let $\bar{y} = (\bar{y}_1, \bar{y}_2) \in P_{\mathbb{B}_l}(\bar{x})$. Then*

(i) $\hat{\alpha}(P_{\mathbb{B}_l}, \bar{x}, \bar{y}) = 0$, *for any $\bar{x} \in Q_I \backslash \mathbb{B}^o_l$ (with $\bar{x} \succcurlyeq_2 \theta$ and $\|\bar{x}\|_l \geq 1$);*
(ii) $\hat{\alpha}(P_{\mathbb{B}_l}, \bar{x}, \bar{y}) = \hat{\alpha}(P_{\mathbb{B}_l}, \bar{x}, \bar{x}) = 1$, *for any $\bar{x} \in \mathbb{B}^o_l$.*

*Proof*. Proof of (i). Let $\bar{x} = (\bar{x}_1, \bar{x}_2) \in Q_I$ with $\bar{x}_1 \geq 0, \bar{x}_2 \geq 0$ and $\bar{x}_1 + \bar{x}_2 \geq 1$. Let $\eta > 0$ be arbitrarily given. Let $\mathbb{B}_l(\bar{x}, \eta)$ and $\mathbb{B}_l(\bar{y}, \eta)$ be the closed balls in $\mathbb{R}^2_l$ with radius $\eta$ centered at $\bar{x}$ and $\bar{y}$, respectively. There are $x = (x_1, x_2) \in \mathbb{B}_l(\bar{x}, \eta)$ satisfying that $x_1 > 0, x_2 > 0$ and $x_1 + x_2 > 1$, and $y = (y_1, y_2) \in \mathbb{B}_l(\bar{y}, \eta)$ with $y = (y_1, y_2) \in P_{\mathbb{B}_l}(x)$. By Proposition 4.1 and Corollary 4.3, there is $y^* = (y^*_1, y^*_2) \in \mathbb{B}_M \backslash \{\theta\}$ (for example, we may take $y^*_1 = y^*_2 = \pm 1$) such that

$$\theta \in \widehat{D}^* P_{\mathbb{B}_l}(\bar{x}, \bar{y})(y^*). \tag{4.51}$$

Then, for this arbitrarily given $\eta > 0$, by (4.51), we have

$$\begin{aligned}
&\inf\{\|z^*\|_M : z^* \in \widehat{D}^* P_{\mathbb{B}_l}(x, y)(w^*), x \in \mathbb{B}_l(\bar{x}, \eta), y \in P_{\mathbb{B}_l}(x) \in \mathbb{B}_l(\bar{y}, \eta), \|w^*\|_M = 1\} \\
&\leq \inf\{\|\theta\|_M : \theta \in \widehat{D}^* P_{\mathbb{B}_l}(x, y)(y^*), x \in \mathbb{B}_l(\bar{x}, \eta), y \in P_{\mathbb{B}_l}(x) \in \mathbb{B}_l(\bar{y}, \eta), \|y^*\|_M = 1\} \\
&= 0.
\end{aligned} \tag{4.52}$$

Proof of (ii). Let $\bar{x} = (\bar{x}_1, \bar{x}_2) \in \mathbb{B}^o_l$. There is $\delta > 0$ such that

$$\mathbb{B}_l(\bar{x}, \eta) \subseteq \mathbb{B}^o_l, \text{ for } \eta < \delta \quad \text{and} \quad \mathbb{B}_l(\bar{x}, \eta) \backslash \mathbb{B}^o_l \neq \emptyset, \text{ for } \eta \geq \delta.$$

For $\bar{x} = (\bar{x}_1, \bar{x}_2) \in \mathbb{B}^o_l$, $P_{\mathbb{B}_l}(\bar{x}) = \bar{x}$, we have

$$\begin{aligned}
&\sup_{0<\eta<\delta} \inf\{\|z^*\|_M : z^* \in \widehat{D}^* P_{\mathbb{B}_l}(x, y)(w^*), x \in \mathbb{B}_l(\bar{x}, \eta), y \in P_{\mathbb{B}_l}(x) \in \mathbb{B}_l(\bar{y}, \eta), \|w^*\|_M = 1\} \\
&= \sup_{0<\eta<\delta} \inf\{\|z^*\|_M : z^* \in \widehat{D}^* P_{\mathbb{B}_l}(x, y)(w^*), x \in \mathbb{B}_l(\bar{x}, \eta), x = P_{\mathbb{B}_l}(x) \in \mathbb{B}_l(\bar{x}, \eta), \|w^*\|_M = 1\} \\
&= \sup_{0<\eta<\delta} \inf\{\|w^*\|_M : w^* \in \widehat{D}^* P_{\mathbb{B}_l}(x, y)(w^*), x \in \mathbb{B}_l(\bar{x}, \eta), x = P_{\mathbb{B}_l}(x) \in \mathbb{B}_l(\bar{x}, \eta), \|w^*\|_M = 1\} \\
&= 1.
\end{aligned} \tag{4.53}$$

Similarly to the proof of (4.52), we have

$$\begin{aligned}
&\sup_{\eta\geq\delta} \inf\{\|z^*\|_M : z^* \in \widehat{D}^* P_{\mathbb{B}_l}(x, y)(w^*), x \in \mathbb{B}_l(\bar{x}, \eta), y \in P_{\mathbb{B}_l}(x) \in \mathbb{B}_l(\bar{y}, \eta), \|w^*\|_M = 1\} \\
&\leq \sup_{\eta\geq\delta} \inf\{\|\theta\|_M : \theta \in \widehat{D}^* P_{\mathbb{B}_l}(x, y)(y^*), x \in \mathbb{B}_l(\bar{x}, \eta), y \in P_{\mathbb{B}_l}(x) \in \mathbb{B}_l(\bar{x}, \eta), \|y^*\|_M = 1\} \\
&= 0.
\end{aligned} \tag{4.54}$$

Then, by the definition of covering constant, part (ii) is proved by (4.53) and (4.54). □

## Appendix

### A1. Mordukhovich Differentiability of the Metric Projection Operator in the 3-d Banach Space

Let $\mathbb{R}^3$ be the ordinary 3-d real vector space. Similarly to the spaces $\left(\mathbb{R}_l^2, \|\cdot\|_l\right)$ and $(\mathbb{R}_M^2, \|\cdot\|_M)$ studied in section 4, let $\left(\mathbb{R}_l^3, \|\cdot\|_l\right)$ denote the 3-dimensional $l_1$-Banach space, in which the $l_1$-norm $\|\cdot\|_l$ on $\mathbb{R}_l^3$ is defined by

$$\|x\|_l = |x_1| + |x_2| + |x_3|, \text{ for any } x = (x_1, x_2, x_3) \in \mathbb{R}_l^3. \tag{5.1}$$

The topological dual space of $\left(\mathbb{R}_l^3, \|\cdot\|_l\right)$ is the Banach space $(\mathbb{R}_M^3, \|\cdot\|_M)$, in which the maximum-norm $\|\cdot\|_M$ on $\mathbb{R}_M^3$ is defined by

$$\|y\|_M = \max\{|y_1|, |y_2|, |y_3|\}, \text{ for any } y = (y_1, y_2, y_3) \in \mathbb{R}_M^3. \tag{5.2}$$

The real pairing between $\mathbb{R}_M^3$ and $\mathbb{R}_l^3$ is written by $\langle\cdot,\cdot\rangle$ and

$$\langle y, x\rangle = x_1y_1 + x_2y_2 + x_3y_3 \text{ for any } x = (x_1, x_2, x_3) \in \mathbb{R}_l^3 \text{ and } y = (y_1, y_2, y_3) \in \mathbb{R}_M^3.$$

$\mathbb{R}_l^3$ and $\mathbb{R}_M^3$ have the same set of elements with $\mathbb{R}^3$. In particular, their null element is denoted by $\theta$ = (0, 0, 0). Let $\mathbb{B}_l$ and $\mathbb{B}_M$ denote the unit balls in $\mathbb{R}_l^3$ and $\mathbb{R}_M^3$, respectively. Actually, $\mathbb{B}_l$ is the closed cube in the space $\mathbb{R}^3$ with vertexes (1, 0, 0), (0, 1, 0), $(0, 0, 1)$, $(-1, 0, 0)$, $(0, -1, 0)$ and (0, 0, $-1$). The set-valued metric projection operator $P_{\mathbb{B}_l}\colon \mathbb{R}_l^3 \rightrightarrows \mathbb{B}_l$ is defined, for any $x = (x_1, x_2, x_2) \in \mathbb{R}_l^3$, by

$$P_{\mathbb{B}_l}(x) = \left\{y \in \mathbb{B}_l \colon \|x - y\|_l = \min_{v \in \mathbb{B}_l} \|x - v\|_l\right\}. \tag{5.3}$$

In the following proposition, similarly to the Proposition 4.1, we find the explicit solutions of the set-valued metric projection operator $P_{\mathbb{B}_l}$ mainly in Octant I. The explicit solutions of the metric projection operator $P_{\mathbb{B}_l}$ can be similarly solved on other seven octants.

Let $O_I$ denote the Octant I in the space $\mathbb{R}^3$, which is the positive cone in $\mathbb{R}^3_l$. Let $\preccurlyeq_3$ denote the partial order on $\mathbb{R}^3_l$ induced by $O_I$, which is defined, for any $x = (x_1, x_2, x_3)$ and $u = (u_1, u_2, u_3) \in \mathbb{R}^2_l$, by

$$x \preccurlyeq_3 u \quad \Leftrightarrow \quad x_i \leq u_i, \text{ for } i = 1, 2, 3.$$

More strictly, we write $$x \prec_3 u \quad \Leftrightarrow \quad x_i < u_i, \text{ for } i = 1, 2, 3.$$

Let $E$, $F$, $G$ be three nonlinear points in $\mathbb{R}^3$. Let $\Delta\{E, F, G\}$ denote the closed triangle in $\mathbb{R}^3$ with vertexes $E$, $F$ and $G$. In particular, we write

$$\Delta := \Delta\{(1, 0, 0), (0, 1, 0), (0, 0, 1)\}.$$

**Proposition 5.1**. *Let* $x = (x_1, x_2, x_3) \in \mathbb{R}^3_l$. *Suppose* $x_i \geq 0$, *for* $i = 1, 2, 3$ *and* $x_1 + x_2 + x_3 \geq 1$, *then* $P_{\mathbb{B}_l}$ *satisfies the following equations*:

$$P_{\mathbb{B}_l}(x) = \left\{y \in \mathbb{B}_l : \|y\|_l = 1, \|x - y\|_l = \|x\|_l - 1, \begin{matrix} 0 \leq y_n \leq x_n \wedge 1, & \text{if } x_n \geq 0, \\ x_n \vee (-1) \leq y_n \leq 0, & \text{if } x_n < 0, \end{matrix} \; n = 1, 2, 3\right\}. \quad (5.3)$$

*In particular*,

(i) *If* $x_i \geq 1$, *for* $i = 1, 2, 3$, *then*

$$P_{\mathbb{B}_l}(x_1, x_2, x_3) = \Delta.$$

(ii) *If* $x_1 = 0, x_2 > 0, x_3 > 0$ *and* $x_2 + x_3 \geq 1$, *then*

$$P_{\mathbb{B}_l}(0, x_2, x_3) = \overline{(0,\ x_2 \wedge 1, 1 - x_2 \wedge 1), (0, 1 - x_3 \wedge 1, x_3 \wedge 1)}.$$

(iii) *If* $x_1 > 0, x_2 = 0, x_3 > 0$ *and* $x_1 + x_3 \geq 1$, *then*

$$P_{\mathbb{B}_l}(x_1, 0, x_3) = \overline{(x_1 \wedge 1, 0, 1 - x_1 \wedge 1), (1 - x_3 \wedge 1, 0, x_3 \wedge 1)}.$$

(iv) *If* $x_1 > 0, x_2 > 0, x_3 = 0$ *and* $x_1 + x_2 \geq 1$, *then*

$$P_{\mathbb{B}_l}(x_1, x_2, 0) = \overline{(x_1 \wedge 1, 1 - x_1 \wedge 1, 0), (1 - x_2 \wedge 1, x_2 \wedge 1, 0)}.$$

*Proof*. The proof of (5.3) follows from the proof of Theorem 3.1, which can indeed be considered as special cases of Theorem 3.1. So, the proof of (5.3) is omitted here and we only directly prove the special case part (i). The proofs of all other parts (ii−iv) are similar to the proof of Proposition 4.1 and they are omitted here too.

Proof of (i). Let $x = (x_1, x_2, x_3) \in \mathbb{R}^3_l$ with $x_i \geq 1$, for $i = 1, 2, 3$. Then, for any $y = (y_1, y_2, y_3) \in \Delta$,

$$\begin{aligned} &\|x - y\|_l \\ &= |x_1 - y_1| + |x_2 - y_2| + |x_3 - y_3| \\ &= x_1 - y_1 + x_2 - y_2 + x_3 - y_3 \end{aligned}$$

$$= x_1 + x_2 + x_3 - (y_1 + y_2 + y_3)$$

$$= x_1 + x_2 + x_3 - 1$$

$$= \|x\|_l - 1. \qquad (5.4)$$

Every $u = (u_1, u_2, u_3) \in \mathbb{B}_l \backslash \Delta$ satisfies that $\|u\|_l \le 1$ and $u_i < 1$, for $i$ = 1, 2, 3. It also satisfies that $u_1 + u_2 + u_3 < 1$. This implies that

$$\|x - u\|_l$$

$$= |x_1 - u_1| + |x_2 - u_2| + |x_3 - u_3|$$

$$= x_1 - u_1 + x_2 - u_2 + x_3 - u_3$$

$$= x_1 + x_2 + x_3 - (u_1 + u_2 + u_3)$$

$$> \|x\|_l - 1, \text{ for any } u = (u_1, u_2, u_3) \in \mathbb{B}_l \backslash \Delta. \qquad (5.5)$$

(i) is proved by (5.4) and (5.5). □

From the results of Theorem 4.2, we may anticipate that the Mordukhovich derivatives of the metric projection operator $P_{\mathbb{B}_l}$ in $\mathbb{R}^3_M$ must be very complicated. Considering the possibility of the paper to exceeds the limit of length, we only find the Mordukhovich derivatives of $P_{\mathbb{B}_l}$ in some special cases in Octant I in the space $\mathbb{R}^3_l$. All other cases can be similarly considered.

**Theorem 5.2**. *Let* $\bar{x} = (\bar{x}_1, \bar{x}_2, \bar{x}_3) \in \mathbb{R}^3_l$ *with* $\|\bar{x}\|_l \ge 1$. *Let* $\bar{y} = (\bar{y}_1, \bar{y}_2, \bar{y}_3) \in P_{\mathbb{B}_l}(\bar{x})$. *Then, we have*

(A) *If* $\bar{x} \in O_I$, *then, for any* $y^* = (y_1^*, y_2^*, y_3^*) \in \mathbb{R}^3_M$,

$$\widehat{D}^* P_{\mathbb{B}_l}(\bar{x}, \bar{y})(y^*) \subseteq -O_I.$$

(B) *If* $\bar{x} \in -Q_I$, *then, for any* $y^* = (y_1^*, y_2^*, y_3^*) \in \mathbb{R}^3_M$,

$$\widehat{D}^* P_{\mathbb{B}_l}(\bar{x}, \bar{y})(y^*) \subseteq O_I.$$

(C) *Suppose* $\bar{x} \in O_I^o \cup (-O_I)^o$. *Let* $y^* = (y_1^*, y_2^*, y_3^*) \in \mathbb{R}^3_M$. *If* $y_1^* = y_2^* = y_3^*$, *then*

$$\theta \in \widehat{D}^* P_{\mathbb{B}_l}(\bar{x}, \bar{y})(y^*).$$

*In particular, if* $\bar{x} \succ_3 (1, 1, 1)$ *and* $\bar{y} \in P_{\mathbb{B}_l}(\bar{x})$ *satisfying* $\theta \prec_3 \bar{y} \prec_3 (1, 1, 1)$, *then*

(I) *the equations* $y_1^* = y_2^* = y_3^*$ *do not hold* $\Longrightarrow$ $\widehat{D}^* P_{\mathbb{B}_l}(\bar{x}, \bar{y})(y^*) = \emptyset$.

(II) $y_1^* = y_2^* = y_3^* \Longrightarrow \widehat{D}^* P_{\mathbb{B}_l}(\bar{x}, \bar{y})(y^*) = \{\theta\}$.

*Proof*. When we consider a special case that $\bar{x} = (\bar{x}_1, \bar{x}_2, \bar{x}_3) \in \mathbb{R}^3_l$ with $\bar{x} \succ_3 \theta$ and $\|\bar{x}\|_l > 1$, there is $\alpha > 0$ such that, for any $u = (u_1, u_2, u_3) \in \mathbb{R}^3_l$,

$$\|u - \bar{x}\|_l < \alpha \quad \Longrightarrow \quad u \succ_3 \theta \text{ and } \|u\|_l > 1.$$

By (i) in Proposition 5.1, this implies that

$$\|u-\bar{x}\|_l < \alpha \quad \Longrightarrow \quad P_{\mathbb{B}_l}(u) = \Delta.$$

Hence, in the case that $\|u-\bar{x}\|_l < \alpha$, for any $v = (v_1, v_2, v_3) \in P_{\mathbb{B}_l}(u)$, by Proposition 5.1, $v$ satisfies that $\theta \preccurlyeq_3 v \preccurlyeq_3 u \wedge 1$ and $\|v\|_l = 1$. Then, for any $m \in \{1, 2, 3\}$, $v_m = 1 - \sum_{n\neq m} v_n$. This implies that $v$ can be rewritten as

$$v = (v_1, v_2, 1-(v_1+v_2)) = (v_1, 1-(v_1+v_3), v_3) = (1-(v_2+v_3), v_2, v_3).$$

Let $y^* = (y_1^*, y_2^*, y_3^*) \in \mathbb{R}_M^3$ and $x^* = (x_1^*, x_2^*, x_3^*) \in \mathbb{R}_M^3$. We calculate

$$\limsup_{\substack{(u,v)\to(\bar{x},\bar{y})\\ v=(v_1,v_2)\in P_{\mathbb{B}_l}(u)}} \frac{\langle x^*, u-\bar{x}\rangle - \langle y^*, v-\bar{y}\rangle}{\|u-\bar{x}\|_l + \|v-\bar{y}\|_M}$$

$$= \limsup_{\substack{(u,v)\to(\bar{x},\bar{y}), \|u-\bar{x}\|_l<\alpha\\ v=(v_1,v_2,v_3)\in P_{\mathbb{B}_l}(u)}} \frac{\langle x^*, u-\bar{x}\rangle - \langle y^*, v-\bar{y}\rangle}{\|u-\bar{x}\|_l + \|v-\bar{y}\|_M}$$

$$= \limsup_{\substack{(u,v)\to(\bar{x},\bar{y}), \|u-\bar{x}\|_l<\alpha\\ v\in P_{\mathbb{B}_l}(u)}} \frac{\sum_{n=1}^{3} x_n^*(u_n-\bar{x}_n) - \sum_{n\neq m}(y_n^*-y_m^*)(v_n-\bar{y}_n)}{\|u-\bar{x}\| + \|v-\bar{y}\|_M}, \text{ for } m = 1, 2, 3$$

$$= \limsup_{\substack{(u,v)\to(\bar{x},\bar{y}), \|u-\bar{x}\|_l<\alpha\\ v\in P_{\mathbb{B}_l}(u)}} \frac{\sum_{n=1}^{3} x_n^*(u_n-\bar{x}_n) + \sum_{n\neq m}(y_m^*-y_n^*)(v_n-\bar{y}_n)}{\|u-\bar{x}\| + \|v-\bar{y}\|_M}, \text{ for } m = 1, 2, 3. \tag{5.6}$$

By (5.6), the proofs of parts (A), (B) and (C) of this theorem are similar with the proofs of (A), (B) and (C) in Theorem 4.2, respectively. They are omitted here.

In particular, suppose $\bar{x} \succ_3 (1, 1, 1)$ and $\bar{y} \in P_{\mathbb{B}_l}(\bar{x})$ with $\theta \prec_3 \bar{y} \prec_3 (1, 1, 1)$. Then, the proofs of (I) and (II) of this theorem are similar to the proofs of (I) and (II) of (a) in (i) of Theorem 4.2, respectively. They are omitted here. □

**Theorem 5.3**. *The covering constant for $P_{\mathbb{B}_l}$ has the following properties*.

(i) *Let $\bar{x} \in O_l \backslash \mathbb{B}_l^o$ ($\bar{x} \succcurlyeq_3 \theta$ and $\|\bar{x}\|_l \geq 1$). Then, for any $\bar{y} \in P_{\mathbb{B}_l}(\bar{x})$*

$$\hat{\alpha}(P_{\mathbb{B}_l}, \bar{x}, \bar{y}) = 0.$$

(ii) *Let $\bar{x} \in \mathbb{B}_l^o$. Then, $P_{\mathbb{B}_l}(\bar{x}) = \bar{x}$ and*

$$\hat{\alpha}(P_{\mathbb{B}_l}, \bar{x}, \bar{x}) = 1.$$

*Proof*. By Theorem 5.2, the proof of this theorem is similar to the proof of Theorem 3.4 and it is omitted here. □

**A2. A Directly Proof of Part (I) of Proposition 4.1**

To directly prove (I) of Proposition 4.1, we consider the following part of Quadrant I

$$\text{(I)} = \{(x_1, x_2) \in \mathbb{R}_l^2 : x_1 \geq 0, x_2 \geq 0 \text{ and } x_1 + x_2 \geq 1\}$$

to the following four parts to precisely prove (a)–(d) in part (I), respectively.

(Part 1). Let $Q = \{x = (x_1, x_2) \in \mathbb{R}_l^2 : x_1 \geq 1 \text{ and } x_2 \geq 1\}$, which is the closed, convex and pointed right-angle cone in the plane $\mathbb{R}^2$ with vertex (1, 1). Then, we have

$$P_{\mathbb{B}_l}(x_1, x_2) = \overline{A_l B_l}, \text{ for any } (x_1, x_2) \in Q. \tag{6.1}$$

Here, $\overline{A_l B_l}$ is the closed segment in the plane $\mathbb{R}^2$ with ending points $A_l$ and $B_l$. Then, we prove (6.1). Let $(x_1, x_2) \in Q$. For any $(y_1, y_2) \in \overline{A_l B_l}$, it satisfies that

$$y_2 + y_1 - 1 = 0,\ 0 \leq y_1 \leq 1 \text{ and } 0 \leq y_2 \leq 1.$$

Then, we calculate

$$\|(x_1, x_2) - (y_1, y_2)\|_l = |x_1 - y_1| + |x_2 - y_2|$$

$$= x_1 - y_1 + x_2 - y_2 = x_1 + x_2 - (y_2 + y_1) = x_1 + x_2 - 1. \tag{6.2}$$

For any $(u_1, u_2) \in \mathbb{B}_l \backslash \overline{A_l B_l}$, it satisfies

$$u_2 + u_1 < 1. \tag{6.3}$$

Then, by (6.3), we calculate

$$\|(x_1, x_2) - (u_1, u_2)\|_l = |x_1 - u_1| + |x_2 - u_2|$$

$$= x_1 - u_1 + x_2 - u_2 = x_1 + x_2 - (u_2 + u_1) > x_1 + x_2 - 1. \tag{6.4}$$

(6.1) is proved by (6.2) and (6.4).

(Part 2). Let $R = \{(x_1, x_2) \in \mathbb{R}_l^2 : 0 \leq x_1 \leq 1, 0 \leq x_2 \leq 1 \text{ and } x_2 + x_1 \geq 1\}$, which is an isosceles right triangle in the plane $\mathbb{R}^2$. Then, we have

$$P_{\mathbb{B}_l}(x_1, x_2) = \overline{(x_1, 1 - x_1), (1 - x_2, x_2)}, \text{ for any } (x_1, x_2) \in R. \tag{6.5}$$

Here, $\overline{(x_1, 1 - x_1), (1 - x_2, x_2)}$ is the closed segment in the plane $\mathbb{R}^2$ with ending points $(x_1, 1 - x_1)$ and $(1 - x_2, x_2)$. Next, we prove (6.5).

Let $(x_1, x_2) \in R$. Then, for any $(y_1, y_2) \in \overline{(x_1, 1 - x_1), (1 - x_2, x_2)}$, it satisfies that

$$y_2 + y_1 - 1 = 0,\ \ 1 - x_2 \leq y_1 \leq x_1 \quad \text{and} \quad 1 - x_1 \leq y_2 \leq x_2. \tag{6.6}$$

Then, similarly to the proof of (6.2), we calculate

$$\|(x_1, x_2) - (y_1, y_2)\|_l = |x_1 - y_1| + |x_2 - y_2|$$

$$= x_1 - y_1 + x_2 - y_2 = x_1 + x_2 - (y_2 + y_1) = x_1 + x_2 - 1. \tag{6.7}$$

For any $(u_1, u_2) \in \mathbb{B}_l \backslash \overline{(x_1, 1 - x_1), (1 - x_2, x_2)}$, there are three cases for consideration.

Case 1. Suppose that $(u_1, u_2) \in \mathbb{B}_l \backslash \overline{A_l B_l}$. That is, $u_2 + u_1 < 1$. Then, by (6.6), we calculate

$$\|(x_1, x_2) - (u_1, u_2)\|_l = |x_1 - u_1| + |x_2 - u_2|$$

$$= x_1 - u_1 + x_2 - u_2 = x_1 + x_2 - (u_2 + u_1) > x_1 + x_2 - 1. \tag{6.8}$$

The following two cases are for $(u_1, u_2) \in \overline{A_l B_l} \backslash \overline{(x_1, 1 - x_1), (1 - x_2, x_2)}$.

Case 2. Suppose that

$$u_2 + u_1 = 1,\ \ x_1 < u_1 \leq 1 \text{ and } \ 0 \leq u_2 < 1 - x_1 < x_2. \tag{6.9}$$

By (6.9), we calculate

$$\|(x_1, x_2) - (u_1, u_2)\|_l = |x_1 - u_1| + |x_2 - u_2|$$

$$= u_1 - x_1 + x_2 - u_2 > u_1 - x_1 + x_2 - 1 + x_1$$

$$> x_1 - x_1 + x_2 - 1 + x_1 = x_1 + x_2 - 1. \tag{6.10}$$

Case 3. Suppose that $u_2 + u_1 = 1,\ \ x_2 < u_2 \leq 1$ and $\ 0 \leq u_1 < 1 - x_2 < x_1$. By (6.9), we calculate

$$\|(x_1, x_2) - (u_1, u_2)\|_l = |x_1 - u_1| + |x_2 - u_2|$$

$$= x_1 - u_1 + u_2 - x_2 > x_1 - 1 + x_2 + \ u_2 - x_2$$

$$> x_1 - 1 + x_2 + \ x_2 - x_2 = x_1 + x_2 - 1. \tag{6.11}$$

(6.7) is proved by (6.8), (6.10) and (6.11).

(Part 3). Let $M = \left\{(x_1, x_2) \in \mathbb{R}_l^2 : x_1 \geq 1, 0 \leq x_2 \leq 1\right\}$. For any $(x_1, x_2) \in M$, we have

$$\|(x_1, x_2) - (1, 0)\|_l = |x_1 - 1| + |x_2| = x_1 + x_2 - 1.$$

Similarly to the proof of (6.5), one can prove that,

$$P_{\mathbb{B}_l}(x_1, x_2) = \overline{(1, 0),\ (1 - x_2, x_2)}, \text{ for any } (x_1, x_2) \in M.$$

(Part 4). Let $L = \left\{(x_1, x_2) \in \mathbb{R}_l^2 : 0 \leq x_1 \leq 1,\ x_2 \geq 1\right\}$. For any $(x_1, x_2) \in L$, we have

$$\|(x_1, x_2) - (0, 1)\|_l = |x_1 - 0| + |x_2 - 1| = x_1 + x_2 - 1.$$

Similarly to the proof of (6.5), one can prove that,

$$P_{\mathbb{B}_l}(x_1, x_2) = \overline{(x_1, 1 - x_1), (0,1)}, \text{ for any } (x_1, x_2) \in L.$$

This completes the directly proof of (I) of Proposition 4.1. □